\documentclass[11pt]{amsart}

\title[Affine Morphisms]{Affine modules and the Drinfeld center}

\author[ \large
P\lowercase{aramita} D\lowercase{as},
S\lowercase{hamindra} G\lowercase{hosh and}
V\lowercase{ed} G\lowercase{upta}]
{\bf \large P\lowercase{aramita} D\lowercase{as},
S\lowercase{hamindra} K\lowercase{umar} G\lowercase{hosh and}
V\lowercase{ed} P\lowercase{rakash} G\lowercase{upta}}

\thanks{Das and Ghosh were supported by K.~U.~Leuven BOF Research Grant OT/08/032, and Gupta was supported by ERC Starting Grant VNALG-200749 and the Institute of Mathematical Sciences (Chennai)}

\date{}

\address{Stat-Math Unit, Indian Statistical Institute, Kolkata, INDIA}
\email{paramita.das@isical.ac.in, shami@isical.ac.in}
\address{Department of Mathematics, Shiv Nadar University, U.P., INDIA}
\email{ved.gupta@snu.edu.in}

\usepackage{amsmath}
\usepackage{amsfonts}
\usepackage{latexsym}
\usepackage{amssymb}
\usepackage{mathrsfs}
\usepackage{amscd}
\usepackage{hyperref}
\usepackage{psfrag}
\usepackage{graphicx}
\usepackage{fullpage}

\newtheorem{thm}{Theorem}[section] 
\newtheorem{lem}[thm]{Lemma}
\newtheorem{cor}[thm]{Corollary} 
\newtheorem{prop}[thm]{Proposition}
\newtheorem{defn}[thm]{Definition}
\newtheorem{rem}[thm]{Remark}

 \newenvironment{pf}{\noindent{\em Proof:}}{\hfill
   $\Box$\vspace*{1mm}}

\newcommand{\comments}[1]{}
\newcommand{\ra}{\rightarrow}
\newcommand{\lra}{\longrightarrow}

\newcommand{\mbb}{\mathbb}
\newcommand{\mcal}{\mathcal}

\newcommand{\N}{\mathbb N}

\newcommand{\C}{\mathbb{C}}

\newcommand{\mscr}{\mathscr}
\newcommand{\vlon}{\varepsilon}

\newcommand{\oline}{\overline}

\newcommand{\vphi}{\varphi}


\begin{document}
\global\long\def\vlon{\varepsilon}
\global\long\def\bt{\bowtie}
\global\long\def\ul#1{\underline{#1}}
\global\long\def\ol#1{\overline{#1}}
\global\long\def\norm#1{\left\|{#1}\right\|}
\global\long\def\os#1#2{\overset{#1}{#2}}
\global\long\def\us#1#2{\underset{#1}{#2}}
\global\long\def\ous#1#2#3{\overset{#1}{\underset{#3}{#2}}}
\global\long\def\t#1{\text{#1}}
\global\long\def\lrsuf#1#2#3{\vphantom{#2}_{#1}^{\vphantom{#3}}#2^{#3}}

\maketitle
\begin{abstract}
Given a finite index subfactor, we show that the {\em affine morphisms at zero level} in the affine category over the planar algebra associated to the subfactor is isomorphic to the fusion algebra of the subfactor as a $*$-algebra.
This identification paves the way to analyze the structure of affine $P$-modules with weight zero for any subfactor planar algebra $P$ (possibly having infinite depth).
Further, for irreducible depth two subfactor planar algebras, we establish an additive equivalence between the category of affine $P$-modules and  the center of the category of $N$-$N$-bimodules generated by $L^2(M)$; this partially verifies a conjecture of Jones and Walker.
\end{abstract}
\section{Introduction}
The standard invariant of a subfactor, which - in certain situations turns out to be a complete invariant - has been described in many seemingly different ways, for instance as a certain category of bimodules (see \cite{Bis97}), as lattices of finite dimensional $C^*$ algebras satisfying certain properties (see \cite{Pop95}), as an algebraic system comprising of graphs, fusion rules and quantum $6j$ symbols (see \cite{Oc2}) or as a planar algebra (see \cite{Jon}). In fact, the theory of planar algebras was initiated by Jones as a tool to study subfactors. The graphical calculus of pictures on a plane turned out to be extremely handy in analyzing the combinatorial data present in a subfactor. Although intimately connected to the theory of subfactors from the outset, planar algebra soon became a subject in its own merit. Moreover, quite recently it has found connections with the theories of random matrices and free probability as well - see \cite{GJS}.

Further, in \cite{Jon01}, Jones introduced the notion of `modules over a planar algebra' or `annular representations', wherein he explicitly obtained all the irreducible modules over the Temperley-Lieb planar algebras for index greater than 4.
Modules over planar algebras have been used in constructing subfactors of index less than 4, namely the subfactors with principal graphs, $E_6$ and $E_8$ - see \cite{Jon01}.
More recently, they have also found application in constructing the Haagerup subfactor - see \cite{Pet}.
Such modules for the group planar algebras were studied by the second-named author in \cite{Gho06} where an equivalence (as additive categories) was established between the category of annular representations over a group planar algebra (that is, planar algebra associated to the fixed point subfactor arising from an outer action of a finite group) and the representation category of a non-trivial quotient of the quantum double of the group over a certain ideal.
The appearance of a non-trivial quotient was due to the fact that the isotopy on annular tangles need not preserve the boundaries of the external and the distinguished internal discs.
On the other hand, {\em affine isotopy} (introduced in \cite{JR06}) does preserve the boundaries of the annulus; in fact, the category of affine modules of a group planar algebra becomes equivalent to the representation category of the quantum double of the group.
Affine modules for the Temperley-Lieb planar algebras were studied in \cite{JR06}.
Certain finiteness results for affine modules of finite depth planar algebras were also established in \cite{Gho}.

The work in this paper was motivated by an attempt to understand the following conjecture of Vaughan Jones and Kevin Walker:

\noindent{\em The category of affine representations of a finite depth subfactor planar algebra is equivalent to  the Drinfeld center of the bimodule category associated to the subfactor.}

This is evident in the case of group planar algebra where the center is equivalent to the representation category of the quantum double of the group; similar results also appear in the world of TQFT's in the work of Kevin Walker.
As a step towards this, we first established an isomorphism between the affine morphisms {\it at zero level} (defined at the beginning of Section \ref{mt}) and the fusion algebra of the bimodule category - as suggested to the second named author by Vaughan Jones and Dietmar Bisch.
Later, this helped us in constructing affine modules with weight zero.
Moreover, we verify the above conjecture in the case of irreducible depth two subfactors.

We now briefly describe the organization of this paper.

Section \ref{prelim} begins with a brief recollection (mainly from \cite{Jon} and \cite{DGG}) of certain basic aspects of planar algebras and their relationship with subfactors and setting up some notations.
For the sake of completeness, in the second part of Section \ref{prelim}, we present a detailed description of the affine category over a planar algebra.

Section \ref{mt} is devoted in proving one of the main theorems in this article, namely Theorem \ref{main-theorem},
the proof of which is divided in three subsections. In the first part, we find a nice spanning set (indexed by the isomorphism class of the irreducible bimodules appearing in the standard invariant of the subfactor) of the space of affine morphisms at zero level.
Here, we crucially use a specific type of affine tangles, namely, the $\Psi^m_{\vlon k , \eta l}$'s and the fact that any affine morphism comes from the action of these affine tangles on regular tangles.
In the second part, we obtain an equivalence relation on planar tangles induced by the effect of affine isotopy.
We use this equivalence relation to show the linear independence of the spanning set, in the last part.

The canonical trace in the fusion algebra induces, via Theorem \ref{main-theorem}, a faithful tracial state on the space of affine morphisms at zero level. In Section \ref{regmod}, we first give a pictorial formulation of this trace.
With this faithful trace at our disposal, we consider the left regular representation of the affine morphisms at zero level, which immediately produces a canonical pair of Hilbert affine $P$-modules (which we call {\em regular}); here $P$ is not assumed to be of finite depth.
Interestingly, in the case of finite depth subfactor planar algebras, it turns out that any weight zero irreducible Hilbert affine $P$-module  is isomorphic to a submodule of one of the above regular  Hilbert affine $P$-modules.
We next analyze the finite von Neumann algebras generated by the affine morphisms at zero level in their GNS representations with respect to the faithful traces considered above. 
Moreover, to every left module over these von Neumann algebras, we uniquely associate a Hilbert affine $P$-module with weight zero; we use Connes fusion techniques and the above regular Hilbert affine $P$-modules at zero level for these constructions.

Section \ref{JWconjirrdep2} deals with the study of Hilbert affine modules over irreducible depth two subfactor planar algebras $P$ which (by the Ocneanu-Szymanski theorem \cite{Sz}) basically arise from actions of finite dimensional Kac algebras, the skein theory of which has been described in \cite{KLS}, \cite{DK}.
We recall the structure maps of the Kac algebras coming from $P$ and the definition of the quantum double $DH$ of a finite dimensional Hopf $*$-algebra $H$ (from \cite{Kas}).
We then construct an explicit isomorphism between the quantum double of  $P_{+2}$  and the affine morphism space at level one, $AP_{+1,+1}$.
Using this isomorphism and the normalized Haar functional on $DP_{+2}$, we build a Hilbert affine $P$-module $V$ which is generated by its $1$ space $V_1 = AP_{+1,+1}$ and contains all irreducible Hilbert affine $P$-modules.
This gives a one-to-one correspondence between the isomorphism classes of irreducible Hilbert affine $P$-modules and that of $V_1 \cong DP_{+2}$.
Thus, we prove the Jones-Walker conjecture in the case of irreducible depth two subfactors. We end this article with some questions related to the monoidal structure on affine modules and the Jones-Walker conjecture in a more general case.
\subsection*{Acknowledgements}
The second named author would like to thank Vaughan Jones for suggesting the problem of Section 3 and also for patiently listening to a preliminary version of a proof of Theorem \ref{main-theorem} during the von Neumann Algebras meeting at Oberwolfach in 2008.
The authors would also like to express sincere gratitude to the referee whose thorough report helped the authors immensely to improve the manuscript to its present form.
\section{Preliminaries}\label{prelim}
\subsection{Some useful notations}\label{usenot}
We will set up some notations about planar algebras which will be used in the forthcoming sections.
Although at first, these notations may seem to be a hindrance towards the readability to some extent, they will greatly reduce the need for drawing elaborate planar pictures again and again.
Understanding and decoding these notations correctly will definitely facilitate smooth reading.
We will not give the definition of planar algebra which can be found in \cite{Jon}; however, we will be consistent with the notation described in \cite{Gho}.
\begin{enumerate}
\item We will consider the natural binary operation on $\{-,+\}$ given by $++ := +$, $+-:=-$, $-+:=-$ and $--:=+$. Notations such as $(-)^l$ have to be understood in this context.
\item We will denote the set of all possible colors of discs in tangles by $\t{Col}:= \left\{ \vlon k: \vlon \in \{+,-\},\, k \in \mbb{N}_0\right\}$ where $\N_0 := \N \cup \{ 0 \}$.
\item In a tangle, we will replace (isotopically) parallel strings by a single strand labelled by the number of strings, and an internal disc with color $\vlon k$ will be replaced by a bold dot with the sign $\vlon$ placed at the angle corresponding to the distinguished boundary components of the disc.
For example,
\psfrag{e}{$\vlon$}
\includegraphics[scale=0.25]{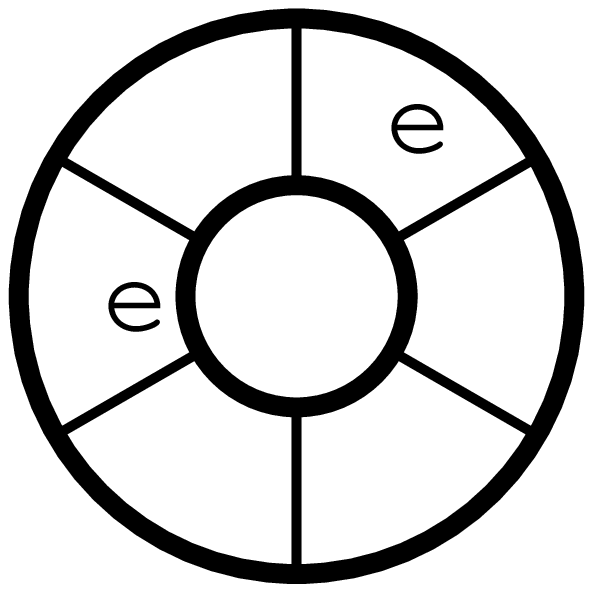}
will be replaced by
\psfrag{2}{$2$}
\psfrag{4}{$4$}
\psfrag{e}{$\vlon$}
\includegraphics[scale=0.25]{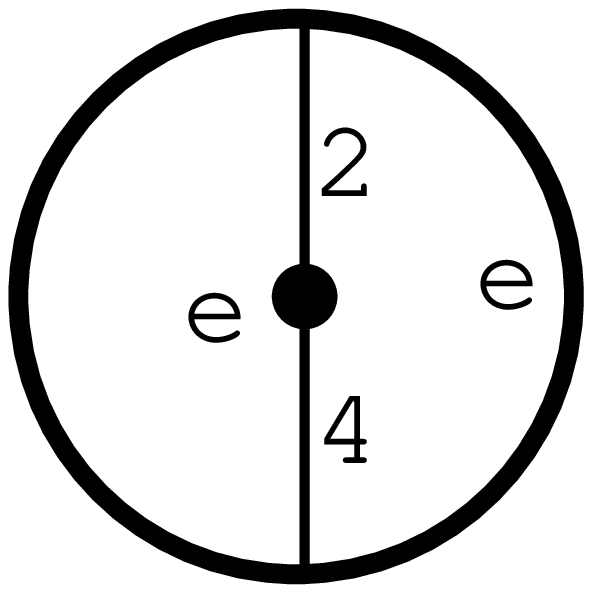}.
In a similar token, if $P$ is a planar algebra, we will replace a $P$-labelled internal disc by a bold dot with the label being placed at the angle corresponding to the distinguished boundary component of the disc; for instance,
\psfrag{x}{$x$}
\psfrag{e}{$\vlon$}
\includegraphics[scale=0.25]{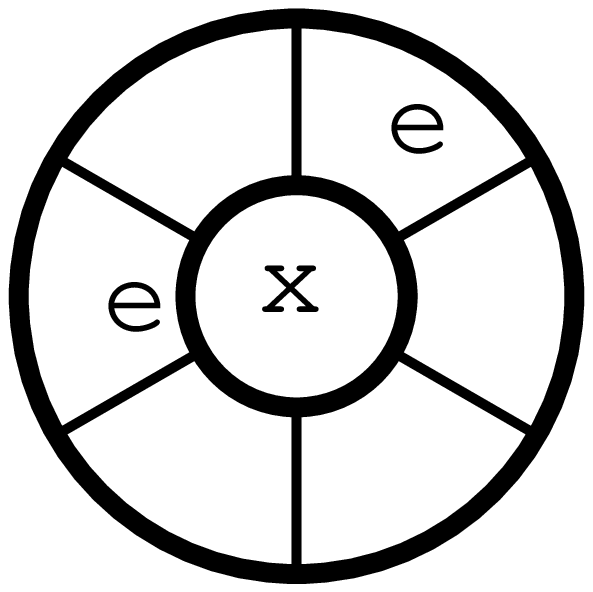}
will be replaced by
\psfrag{x}{$x$}
\psfrag{2}{$2$}
\psfrag{4}{$4$}
\psfrag{e}{$\vlon$}
\includegraphics[scale=0.25]{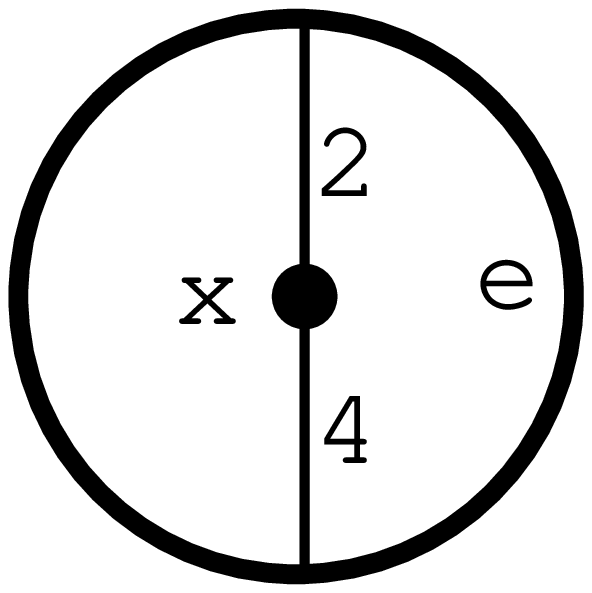}
where $x \in P_{\vlon 3}$.
\comments{
$\!\!\!\!\! \psfrag{replace}{will be replaced by}
\psfrag{2}{$2$}
\psfrag{4}{$4$} \psfrag{e}{$\vlon$}
\vcenter{\includegraphics[scale=0.25]{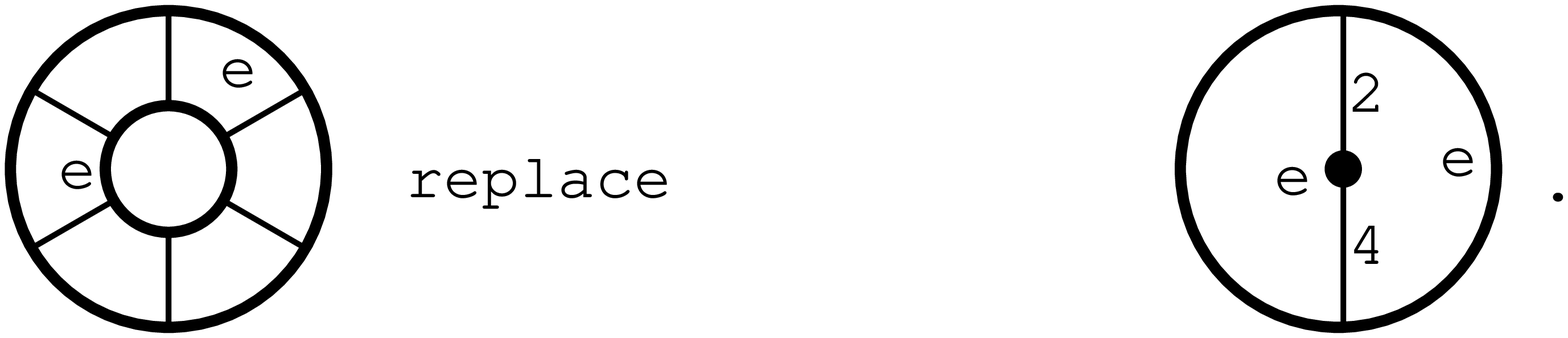}}$

$\!\!\!\!\! \psfrag{replace}{will be replaced by}
\psfrag{x}{$x$}
\psfrag{2}{$2$}
\psfrag{4}{$4$} \psfrag{e}{$\vlon$}
\vcenter{\includegraphics[scale=0.25]{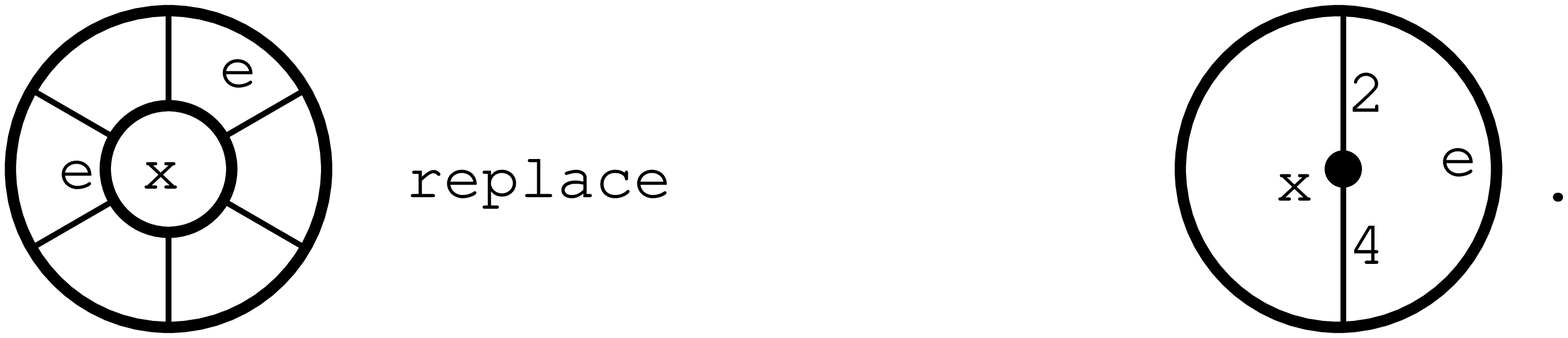}}$
}
We will reserve alphabets like $x, y, z$ to denote elements of $P$, $\vlon, \eta, \nu$ to denote a sign, and $k, l, m$ to denote a natural number to avoid confusion.
It should be clear from the context what a bold dot or a string in a picture is labelled by.
\item We set some notation for a set of `generating tangles' (that is, tangles which generate all tangles by composition) in Figure \ref{tangles}.
\begin{figure}[h]
\psfrag{mult}{Multiplication tangle}
\psfrag{unit}{Unit tangle}
\psfrag{identity}{Identity tangle}
\psfrag{jonesproj}{Jones projection tangle}
\psfrag{rightinclusion}{Right inclusion tangle}
\psfrag{leftinclusion}{Left inclusion tangle}
\psfrag{rightcondexp}{Right conditional expectation tangle}
\psfrag{leftcondexp}{Left conditional expectation tangle}
\psfrag{k}{$k$}
\psfrag{2k}{$2k$}
\psfrag{e}{$\vlon$}
\psfrag{-e}{$-\vlon$}
\psfrag{M}{$ M_{\vlon k} =$}
\psfrag{m}{$: (\vlon k,\vlon k)\rightarrow \vlon k$}
\psfrag{1ek}{$1_{\vlon k} =$}
\psfrag{1}{$: \emptyset \rightarrow \vlon k$}
\psfrag{RI}{$RI_{\vlon k} =$}
\psfrag{ri}{$: \vlon k \rightarrow \vlon (k+1)$}
\psfrag{LI}{$LI_{\vlon k}=$}
\psfrag{li}{$: \vlon k \rightarrow -\vlon (k+1)$}
\psfrag{Id}{$I_{\vlon k}=$}
\psfrag{id}{$: \vlon k \rightarrow \vlon k$}
\psfrag{E}{$E_{\vlon (k+1)}=$}
\psfrag{jp}{$: \emptyset \rightarrow \vlon (k+2)$}
\psfrag{RE}{$RE_{\vlon (k+1)} =$}
\psfrag{re}{$: \vlon(k+1) \rightarrow \vlon k$}
\psfrag{LE}{$LE_{\vlon (k+1)} =$}
\psfrag{le}{$: \vlon(k+1) \rightarrow -\vlon k$}
\includegraphics[scale=0.25]{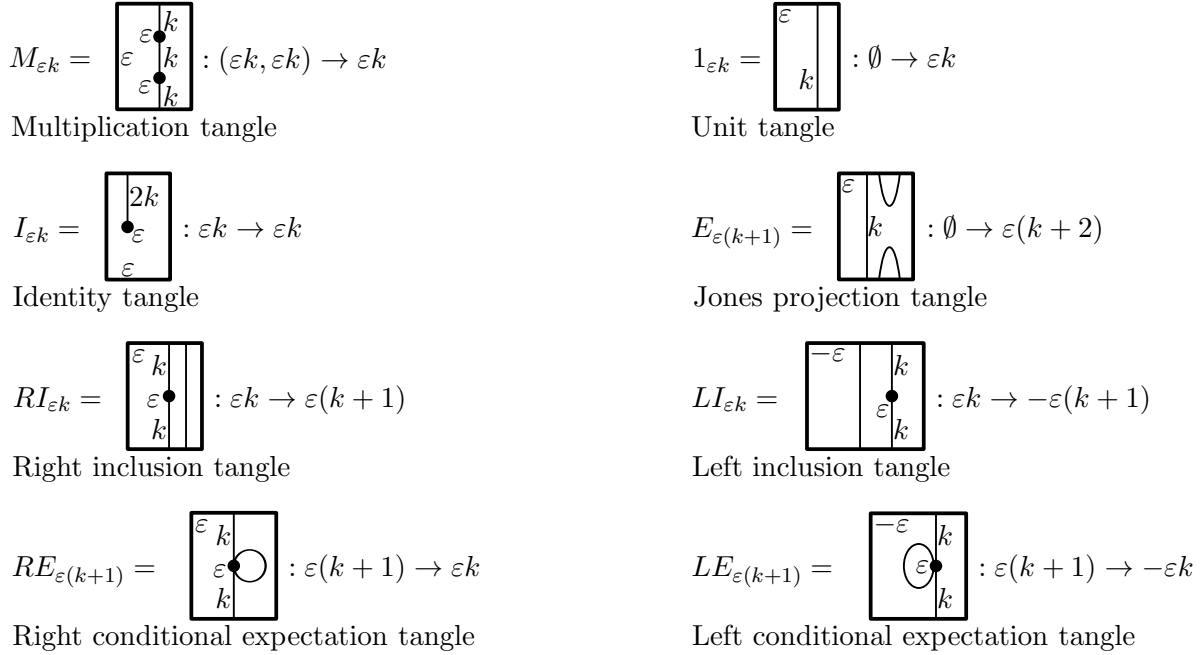}
\caption{Generating tangles.}
\label{tangles}
\end{figure}
\item ${\mcal T}_{\vlon k}$ (resp., ${\mcal T}_{\vlon k} (P)$) will denote the set of tangles (resp., $P$-labelled tangles) which has $\vlon k$ as the color of the external disc; ${\mcal P}_{\vlon k} (P)$ will be the vector space ${\mcal T}_{\vlon k} (P)$ as a basis.
The action of $P$ induces a linear map ${\mcal P}_{\vlon k} (P) \ni T \os{P}{\longmapsto} P_T \in P_{\vlon k}$.
\end{enumerate}  


We now recall (from \cite{Jon}) the notion of the {\em $n$-th cabling} of a planar algebra $P$ with modulus $(\delta_- , \delta_+)$, denoted by $c_n (P)$.
For a tangle $T$, let $c_n (T)$ be the tangle obtained by (a) replacing every string by $n$ many strings parallel to it, and (b) putting $n$ consecutive caps on the distinguished boundary component of every negatively signed (internal or external) disc and around the minus sign which is then replaced by a plus sign.

\noindent {\em Vector spaces:} For all colors $\vlon k$, $c_n
(P)_{\vlon k} := \t{Range} \,(P_{c_n(I_{\vlon k})})$.

\noindent {\em Action of tangles:} For all tangles $T$, $c_n (P)_T :=
\left[ \ous{n}{\prod}{l=1} \delta_{(-)^l} \right]^{- w } P_{c_n (T)}$
where $w$ is the number of negatively signed internal disc(s) of $T$.

Note that $c_1(P)$ is isomorphic to $P$, $c_m (c_n (P)) = c_{mn} (P)$
and $c_n (P)$ has modulus $\left( \ous{n}{\prod}{l=1} \delta_{(-)^l}
  \; , \; \ous{n}{\prod}{l=1} \delta_{(-)^{l+1}} \right)$.
\subsection{Planar algebras and subfactors}
In this section, we will recall certain basic facts about subfactors and its interplay with planar algebras.
For the rest of this section, let $M_{-1} := N \subset M
=: M_0$ be a subfactor with $\delta^2 : = [M:N] < \infty$ ($\delta
>0$) and $\{M_k\}_{k \geq 1}$ be a tower of basic constructions with
$\{e_k \in \mscr{P} (M_k) \}_{k \geq 1}$ being a set of Jones
projections.  For each $k \geq 1$, set $e_{[-1,k]} := \delta^{k
  (k+1)}(e_{k+1}e_k \cdots e_1) (e_{k+2} e_{k+1}\cdots e_2)\cdots
(e_{2k+1} e_{2k}$ $\cdots e_{k+1}) \in N^\prime \cap M_{2k+1}$, $
e_{[0, k]} :=\delta^{k (k-1)} (e_{k+1}e_k \cdots e_2) (e_{k+2}
e_{k+1}\cdots e_2)\cdots (e_{2k} e_{2k-1} \cdots e_{k+1}) \in M^\prime
\cap M_{2k}$ and $v_k := \delta^k e_k e_{k-1} \cdots e_1 \in N^\prime
\cap M_k$.  Then, the tower of $II_1$ factors $N \subset M_k\subset
M_{2k+1}$ (resp., $M \subset M_k \subset M_{2k}$) is an instance of
basic construction with ${e_{[-1,k]}}$ (resp., ${e_{[0,k]}}$) as Jones
projection, that is, there exists an isomorphism $ {\varphi_{-1, k}} :
M_{2k+1} {\lra} \mcal{L}_N(L^2(M_k))$ (resp., ${\varphi_{0, k}} :
M_{2k} {\lra} \mcal{L}_M(L^2(M_k))$) given by
\begin{eqnarray*}
\varphi_{-1, k}(x_{2k+1}) \hat{x}_{k} & = & 
\delta^{2(k+1)} E_{M_k}(x_{2k+1}x_{k}e_{[-1,k]})^{\widehat{}}  \\
\text{(resp., } \varphi_{0,k} (x_{2k}) \hat{x}_{k} &= &
\delta^{2k} E_{M_k}(x_{2k}x_{k}e_{[0,k]})^{\widehat{}} \text{ )}
\end{eqnarray*}
for all $x_i \in M_i$, $i = k,\, 2k,\, 2k+1$, which is identity
restricted to $M_k$ and sends $e_{[-1,k]}$ (resp., $e_{[0,k]}$) to the
projection with range $L^2(N)$ (resp., $L^2(M)$).  Also, $\varphi_{-1,
  k} (M_i'\cap M_{2k+1}) =\,_{M_i}\mcal{L}_N(L^2(M_k))$ (resp.,
$\varphi_{0, k}(M_i'\cap M_{2k}) = \, _{M_i}\mcal{L}_M(L^2(M_k))$) and
$\vphi_{0,k} = \left. \vphi_{-1,k} \right|_{M_{2k}}$ for all $k \geq 0
$, $-1 \leq i \leq k$.

We now state the `extended Jones' theorem' which provides an important
link between finite index subfactors and planar algebras. This was
first established for extremal finite index subfactors in
\cite{Jon}. Later, it was extended to arbitrary finite index
subfactors in \cite{Bur, JP, DGG}. As mentioned above, we will follow
the set up of \cite{DGG}.
\begin{thm}\label{jones-theorem}
  $P$ defined by $P_{\vlon k} = N' \cap M_{k-1} $ or $M' \cap M_{k}$
  according as $\vlon = +$ or $-$, has a unique unimodular bimodule
  planar algebra structure with the $\ast$-structure given by the
  usual $\ast$ of the relative commutants such that for each $k \in
  \N_0$,
\begin{enumerate}
\item the action of multiplication tangles is given by the usual
multiplication in the relative commutants,
\item the action of the left inclusion tangle $LI_{-k}$ is given by
  the usual inclusion $M' \cap M_k \subset N'\cap M_k$,
\item the action of the right inclusion tangle $RI_{+k}$ is given by
  the usual inclusion $ M_{k-1} \subset M_{k}$,
\item $P_{E_{+(k+1)}} = \delta e_{k+1}$,
\item $P_{LE_{+(k+1)}} = \delta^{-1} \us{i}{\sum} b_i^* x b_i $ for
  all $ x \in P_{+(k+1)}$,
\end{enumerate}
where $\{ b_i\}_i$ is a left Pimsner-Popa basis for the subfactor $N
\subset M$.  ($P$ will be referred as the planar algebra associated to
the tower $\{M_k\}_{k \geq -1}$ with Jones projections $\{e_k\}_{k
  \geq 1}$.)
\end{thm}
\begin{rem}\label{jones-theorem-rem}
  Apart from the action of the tangles given in conditions (1) - (5),
  it is also worth mentioning the actions of a few other useful
  tangles, namely,

(a) $P_{RE_{+k}} = \delta \left. E^{M_{k-1}}_{M_{k-2}}
\right|_{P_{+k}}$,

(b) $P_{TR^r_{+k}} = \delta^k \left. tr_{M_{k-1}} \right|_{P_{+k}}$,

(c) $\delta^{-k} P_{TR^l_{+{2l}}}$ (resp., $\delta^{-k}
P_{TR^l_{+(2l-1)}}$) is given by the trace on $P_{+2l} = N' \cap
M_{2l-1}$ (resp., $P_{+(2l-1)} = N' \cap M_{2l-2}$) induced by the
canonical trace on ${_N} {\mcal L} ( L^2 (M_{l-1}) )$ via the map
$\vphi_{-1,l-1}$ (resp., $\vphi_{0,l-1}$) where $TR^l_{\vlon k}$ (resp., $TR^r_{\vlon k}$) denotes the left (resp., right) trace tangle as described in Figure \ref{trace-tangles}.
\begin{figure}[h]
\psfrag{right}{Right trace tangle}
\psfrag{left}{Left trace tangle}
\psfrag{k}{$k$}
\psfrag{trr}{$TR^r_{\vlon k} :=$}
\psfrag{e}{$\vlon$}
\psfrag{-e}{$(-)^k \vlon$}
\psfrag{trl}{$TR^l_{\vlon k} :=$}
\psfrag{r=}{$: \vlon k \rightarrow \vlon 0$}
\psfrag{l=}{$:\vlon k \ra (-)^k \vlon 0$}
\includegraphics[scale=0.25]{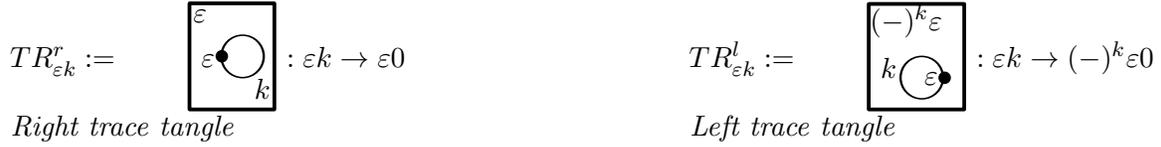}
\caption{Trace tangles.}
\label{trace-tangles}
\end{figure}
\end{rem}
\begin{cor}\label{dualpa}
  (a) $P_{E'_{-k}} (y) = \delta \us{i}{\sum} b^*_i e_1 y e_1 b_i$ for
  all $y \in P_{-k} = M' \cap M_k$, where $E'_{-k} = LI_{+(k-1)} \circ
  LE_{-k}$ and $\{ b_i \}_i$ is a left Pimsner-Popa basis for $N
  \subset M$,

  (b) the $n$-th dual of $P$, $\lambda_n (P) =$ the planar algebra
  associated to the tower $\{ M_{k+n} \}_{k \geq -1}$ with Jones
  projections $\{ e_{k+n} \}_{k \geq 1}$.
\end{cor}
If $e_{[l,k+l]}$ denotes the projection obtained by replacing each
$e_\bullet$ in the defining equation of $e_{[0,k]}$ (as above), by
$e_{l+\bullet}$, then $M_l \subset M_{k+l} \subset M_{2k+l}$ is an
instance of basic construction with $e_{[l,k+l]}$ as Jones projection.
\begin{rem}\label{cabling}
  An easy consequence of Corollary \ref{dualpa} (b) and Theorem
  \ref{jones-theorem} is $c_n(P) =$ the planar algebra associated to
  the tower $\{ M_{n(k+1)-1} \}_{k \geq -1}$ with Jones projections
  $\{e_{[n(k-1)-1,nk-1]} \}_{k \geq 1}$.
\end{rem}
\begin{prop}\label{rotcont} 
  If $J_k$ denotes the canonical conjugate-linear unitary operator on
  $L^2 (M_k)$ and $R^m_{\vlon n}$ denotes the tangle
  \psfrag{e}{$\vlon$} \psfrag{m}{$m$} \psfrag{2n-m}{$2n-m$}
  \psfrag{-e}{$(-)^m \vlon$}
\includegraphics[scale=0.2]{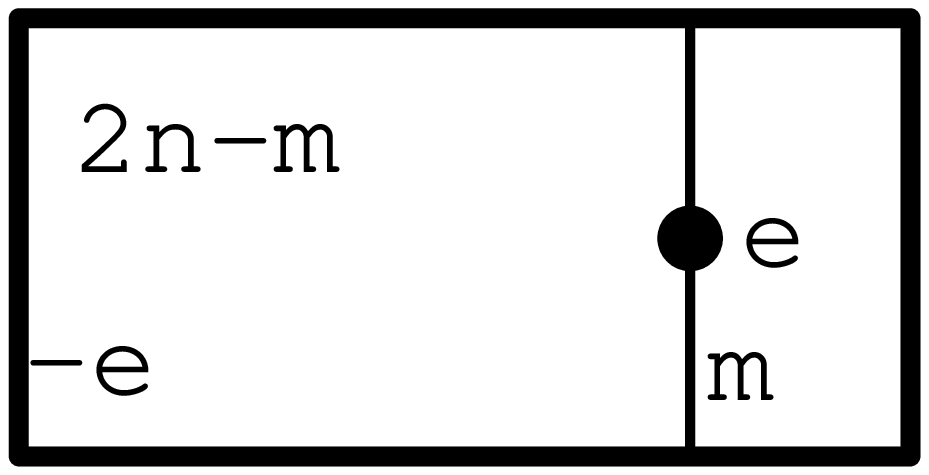}
, then for all $k \geq 0$, we have:

(a) $\vphi_{-1,k} \left( P_{R^{2k+2}_{+(2k+2)}} (x) \right) = J_k
\vphi_{-1,k} (x^*) J_k \;$ and $\; \t{Range} \,
\vphi_{-1,k} \left( P_{R^{2k+2}_{+(2k+2)}} (p) \right) \os{N \text{-}
  N}{\cong} \oline{\t{Range} \,\vphi_{-1,k} (p)}$ for all $x \in P_{+(2k+2)}$, $p \in {\mathscr P}
(P_{+(2k+2)})$,

(b) $\vphi_{-1,k} \left( P_{R^{2k+1}_{+(2k+1)}} (x) \right) = J_k
\vphi_{0,k} (x^*) J_k \;$ and $\; \t{Range} \,
\vphi_{-1,k} \left( P_{R^{2k+1}_{+(2k+1)}} (p) \right)$$\os{M \text{-}
  N}{\cong}$$\oline{\t{Range} \,\vphi_{-1,k} (p)}$ for all $x \in P_{+(2k+1)}$, $p \in {\mathscr P}
(P_{+(2k+1)})$,

(c) $\vphi_{0,k} \left( P_{R^{2k+1}_{-(2k+1)}} (x) \right) = J_k
\vphi_{-1,k} (x^*) J_k \;$ and $\; \t{Range} \,
\vphi_{0,k} \left( P_{R^{2k+1}_{-(2k+1)}} (p) \right) \os{N \text{-}
  M}{\cong} \oline{\t{Range} \,\vphi_{-1,k} (p)}$ for all $x \in P_{-(2k+1)}$, $p \in {\mathscr P}
(P_{-(2k+1)})$,

(d) $\vphi_{0,k} \left( P_{R^{2k}_{-2k}} (x) \right) = J_k \vphi_{0,k}
(x^*) J_k \;$ and $\; \t{Range} \,\vphi_{0,k} \left(
  P_{R^{2k}_{-2k)}} (p) \right) \os{M \text{-} M}{\cong} \oline{\t{Range} \,
  \vphi_{0,k} (p)}$ for all $x \in P_{-2k}$, $p \in {\mathscr P} (P_{-2k})$.
\end{prop}
\begin{pf}
  The isomorphism in the second part in each of (a), (b), (c) and (d),
  follows from the first part using \cite[Proposition $3.11$]{Bis97}.
  For the first parts, it is enough to establish only for (a) because
  all others can be deduced using conditions (2) and (3) in Theorem
  \ref{jones-theorem}, and the relation $\vphi_{0,k} =
  \left. \vphi_{-1,k} \right|_{M_{2k}}$.

  First, we will prove part (a) for $k=0$.  Note that if $\{b_i\}_i$
  is a left Pimsner-Popa basis for $N \subset M$, then
 \[
P_{R^2_{+2}} (x) = P_{R^1_{-2}} \left( P_{R^1_{+2}} (x) \right) =
P_{R^1_{-2}} \left( \delta \us{i}{\sum} b^*_i x e_2 e_1 b_i \right) =
\delta^4 \us{i}{\sum} E_{M_1} \left( e_2 e_1 b^*_i x e_2 e_1 b_i \right)
\]
where we use the conditions of Theorem \ref{jones-theorem} in a
decomposition of the rotation tangle into the generating ones $R^1_{+2}
= LE_{+3} \circ M_{+3} ( RI_{+2} , M_{+3} (E_{+2} , RI_{+2} \circ
E_{+1}) )$ (resp., $RE_{+3} \circ M_{+3} (M_{+3} (E_{+2} , RI_{+2}
\circ E_{+1}) , LI_{-2})$) for establishing the second (resp., third)
equality. For $y \in M$, note that
\begin{align*}
\vphi_{-1,0} \left( P_{R^2_{+2}} (x) \right) \hat{y} & =
\delta^6 \sum_i E_M ( e_2 e_1 b_i^* x e_2 e_1 b_i y e_1 ) \hat{}
= \delta^6 E_M ( e_2 e_1 y x
e_2 e_1 ) \hat{} = \delta^2 E_M (e_1 y x)\hat{}\\
& = J_0 \vphi_{-1, 0}
(x^*) J_0 \hat{y}.
\end{align*}
Now, let $k > 0$. Using the above and
Remark \ref{cabling}, we obtain
\[
\vphi_{-1,k} \left( P_{R^{2k+2}_{+(2k+2)}} (x) \right) = \vphi_{-1,k} \left( {c_{k+1}(P)}_{R^2_{+2}} (x) \right) = J_k \vphi_{-1 , k} (x^*) J_k.
\]
\end{pf}
\vspace{2mm}

We will make repeated use of the following standard facts, whose proof
can be found in \cite{Bis97}.
\begin{lem} \cite{Bis97}\label{e-even-odd}
For each $k \geq 0$, and $ X \in \{N, M\}$, we have:
\begin{enumerate}
\item $\t{Range} \,\varphi_{-1, k} (p) \os{X-N}{\cong} \t{Range} \,\varphi_{-1, k+1}
  (pe_{2k+3})$ for all $p \in \mscr{P} (X' \cap M_{2k +1})$,
\item $\t{Range} \,\varphi_{0, k} (p) \stackrel{X-M}{\cong} \t{Range} \,\varphi_{0,
    k+1} (pe_{2k+2})$ for all $p \in \mscr{P} (X' \cap M_{2k})$.
\end{enumerate}
\end{lem}
\noindent From this, one can easily deduce the following.
\begin{cor}\label{p-q-equivalence}
For $k > l \geq 0$ and $X \in \{N, M\}$, the following holds:
\begin{enumerate}
\item For all $p \in \mscr{P} (X' \cap M_{2k +1})$ and $q \in \mscr{P}
  (X' \cap M_{2l + 1})$ satisfying $\t{Range} \,\varphi_{-1, k} (p)
  \stackrel{X-N}{\cong} \t{Range} \,\varphi_{-1, l} (q)$, $p$ is
  MvN-equivalent to $q e_{2l+3}\cdots e_{2k+1}$ in $X' \cap M_{2k +
    1}$.
\item For all $p \in \mscr{P} (X' \cap M_{2k})$ and $q \in \mscr{P}
  (X' \cap M_{2l})$ satisfying $\t{Range} \,\varphi_{0, k} (p)
  \stackrel{X-M}{\cong} \t{Range} \,\varphi_{0, l} (q)$, $p$ is MvN-equivalent
  to $q e_{2l+2}\cdots e_{2k}$ in $X' \cap M_{2k}$.
\end{enumerate}
\end{cor}
\subsection{Affine Category over a  Planar Algebra}\label{affrep}
\begin{figure}[h]
\psfrag{AR}{$AR_{\vlon k}$}
\psfrag{y}{$\eta$}
\psfrag{2k-1}{$2k-1$}
\psfrag{2k}{$2k$}
\psfrag{2l}{$2l$}
\psfrag{m}{$m$}
\psfrag{e}{$\varepsilon$}
\psfrag{-e}{$-\varepsilon$}
\psfrag{psi}{$\Psi^m_{\varepsilon k, \eta l}$}
\psfrag{1ek}{$A1_{\varepsilon k}$}
\includegraphics[scale=0.30]{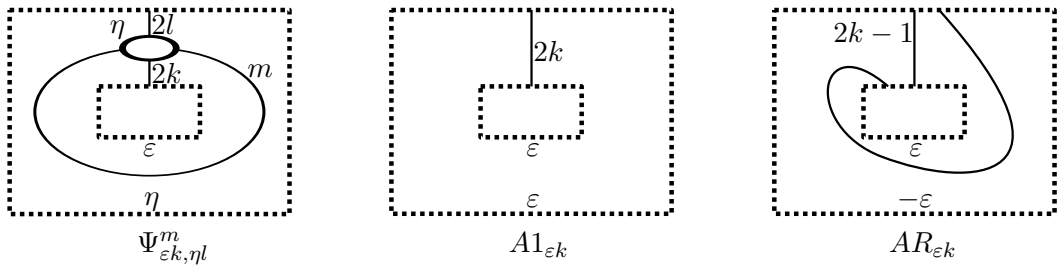}
\caption{Some useful affine tangles. ($\varepsilon , \eta \in \{+, -
  \}, k, l \in \N_0, m \in (2\N_0 + \delta_{\vlon \neq
    \eta})$) \label{aff-tangles}}
\end{figure}

In this subsection, for the sake of self containment, we recall (from
\cite{Gho}) in some detail what we mean by the \textit{affine category
  over a planar algebra} and the corresponding {\em affine morphisms}
(with slight modifications).
\begin{defn}
  For each $ \varepsilon , \eta \in \{ +, - \}$ and $k, l \geq 0$, an
  $(\varepsilon k , \eta l)$-affine tangular picture consists of the
  following:
\begin{itemize}
\item finitely many (possibly none) non-intersecting subsets $D_{1},
  \cdots , D_{b}$ (referred as disc(s)) of the interior of the
  rectangular annular region $RA:= [-2, 2]\times [-2, 2] \setminus
  (-1, 1)\times (-1, 1)$, each of which is homeomorphic to the unit
  disc and has even number of marked points on its boundary, numbered
  clockwise,
\item non-interescting paths (called strings) in $RA \setminus \left[
    \ous{b}{\sqcup}{i=1} \mathrm{Int} (D_i) \right]$, which are either
  loops or meet the boundaries of the discs or $RA$ exactly at two
  distinct points in $\{ (\frac{i}{2k}, 1): 0 \leq i \leq 2k-1\}
  \sqcup \{ (\frac{j}{2l}, 2) : 0 \leq j \leq 2l-1\} \sqcup
  \{\text{marked points on the discs} \}$ in such a way that every
  point in this set must be an endpoint of a string,
\item a checker-board shading on the connected components of $Int(RA)
  \setminus \left[ \left( \ous{b}{\sqcup}{i=1} D_i\right)\, \cup
    \{\mathrm{strings}\, \} \right]$ such that the component near the
  point $(0,-1)$ (resp., $(0,-2)$) is unshaded or shaded according as
  $\vlon$ (resp., $\eta$) is $+$ or $-$.

\end{itemize}
\end{defn}
\begin{defn}\label{aff-isotopy}
  An affine isotopy of an affine tangular picture is a map $\varphi :
  [0, 1] \times RA \ra RA$ such that
\begin{enumerate}\vspace*{-2mm}
\item $\varphi(t, \cdot)$ is a homeomorphism of $RA$, for all $t \in
  [0, 1]$;
\item $\varphi(0, \cdot) = id_{RA}$; and
\item $\left. \varphi(t, \cdot) \right|_{\partial (RA)} = id_{\partial
    (RA)}$ for all $t \in [0, 1]$.
\end{enumerate}
\end{defn}
Two affine tangular pictures are said to be {\em affine isotopic} if one can be obtained from the other using an affine isotopy preserving checker-board shading and the distinguished boundary components of the discs.
It may be noted here that condition (3) in Definition \ref{aff-isotopy} distinguishes affine isotopy from {\em annular isotopy} (see \cite{JR06}, \cite{Gho}).
\begin{defn}
An $(\varepsilon k , \eta l )$-affine tangle is the affine isotopy class of an $(\varepsilon k , \eta l )$-affine tangular picture.
\end{defn}
Time and again, for the sake of convenience, we will abuse terminology by referring to an affine tangular picture as an affine tangle (corresponding to its affine isotopy class) and the figures might not be sketched to the scale but are clear enough to avoid any ambiguity.
In Figure \ref{aff-tangles}, we draw a specific affine tangle called $\Psi_{\varepsilon k, \eta l}^m$, where the labels next to strings have the same significance as that explained in (3) of Section \ref{usenot}.
This affine tangle will play an important role in the following discussions.

\noindent {\bf Notations:} For each $ \varepsilon , \eta \in \{ +, - \}$ and $k, l \geq 0$, let
\begin{itemize}
\item $\mathcal{AT}_{\varepsilon k, \eta l}$ denote the set of all $(\varepsilon k, \eta l)$-affine tangles, and
\item $\mcal{A}_{\varepsilon k, \eta l}$ denote the complex vector space with $\mcal{AT}_{\varepsilon k, \eta l}$ as a basis.
\end{itemize}
The composition
of affine tangles $T \in \mathcal{AT}_{\varepsilon k, \eta l}$ and $S
\in \mathcal{AT}_{\xi m, \varepsilon k}$ is given by $T \circ S :=
\frac{1}{2}(2T \cup S) \in \mcal{AT}_{\xi m , \eta l}$; this
composition is linearly extended to the level of the vector spaces $\mcal{A}_{\vlon k, \eta l}$'s.

\begin{rem} \label{Psi-remark}
For each $A \in \mathcal{AT}_{\varepsilon k, \eta l}$, there is an $m \in (2\mbb{N}_0 + \delta_{\vlon \neq \eta})$ and a  $T \in \mcal{T}_{\eta(k + l +m)}$ such that $A = \Psi_{\varepsilon k, \eta l}^m (T)$ where $\Psi_{\varepsilon k, \eta l}^m (T)$ is the isotopy class of the affine tangular picture obtained by inserting $T$ in the disc of $\Psi_{\varepsilon k, \eta l}^m$.
\end{rem}
In the above remark, the $m$ can be chosen as large as one wants and the insertion method extends linearly to a linear map $\Psi_{\vlon k, \eta l}^m : \mcal{P}_{\eta (k + l + m)} \ra \mcal{A}_{\varepsilon k, \eta l}$, and for each $A \in \mathcal{A}_{\varepsilon k, \eta l}$, there is an $m \in \mbb{N}_0 $ and an $X \in \mcal{P}_{\eta(k + l +m)}$ such that $A = \Psi_{\varepsilon k, \eta l}^m (X)$.
Let $P$ be a planar algebra.
An $( \varepsilon k, \eta l)$-affine tangle is said to be $P$-{\em labelled} if each disc is labelled by an element of $P_{\nu m}$ where $\nu m$ is the color of the corresponding disc.
Let $\mathcal{AT}_{\varepsilon k, \eta l}(P)$ denote the collection of all $P$-{\em labelled} $(\varepsilon k, \eta l)$-affine tangles, and let $\mathcal{A}_{\varepsilon k, \eta l}(P)$ be the vector space with $\mathcal{AT}_{\varepsilon k, \eta l}(P)$ as a basis.
Composition of $P$-{\em labelled} affine tangles also makes sense as above and extends to their complex span.
Note that $\Psi_{\varepsilon k, \eta l}^m$ also induces a linear map from $\mathcal{P}_{\eta (k+l+m)}(P)$ into $\mathcal{A}_{\varepsilon k, \eta l}(P)$.
Moreover, from \ref{Psi-remark}, we may conclude that, for each $A \in \mathcal{A}_{\varepsilon k, \eta l}(P)$, there is an $m \in (2\mbb{N}_0 + \delta_{\vlon \neq \eta})$ and an $X \in \mcal{P}_{\eta(k + l +m)}(P)$ such that $A = \Psi_{\varepsilon k, \eta l}^m (X)$.

\vspace{2mm}
Now, consider the set $\mathcal{W}_{\varepsilon k,\eta l} := \us{m \in
  \N_0}{\cup} \left\{ \Psi _{\varepsilon k, \eta l }^{m}(X) : X \in
  \mathcal{P}_{\eta (k+l+m)}(P) \text{ s.t. } P_X = 0 \right\}$. It is
straight forward - see \cite{Gho} - to observe that
$\mathcal{W}_{\varepsilon k,\eta l}$ is a vector subspace of
$\mathcal{A}_{\varepsilon k, \eta l}(P)$.

Define the category $AP$ by:
\begin{itemize}
\item $\text{ob}(AP) := \{\varepsilon k : \varepsilon \in \{+,-\}, k\in
  \mathbb{N}_{0}\}$,
\item $\text{Mor}_{AP}(\varepsilon k, \eta l) :=
  \frac{\mathcal{A}_{\varepsilon k,\eta l}
    (P)}{\mathcal{W}_{\varepsilon k, \eta l}}$ (also denoted by
  $AP_{\varepsilon k,\eta l}$),
\item composition of morphisms is induced by the composition of
  $P$-{\em labelled} affine tangles (see \cite{Gho}),
\item the identity morphism of $\varepsilon k$ is given by the class
  $[A1_{\varepsilon k}]$ - Figure \ref{aff-tangles}.
\end{itemize}
$AP$ is a $\mbb C$-linear category and is called the \textit{affine category over $P$} and the morphisms in this category are called {\em affine morphisms}.

For $\vlon k, \eta l \in \t{Col}$ and $m \in (2\N_0 + \delta_{\vlon \neq \eta})$, consider the composition map
\[
\psi_{\vlon k , \eta l}^{m} \; : \; P_{\eta (k+l+m)} \stackrel{ I_{\eta (k+l+m)}}{\lra} \mcal{P}_{\eta (k+l+m)}(P) \stackrel{ \Psi^m_{\vlon k, \eta l}}{\lra} \mcal{A}_{\vlon k, \eta l}(P) \stackrel{q_{\vlon k, \eta l}}{\ra} AP_{\vlon k, \eta l}
\]
where the map $I_{\eta (k+l+m)} : P_{\eta (k+l+m)} \ra \mcal{P}_{\eta (k+l+m)}(P)$ is obtained by labelling the internal disc of the identity tangle $I_{\eta (k+l+m)}$ (defined in Figure \ref{tangles}) by a vector in $P_{\eta (k+l+m)}$, and $q_{\vlon k , \eta l}: \mcal{A}_{\vlon k, \eta l}(P) \ra AP_{\vlon k, \eta l}$ is the quotient map. Note that $\psi^m_{\vlon k, \eta l}$ is indeed linear, although $I_{\eta (k+l+m)}$ is not.
\begin{rem}\label{psi-remark}
For each $a \in AP_{\vlon k, \eta l}$, there exists $m \in (2\mbb{N}_0 + \delta_{\vlon \neq \eta})$ and  $x \in P_{\eta(k + l +m)}$ such that $a = \psi_{\vlon k, \eta l}^m (x)$.
\end{rem}

\noindent {\bf $*$-structure:} If $P$ is a $\ast $-planar algebra, then each $\mathcal{P}_{\varepsilon k }(P)$ becomes a $\ast $-algebra where $\ast $ of a {labelled} tangle is given by $\ast $ of the unlabelled tangle whose internal discs are labelled with $\ast$ of the labels.
Further, one can define $\ast$ of an affine tangular picture by reflecting it inside out such that the reflection of the distinguished boundary segment of any disc becomes the same for the disc in the reflected picture; this also extends to the $P$-{\em labelled} affine tangles as in the case of $P$-{\em labelled} tangles.
Clearly, $\ast$ is an involution on the space of $P$-labelled affine tangles, which can be extended to a conjugate linear isomorphism $\ast : \mcal{A}_{\varepsilon k, \eta l} (P) \rightarrow \mcal{A}_{\eta l ,  \varepsilon k} (P)$ for all colours $\varepsilon k$, $\eta l$.
Moreover, it is readily seen that $\ast \left( \mathcal{W}_{\varepsilon k, \eta l}\right) =\mathcal{W}_{\eta l,\varepsilon k}$; so the category $AP$ inherits a $\ast$-category structure.
\begin{defn}
Let $P$ be a planar algebra.
\begin{enumerate}
\item A $\C$-linear functor $V$ from $AP$ to $\mcal{V}ec$ (the category of complex vector spaces), is said to be an `affine $P$-module', that is, there exists a vector space $V_{\vlon k}$ for each $\vlon k \in \text{Col}$ and a linear map $AP_{\vlon k , \eta l} \ni a \os{V}{\longmapsto} V_a \in \text{Mor}_{\mcal{V}ec} (V_{\vlon k} , V_{\eta l}) $ for every $\vlon k , \eta l \in \text{Col}$ such that compositions and identities are preserved. ($V_a$ will be referred as the action of the affine morphism $a$.)

\item Further, for a $*$-planar algebra $P$, $V$ is called a `$*$-affine $P$-module' if it preserves the $*$-structure of $AP$, that is,
$\langle \xi , a \eta \rangle = \langle a^* \xi , \eta \rangle$ for all affine morphisms $a$, and $\xi$ and $\eta$ in appropriate $V_{\vlon k}$'s.
\item A $*$-affine $P$-module $V$ will be called `Hilbert affine $P$-module' if $V_{\vlon k}$'s are Hilbert spaces.
\end{enumerate}
\end{defn}
An affine module is said to be {\em bounded} (resp., {\em locally finite}) if the actions of the affine morphisms are bounded operators with respect to the norm coming from the inner product (resp., $V_{\vlon k}$'s are finite dimensional).
By closed graph theorem, a Hilbert affine $P$-module is automatically bounded; conversely, every bounded $*$-affine $P$-module can be completed to a Hilbert affine module.

Below, we give a list of some standard facts on Hilbert affine $P$-modules for a $*$-planar algebra $P$ with modulus. The proofs of the facts, as stated here, are straight-forward exercises (for analogous statements on annular tangles see \cite{Jon01}).
If $V$ is a Hilbert affine $P$-module and $S_{\vlon k}\subset V_{\vlon k}$ for $\vlon k\in\mbox{Col}$, then one can consider the `submodule of $V$ generated by $S=\us{\vlon k\in\mbox{Col}}{\coprod}S_{\vlon k}$' (denoted by $[S]$) given by $\left\{ \left[S\right]_{\eta l}:=\ol{\mbox{span}\left\{ \us{\varepsilon k\in\mbox{Col}}{\cup}AP_{\vlon k,\eta l}\left(S_{\varepsilon k}\right)\right\} }^{\parallel\cdot\parallel}\right\} _{\eta l}$ which
is also the smallest submodule of $V$ containing $S$.
\begin{rem}\label{irrmod}
Let $V$ be a Hilbert affine $P$-module and $W$ be an $AP_{\vlon k,\vlon k}$-submodule
of $V_{\vlon k}$ for some $\vlon k\in\mbox{Col}$. Then,
\begin{enumerate}
\item $V$ is irreducible if and only if $V_{\vlon k}$ is irreducible $AP_{\vlon k,\vlon k}$-module
for all $\vlon k\in\mbox{Col}$ if and only if $\left[v\right]=V$
for all $0\neq v\in V$.
\item $W$ is irreducible $\Leftrightarrow$ $\left[W\right]$ is an irreducible
submodule of $V$.\end{enumerate}
\end{rem}
\begin{rem}\label{affmodmor}
If $V$ and $W$ are Hilbert affine $P$-modules for which there exists an $\vlon k\in\mbox{Col}$ such that $V=\left[V_{\vlon k}\right]$ and there exists an $AP_{\vlon k,\vlon k}$-linear isometry $\theta:V_{\vlon k}\ra W_{\vlon k}$ , then $\theta$ extends uniquely to an isometry (of Hilbert affine $P$-modules) $\tilde{\theta}:V\ra W$.
\end{rem}
For $\vlon = \{ +, - \}$, we will also consider \emph{Hilbert $\vlon$-affine
$P$-module} $V$ consisting of the Hilbert spaces $V_{\pm0},V_{1},V_{2},\ldots$ equipped with a $*$-preserving action of affine morphisms as follows: 
\[
\left.
\begin{tabular}{l}
$AP_{\vlon k,\vlon l} \times V_k \ra V_l$\\
$AP_{\vlon k,\eta 0} \times V_k \ra V_{\eta 0}$\\
$AP_{\eta 0, \vlon l} \times V_{\eta 0} \ra V_l$\\
$AP_{\eta 0, \nu 0} \times V_{\eta 0} \ra V_{\nu 0}$
\end{tabular}
\right\} \t{ for all } k, l \in \N,\; \eta, \nu \in \{\pm\}.
\]
\begin{rem}\label{eaff}
The restriction map from the set of isomorphism classes of Hilbert
affine $P$-modules to that of the Hilbert $\vlon$-affine $P$-modules,
is a bijection.
\end{rem}
To see this, consider an irreducible Hilbert $+$-affine $P$-module $V$.
Define $\t{Ind} \, V_{\vlon k}:=V_{k}$ and $\t{Ind} \, V_{\vlon0}:=V_{\vlon0}$
(as Hilbert spaces) and the action of affine morphisms by  $AP_{\vlon k,\eta l}\times \t{Ind} \, V_{\vlon k} \ni (a,v) \longmapsto \left(r_{\eta l}^{-1}\circ a\circ r_{\vlon k}\right)\cdot v\in \t{Ind} \, V_{\eta l}$
where $r_{\vlon k}=\left\{ \begin{array}{ll}
A1{}_{\vlon k}, & \mbox{ if }k=0\mbox{ or }\vlon=+,\\
AR_{\vlon k}, & \mbox{ otherwise;}\end{array}\right.$\\
$A1{}_{\vlon k}$ and $AR_{\vlon k}$ being the affine tangles mentioned in Figure \ref{aff-tangles}.

\vspace{2mm}
For every affine $P$-module $V$,$\mbox{dim}\left(V_{+k}\right)=\mbox{dim}\left(V_{-k}\right)$ for all $k\geq1$ and it increases as $k$ increases. This motivates the following definition:
\begin{defn}
The `weight of $V$' is defined to be the smallest non-negative integer $k$ such that $V_{+k}$ or $V_{-k}$ is nonzero.
\end{defn}
\section{Affine morphisms at zero level}\label{mt}
In this section, we will be interested in understanding the
{\em affine morphisms at zero level} of a $\ast$-planar algebra $P$,
that is, in the space
\[
AP_{0, 0} :=
\begin{bmatrix}
AP_{+ 0, + 0} & AP_{- 0, + 0} \\
AP_{+0, - 0} & AP_{- 0, - 0}
\end{bmatrix},
\] 
which has a natural $*$-algebra structure induced by matrix
multiplication with respect to composition of affine morphisms and the
$\ast$ as discussed before.  On the other hand, given a finite index
subfactor $N \subset M$, for each $\vlon, \eta \in \{+, -\}$, we set $
V_{\vlon, \eta} := \{ \text{isomorphism classes of irreducible }
X_\eta \text{-} X_\vlon$ bimodules appearing in the standard
invariant$\} = \{$isomorphism classes of irreducible sub-bimodules of
${_{X_\eta}} L^2(M_k)_{X_\vlon} \text{ for some } k \in \N_0 \}$ where
$X_+$ (resp., $X_-$) denotes $N$ (resp., $M$).  Then, the usual matrix
multiplication with respect to appropriate relative tensor products
and the matrix adjoint with respect to the contragradients of
bimodules induce a natural $\ast$-algebra structure on the space
\[
\mcal{F}_{N \subset M} := \begin{bmatrix} \mbb{C}V_{+, +} &
  \mbb{C}V_{-, +}\\ \mbb{C}V_{+, -} & \mbb{C}V_{-, -}
\end{bmatrix}.
\] 
We will aim to prove the following:
\begin{thm}\label{main-theorem}
Let $N \subset M$ be a finite index subfactor and $P$ be its associated planar algebra. Then,
\[
AP_{0, 0}  \cong \mcal{F}_{N \subset M}
\] as $\ast$-algebras.
\end{thm}
\subsection{A spanning set for $AP_{0, 0}$}
In this subsection, $P$ will always denote the planar algebra associated to the tower of basic construction $\{ M_k \}_{k \in \N}$ of a finite index subfactor $N \subset M$ with Jones projections $\{ e_k \}_{k \in \N}$, and $\psi^m_{\vlon , \eta}$ will denote the linear map $\psi^m_{\vlon 0 , \eta 0}$ introduced right before Remark \ref{psi-remark}.
We first list some elementary yet useful properties of the $\psi$-maps.
\begin{lem}\label{psi-faithful}
For $\vlon, \eta \in \{+, -\}$ and $k \in (2\N_0 + \delta_{\vlon \neq \eta})$, $\psi_{\vlon 0, \eta 0}^{k} (p) \neq 0$ for all nonzero $p \in \mscr{P}(P_{\eta k})$.
\end{lem}
\begin{pf}
Let $\omega_{\vlon, \eta} : \mcal{A}_{\vlon 0, \eta 0} (P) \ra \mcal{P}_{\eta 0} (P)$ be the map defined by sending an affine tangle $[A] \in {\mcal AT}_{\vlon 0 , \eta 0}$ to the tangle obtained by ignoring the internal rectangle in $A$.
Note that  ${\mcal W}_{\vlon 0, \eta 0} \subset \t{ker}\, (P \circ \omega_{\vlon, \eta})$; thus, $P \circ \omega_{\vlon, \eta}$ induces a linear map  $\omega'_{\vlon, \eta} : AP_{\vlon 0, \eta 0} \ra P_{\eta 0} \cong \mbb{C}$.
Clearly, $\omega'_{\vlon, \eta} \circ \psi_{\vlon 0, \eta 0}^k = P_{TR^r_{\eta k}}$.
This proves the lemma.
\end{pf}
\begin{lem}\label{psi-properties} Let $\vlon, \eta \in \{+, -\}$ and
  $k \in (2\N_0 + \delta_{\vlon \neq \eta})$.
\begin{enumerate}
\item The map $\psi^k_{\vlon, \eta}$ is tracial, (that is,
  $\psi^k_{\vlon, \eta} (xy) = \psi^k_{\vlon, \eta} (yx)$ for all $x ,
  y \in P_{\eta k}$) and hence, factors through the center of $P_{\eta
    k}$.
\item $\psi^k_{\vlon, \eta} (x) = \psi^{k+2}_{\vlon, \eta} (x
  e_{(k+1+\delta_{\eta = -})})$ for all $x \in P_{\eta k}$.
\end{enumerate}
\end{lem}
\begin{pf}
  Both follow from simple application of affine isotopy and also using
  the relation between the Jones projections and the Jones projection
  tangles, in (2).
\end{pf}

From Corollary \ref{p-q-equivalence} and Lemma \ref{psi-properties},
we deduce the following where, for convenience, we use $\vphi_{\vlon
  k}$ to denote $\vphi_{-1 , \frac{k}{2}-1 }$ or $\vphi_{0 ,
  \frac{k-1}{2}}$ (resp., $\vphi_{0 , \frac{k}{2}}$ or $\vphi_{-1 ,
  \frac{k-1}{2}}$) according as $k$ is even or odd if $\vlon = +$
(resp., $\vlon = -$).
\begin{cor}\label{psi-corollary}
  Let $\vlon,\eta \in \{+, -\}$ and $k, l \in (2\mbb{N}_0 +
  \delta_{\vlon \neq \eta})$.  Then, for all $p \in \mscr{P}(P_{\eta
    k})$ and $q \in \mscr{P}(P_{\eta l})$ satisfying $\t{Range} \,
  \varphi_{\eta k} (p) \stackrel{X_{\eta} \text{-} X_{\vlon}}{\cong}
  \t{Range} \, \varphi_{\eta l} (q)$, we have $\psi_{\vlon , \eta}^{k} (p) =
  \psi_{\vlon , \eta}^{l} (q)$.
\end{cor}
\begin{defn}
  Weight of a projection $p \in P_{\vlon k}$ for even (resp., odd)
  $k$, denoted by $wt(p)$, is defined to be the smallest even (resp.,
  odd) integer $l$ such that there exists a projection $q \in P_{\vlon
    l}$ satisfying $\t{Range} \, (\varphi_{\vlon k}(p)) \cong \t{Range} \,
  (\varphi_{\vlon l}(q))$ as $X_{\vlon}$-$X_{(-)^k \vlon}$-bimodules.
\end{defn}
Let ${\mcal S}_{\vlon k}$ be a maximal set of non-equivalent minimal
projections in $P_{\vlon k}$ with weight $k$ for all colors $\vlon k$.
\begin{rem}\label{spanset}
  In view of Remark \ref{psi-remark}, Lemma \ref{psi-properties} and
  Corollary \ref{psi-corollary}, we observe that $AP_{\vlon 0, \eta
    0}$ is spanned by the set $\us{k \in (2\N_0 + \delta_{\vlon \neq
      \eta})}{\cup} \left\{ \psi^k_{\vlon, \eta} (p) : p \in {\mcal
      S}_{\eta k} \right\}$ for $\vlon , \eta \in \{+, -\}$.
\end{rem}
We shall, in fact, see that these sets are linearly independent and
hence form bases.

\subsection{Equivalence on tangles induced by affine isotopy}\label{eqtang} $ $

\begin{figure}[h!]
\psfrag{k}{$k$}
\psfrag{e}{$\eta$}
\psfrag{X}{$X$}
\psfrag{Y}{$Y$}
\psfrag{s}{$\sim$}
\psfrag{k-i-2}{$k\!-\!i\!-\!2$}
\psfrag{i}{$i$}
\psfrag{(i)}{$(i)$}
\psfrag{(ii)}{$(ii)$}
\includegraphics[scale=0.2]{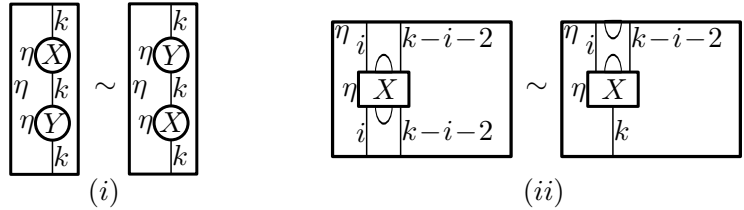}
\caption{Equivalence relation $\sim\,$. ($\eta \in \{+, - \}$, $k \in
  \N_0$, $0 \leq i \leq k-2$, $X , Y \in \mcal{T}_{\eta
    k}(P)$)}\label{eq-relation}
\end{figure}

For $\vlon, \eta \in \{+,- \}$, set ${\mcal T}_{\vlon , \eta} := \us{l
  \in \N_0}{\sqcup} {\mcal T}_{\eta (2l + \delta_{\vlon \neq \eta})}
(P)$.  Define the equivalence relation $\sim$ on $\mcal{T}_{\vlon ,
  \eta}$ generated by the relations given by the pictures in Figure
\ref{eq-relation}.
The following topological lemma involving this equivalent relation will be crucial in the forthcoming section.
\begin{lem}\label{psi-equivalence}
If $\vlon, \eta \in \{+,-\}$, $k , l \in (2\N_0 + \delta_{\vlon \neq \eta})$ and $S \in {\mcal T}_{\eta k} (P) , T \in {\mcal T}_{\eta l} (P)$, then $\Psi_{\vlon 0, \eta 0}^{k} (S) = \Psi_{\vlon 0, \eta 0}^{l} (T)$ if and only if $X \sim Y$.
\end{lem}
\begin{pf}
If $S \sim T$ either by relation (i) or (ii) in Figure \ref{eq-relation}, then using affine isotopy, we easily see that  $\Psi_{\vlon 0, \eta 0}^{k} (U) = \Psi_{\vlon 0, \eta 0}^{l} (V)$.
For the `only if' part, consider pictures $\hat{S}$ and $\hat{T}$ in
the isotopy class of $S$ and $T$ respectively, and
$\hat{\Psi}_{\vlon 0, \eta 0}^{m}$ as in Figure \ref{aff-tangles} to
represent $\Psi_{\vlon 0, \eta 0}^m$ for $m = k, l$.  Since
$\Psi_{\vlon 0, \eta 0}^{k} (S) = \Psi_{\vlon 0, \eta 0}^{l} (T)$,
we have an affine isotopy $\varphi : [0,1] \times\, RA \ra RA$ (as
in Definition \ref{aff-isotopy}) such that $\varphi \left( 1 ,
\hat{\Psi}_{\vlon 0, \eta 0}^{k} (\hat{X}) \right) =
\hat{\Psi}_{\vlon 0, \eta 0}^{l} (\hat{Y})$.  Let $p$ be the
straight path in $RA$ joining the points $(0, -1)$ and $(0,-2)$ and
suppose $\tilde{p} := \varphi\left(1 , p \right)$ which is also a
simple path in $RA$ joining the same two points.  Note that cutting
$\hat{\Psi}_{\vlon 0, \eta 0}^{k} (\hat{X}) $ (resp.,
$\hat{\Psi}_{\vlon 0, \eta 0}^{l} (\hat{Y}) $) along the path $p$
and straightening gives $\hat{X}$ (resp., $\hat{Y}$), as shown in
Figure \ref{cut}.
\begin{figure}[h!]
\psfrag{p}{$p$}
\psfrag{k}{$k$}
\psfrag{=}{$\!_{=}$}
\psfrag{e}{${\vlon}$}
\psfrag{y}{$\eta$}
\psfrag{X}{$\hat{X}$}
\psfrag{cut}{$\ous{\text{cutting}}{\longmapsto}{\text{along } p}$}
\includegraphics[scale=0.25]{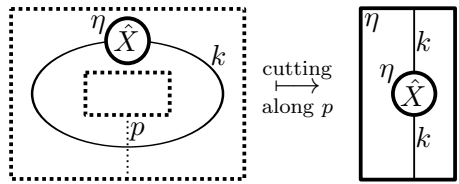}
\caption{Cutting along a simple path}\label{cut}
\end{figure}
Further, since $\varphi$ is an affine isotopy, even if we cut $\hat{\Psi}_{\vlon 0, \eta 0}^{l} (\hat{Y}) $ along $\tilde{p}$, we still obtain $\hat{X}$ (upto planar isotopy).
Let $A_0$ denote the affine tangular picture $\hat{\Psi}_{\vlon 0, \eta 0}^{l} (\hat{Y})$ and $\mcal{SP} (A_0)$ denote the set of those simple paths in $RA$ with end points $(0, -1)$ and $(0, -2)$ such that they (a) do not meet any disc in $A_0$, (b) intersect the set of strings discretely and non-tangentially, and (c) are equivalent to the straight path $p$ via some affine isotopy.
Clearly, $p, \tilde{p} \in \mcal{SP}(A_0)$.

Analogous to the equivalence relation $\sim$ on $\mcal{T}_{\vlon ,
  \eta}$, we consider a equivalence relation $\sim$ on $\mcal{SP}
(A_0)$ generated by the local moves as shown in Figure
\ref{local-moves}.
\begin{figure}[h]
\psfrag{s}{$\sim$}
\psfrag{i}{$(i)$}
\psfrag{k}{$k$}
\psfrag{l}{$l$}
\psfrag{ii}{$(ii)$}
\psfrag{iii}{$(iii)$}
\psfrag{iii'}{$(iii)'$}
\psfrag{1}{$(i)$}
\psfrag{p}{$p_1$}
\psfrag{p2}{$p_2$}
\psfrag{p3}{$_{p_3}$}
\psfrag{x}{$x$}
\includegraphics[scale=0.2]{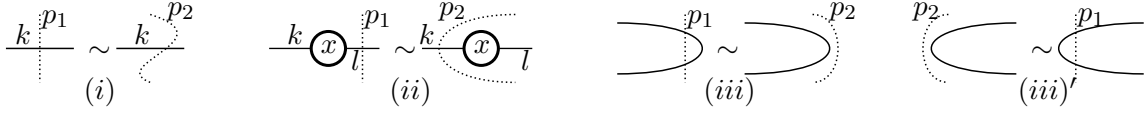}
\caption{Equivalence relation on $\mcal{SP}(A_0)$. ($k , l \in \N_0$ such that $(k+l) \in 2 \N_0$, $x \in P_{\pm (\frac{k+l}{2})}$)}\label{local-moves}
\end{figure}
Note that cuts along two paths related by move $(i)$ give same
labelled tangles (upto tangle isotopy); and cuts along paths related
by moves $(ii)$ and $\{(iii), (iii)'\}$ correspond to equivalence
relations $(i)$ and $(ii)$ of Figure \ref{eq-relation}, respectively.
Thus, it is enough to show that the paths $p$ and $\tilde{p}$ are
equivalent under this relation which will imply $X \sim Y$.  It is not
hard to prove that $p$ can obtained from $\tilde{p}$ by applying
finitely many moves of the above types.  We will not give a complete
proof of this fact here; however, one can extract a detailed proof
from the strategy used in the proof of \cite[Proposition $2.8$]{Gho06}
which proves the same type of statement but for `annular tangles'
where the isotopy has no restriction on the internal and external
boundaries as in affine isotopy.  So, one has to make necessary
modifications, namely, ignoring the rotation move in \cite{Gho06} but
even this is not an issue for us because we are working with affine
morphism from $\vlon 0$ to $\eta 0$ and no strings are attached to the
boundary of $RA$.  This completes the proof of the lemma.
\end{pf}
\subsection{Proof of Theorem \ref{main-theorem}}
We first set up the following notation:\\
\noindent For $p \in \mscr{P}_{min} (\mcal{Z}(P_{\vlon l}))$ and $\eta
= (-)^l \vlon$, we write $v^p_{\eta , \vlon} \in V_{\eta, \vlon}$ for
the isomorphism class of the $X_{\vlon}$-$X_{\eta}$ bimodule $\t{Range} \,
\varphi_{\vlon l}(p_0)$ for any $p_0 \in \mscr{P}_{min}(P_{\vlon l})$
with $p_0 \leq p$.\\
\noindent Now, for $\vlon, \eta \in \{+, -\}$ and $k \in (2\N_0 +
\delta_{\vlon \neq \eta})$, consider the map
\[
P_{\eta k} \ni x \os{\gamma^k_{\vlon , \eta}}{\longmapsto} \us{{p \in
    \mscr{P}_{min} (\mcal{Z} (P_{\eta k}))}}{\sum} {\sqrt{dim\, (p
    P_{\eta k})}} \; \left[ \dfrac{tr_{M_{(k-\delta_{\eta = +})}}
    (xp)}{tr_{M_{(k-\delta_{\eta = +})}} (p)} \right] \; v^p_{\vlon ,
  \eta} \in \mbb{C} V_{\vlon , \eta}.
\]
\begin{rem}\label{gamma-range}
  The above definition directly implies $\gamma^k_{\vlon, \eta} (p_0)
  = v^p_{\vlon , \eta}$ for all $p \in \mscr{P}_{min} (
  \mcal{Z}(P_{\eta k}))$ and $p_0 \in \mscr{P}_{min} (P_{\eta k})$
  satisfying $p_0 \leq p$.
\end{rem}
\begin{lem}\label{gamma-properties} 
  If $\vlon, \eta \in \{+, -\}$ and $k \in (2\N_0 + \delta_{\vlon \neq
    \eta})$, then
\begin{enumerate} 
\item $\gamma^k_{\vlon, \eta}$ is tracial,
\item $\gamma^k_{\vlon, \eta} (x) = \gamma^{k+2}_{\vlon, \eta}(x
  e_{(k+1+\delta_{\eta = -})})$ for all $x \in P_{\eta k}$.
\end{enumerate}
\end{lem}
\begin{pf}
  Note that any partial isometry in $P_{\eta k}$ with orthogonal
  initial and final projections, is in the kernel of $\gamma^k_{\vlon
    , \eta}$; this along with Remark \ref{gamma-range} imply (1).

  For (2), let $\left\{ e^p_{i, j} : p \in
    \mscr{P}_{min}(\mcal{Z}(P_{\eta k})), 1 \leq i, j \leq \sqrt{dim(
      pP_{\eta k})} \right\}$ be a system of matrix units for $P_{\eta
    k}$.  Fix a $p \in \mscr{P}_{min}(\mcal{Z}(P_{\eta k}))$.  Then,
  by (1), $\gamma^{k}_{\vlon, \eta}(e^p_{i, j}) = 0 =
  \gamma^{k+2}_{\vlon, \eta}(e^p_{i, j}\, e_{k+1})$ for all $1 \leq i
  \neq j \leq \sqrt{dim (pP_{\eta k})}$.  It is easy to check that
  $e^p_{i,i} e_{(k+1 + \delta_{\eta = -})}$ is a minimal projection;
  let $\tilde{p}$ be its central support in $P_{\eta (k+2)}$.  By
  Remark \ref{gamma-range} and Lemma \ref{e-even-odd}, we have
  $\gamma^{k}_{\vlon, \eta} (e^p_{i,i}) = v^p_{\vlon , \eta} =
  v^{\tilde{p}}_{\vlon , \eta} = \gamma^{k+2}_{\vlon, \eta} (e^p_{i,i}
  e_{(k+1 + \delta_{\eta = -})})$.
\end{pf}
\begin{cor}\label{gamma-equivalence}
  If $\vlon, \eta \in \{+, -\}$, $k, l \in (2\N_0 + \delta_{\vlon \neq
    \eta})$, $S \in \mcal{T}_{\eta k}(P)$ and $T \in \mcal{T}_{\eta
    l}(P)$ such that $S \sim T$, then $\gamma^k_{\vlon, \eta} (P_S) =
  \gamma^l_{\vlon, \eta} (P_T)$.
\end{cor}
\begin{pf}
If $S \sim T$ by relation $(i)$, as shown in Figure \ref{eq-relation}, then part (1) of Lemma \ref{gamma-properties} does the job.
Suppose $S$ and $T$ denote the tangles on the left and the right sides of relation $(ii)$ in Figure \ref{eq-relation} respectively, and let $Z=
\psfrag{y}{$\eta$}
\psfrag{i}{$i$}
\psfrag{k-i-2}{$k-i-2$}
\includegraphics[scale=0.2]{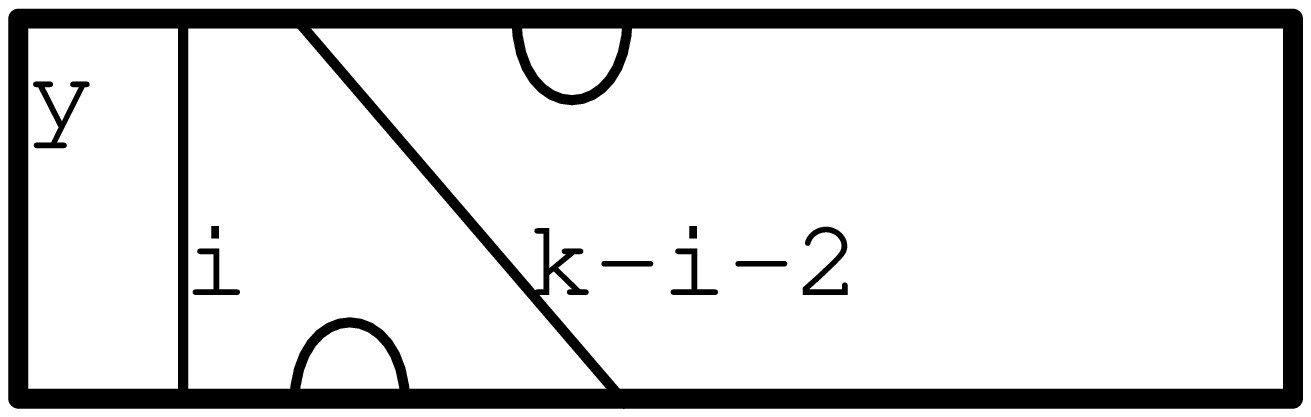}$.
Then, we have 
\[
\gamma^{k-2}_{\vlon, \eta} (P_{S})
= \gamma^{k}_{\vlon, \eta} (P_{S} \; e_{(k - 1 + \delta_{\eta = -} )})
= \delta^{-1} \gamma^{k}_{\vlon, \eta} (P_Z P_X P_{Z^*})
= \delta^{-1} \gamma^{k}_{\vlon, \eta} (P_{Z^*} P_Z P_X)
= \gamma^{k}_{\vlon, \eta} (P_T)
\]
where we use parts (2) and (1) of Lemma \ref{gamma-properties} to
obtain the first and third equalities.
\end{pf}

We are now just one step away from establishing the required isomorphism.
For $\vlon , \eta \in \{ +, - \}$, consider the map
\[
\mcal{AT}_{\vlon 0, \eta 0}(P) \ni A \os{ \Lambda_{ \vlon , \eta } }{\longmapsto} \gamma_{\vlon, \eta }^k (P_T) \in \mbb{C}V_{\vlon , \eta}
\]
where (by Remark \ref{Psi-remark}) $A = \Psi_{\vlon 0, \eta0}^{k} (T)$ for some $k \in (2 \N_0 + \delta_{\vlon \neq \eta})$ and $T \in
\mcal{T}_{\eta k}(P)$. 
$\Lambda_{\vlon, \eta}$ is indeed a well-defined map due to Corollary \ref{gamma-equivalence}.
Extend this map linearly to $\Lambda_{\vlon, \eta} : \mcal A_{\vlon 0 , \eta 0} (P) \ra \C V_{\vlon , \eta}$.
Note that $\Lambda_{\vlon, \eta} (A) = \gamma_{\vlon, \eta }^k (P_X)$ whenever $A = \Psi_{\vlon 0, \eta0}^{k} (X)$ and $X \in \mcal{P}_{\eta k}(P)$; this implies
$ \mcal{W}_{\vlon 0, \eta 0} \subset \t{ker}\,
\Lambda_{\vlon, \eta} $.
Thus, each $\Lambda_{\vlon, \eta}$ induces a linear map $\lambda_{\vlon, \eta} : AP_{\vlon 0, \eta 0} \lra \mbb{C} V_{\vlon , \eta}$, that is, $\Lambda_{\vlon , \eta} = \lambda_{\vlon , \eta} \circ q_{\vlon , \eta}$.

\noindent {\em Proof of Theorem \ref{main-theorem}:}
Define $\lambda :=
\begin{bmatrix}
\lambda_{+,+} & \lambda_{-,+}\\
\lambda_{+,-} & \lambda_{-,-}
\end{bmatrix}$.  We will show that $\lambda : AP_{0,0} \lra
\mcal{F}_{N \subset M}$ is a $\ast$-algebra isomorphism.  Clearly,
$\lambda$ is linear. Now, for $\vlon , \eta \in \{+, - \}$, $k \in (2
\N_0 + \delta_{\vlon \neq \eta})$ and $p \in {\mcal S}_{\eta k}$
(defined before Remark \ref{spanset}), let $\tilde{p}$ denote the
central support of $p$ in $P_{\eta k}$.  Note that $\lambda_{\vlon ,
  \eta} (\psi_{\vlon , \eta} (p)) = \Lambda_{\vlon , \eta}
(\Psi_{\vlon , \eta} (I_{\eta k} (p))) = \gamma_{\vlon , \eta} (p) =
v^{\widetilde{p}}_{\vlon , \eta} \in V_{\vlon , \eta}$ where the first
two equalities follow easily unravelling the definitions and the last
one comes from Remark \ref{gamma-range}.  On the other hand, from
Corollary \ref{p-q-equivalence} and definition of $V_{\vlon , \eta}$,
we get $\left\{v^{\tilde{p}}_{\vlon , \eta} : k \in (2\N_0 +
  \delta_{\vlon \neq \eta}), p \in {\mcal S}_{\eta k} \right\} =
V_{\vlon , \eta}$. This and Remark \ref{spanset}, imply that $\lambda_{\vlon , \eta}$ is injective as well as surjective.

A closer look at the $\ast$-structures of $\mcal{F}_{N \subset M}$
(resp., $AP_{0,0}$) reveals $\left[ v^{ \widetilde{p} }_{\vlon ,
    \eta} \right]^* = v^{ \widetilde{q}}_{\eta , \vlon}$ using
Proposition \ref{rotcont} (resp., $\left[ \psi_{\vlon , \eta} (p)
\right]^* = \psi_{\eta , \vlon} (q)$) where $q = P_{R^k_{\eta k}} (p)$
for all $\vlon , \eta \in \{+ , -\}$, $k \in (2 \N_0 + \delta_{\vlon
  \neq \eta})$ and $p \in {\mcal S}_{\eta k}$. Hence, $\lambda$ is
$*$-preserving.

It remains to show that $\lambda$ is an algebra homomorphism.  Note
that for $\vlon , \eta , \nu \in \{+ , -\}$, $k \in (2\N_0 +
\delta_{\nu \neq \eta})$, $l \in (2\N_0 + \delta_{\eta \neq \vlon})$,
$x \in P_{\nu k}$ and $y \in P_{\eta l}$, we have $\psi_{\vlon , \nu}
\left( P_{H_{\nu k , \eta l}} (x , y) \right) = \psi_{\eta , \nu} (x)
\circ \psi_{\vlon , \eta} (y)$ where the tangle $H_{\nu k , \eta l}$
is given by \psfrag{y}{$\eta$} \psfrag{k}{$k$} \psfrag{l}{$l$}
\psfrag{n}{$\nu$}
\includegraphics[scale=0.2]{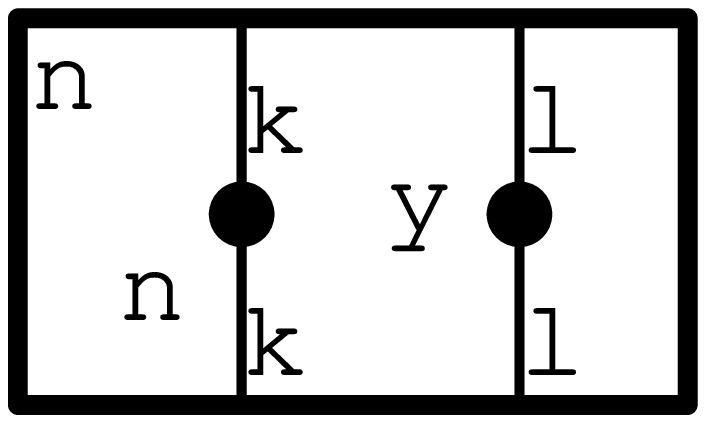}.
So, one needs to
check $\t{Range} \, \vphi_{\nu (k+l)} (P_{H_{\nu k , \eta l}} (p , q)) \cong
\t{Range} \, \vphi_{\nu k} (p) \us{X_\eta}{\otimes} \t{Range} \, \vphi_{\eta l} (q)$
as $X_\nu$-$X_\vlon$-bimodules where $p \in {\mathscr P} (P_{\nu k})$
and $q \in {\mathscr P} (P_{\eta l})$.  One way of seeing this is by
translating some results in \cite[Theorem $4.6$]{Bis97} in the
language of planar algebras.  However, this isomorphism comes for free
from the isomorphism between $P$ and the normalized bimodule planar
algebra associated to ${_N} L^2 (M)_M$, established in the proof of
\cite[Theorem $5.4$]{DGG}.

Hence, $\lambda$ is a $*$-algebra isomorphism.\qed
\section{Affine modules with zero weight}\label{regmod}
In this section, we will analyze the affine $P$-modules with weight zero for any subfactor planar algebra $P$ (possibly having infinite depth).

Throughout this section, $\vlon$ will denote an element of $\{\pm\}$ and $P$ will continue to be the planar algebra associated to a finite index subfactor $N \subset M$.
\comments{
Set $I^\vlon := AP_{- \vlon 0 , \vlon 0} \circ AP_{\vlon 0 , -\vlon 0}$.
Note that $I^\vlon$ is a $*$-closed two-sided ideal in $AP_{\vlon 0,\vlon 0}$.
Given any $*$-affine $P$-module $V$ with weight zero, one can consider the submodule $V^0$ generated by $I^+ (V_{+ 0})$ which is the same as the one generated by $I^- (V_{- 0})$.
Let $V^\vlon$ be the submodule of $W:= \left( V^0 \right)^\perp$ 

\begin{rem}
Given any $*$-affine $P$-module $V$, one can define submodules
$V^+$,
$V^0 := $, and $V^-$ satisfying:

(i) $V^\vlon_{-\vlon 0} = \{0\}$

(ii)
\end{rem}
}
Let us consider the trace  on the algebra $\C V_{\vlon,\vlon}$ (introduced in Section \ref{mt}) given by $V_{\vlon,\vlon} \ni v \os{\omega_\vlon}{\ra} \delta_{v = 1_\vlon} \in \C$ where $1_\vlon$ is the isomorphism class of the trivial bimodule in $V_{\vlon , \vlon}$.
Clearly, $\omega_\vlon$ is positive definite.
By the isomorphism in Theorem \ref{main-theorem}, $\omega_\vlon$  induces a positive definite trace on $AP_{\vlon 0 , \vlon 0}$.
In the following lemma, we present a pictorial interpretation of $\omega_\vlon$.
\begin{lem}\label{oe}
For all $k \in \N_0$ and $x \in P_{\vlon 2k}$, we have
\[
\omega_\vlon \left( \psi^{2k}_{\vlon 0 , \vlon 0} (x) \right) = \delta^{-k} \us{\alpha}{\sum} P_{
\psfrag{w}{$w_\alpha$}
\psfrag{w*}{$w^*_\alpha$}
\psfrag{2k}{$2k$}
\psfrag{x}{$x$}
\includegraphics[scale=0.25]{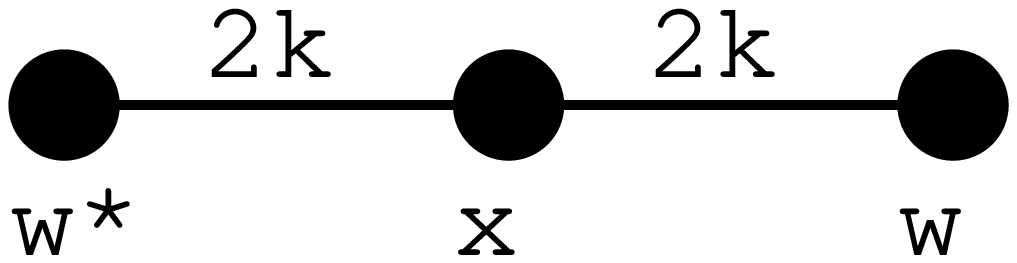}
} \in P_{\vlon 0} \cong \C
\]
where $\left\{ w_\alpha \right\}$ is any orthonormal basis of $P_{\vlon k}$ with repect to the canonical trace (that is, the normalized picture trace).
\end{lem}
\begin{pf}
Let $\left\{ E^i_{\alpha , \beta} : 0 \leq i \leq n , 1 \leq \alpha , \beta \leq d_i \right\}$ be a system of matrix units of the finite dimensional $C^*$-algebra $P_{\vlon 2k}$ where $i$ gives the indexing of the matrix summands and $d_i$ is the order of $i$-th summand; further, let us assume the $0$-th summand is the one whose minimal projections correspond to $1_\vlon \in V_{\vlon , \vlon}$.
Now, there exist scalars $x^i_{\alpha , \beta}$ such that $x = \us{i}{\sum} \us{\alpha, \beta}{\sum} x^i_{\alpha , \beta} E^i_{\alpha , \beta}$.
So, by the isomorphism in Theorem \ref{main-theorem} and definition of $\omega_\vlon$, we get $\omega_\vlon \left( \psi^{2k}_{\vlon 0 , \vlon 0} (x) \right) = \us{\alpha}{\sum} x^0_{\alpha , \alpha}$.

Consider the minimal projection $p := \delta^{-k} P \,_{
\psfrag{k}{$k$}
\psfrag{e}{$\vlon$}
\includegraphics[scale=0.25]{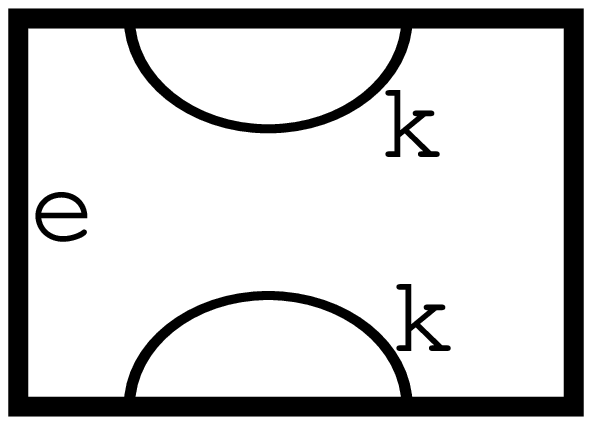}
}$ in $P_{\vlon 2k}$, which also corresponds to $1_\vlon \in V_{\vlon , \vlon}$.
Let $v_\alpha \in P_{\vlon 2k}$ such that $v_\alpha v^*_\alpha = E^0_{\alpha , \alpha}$ and $v^*_\alpha v_\alpha = p$.
Set $w_\alpha := P \,_{
\psfrag{va}{$v_\alpha$}
\psfrag{2k}{$2k$}
\psfrag{e}{$\vlon$}
\includegraphics[scale=0.25]{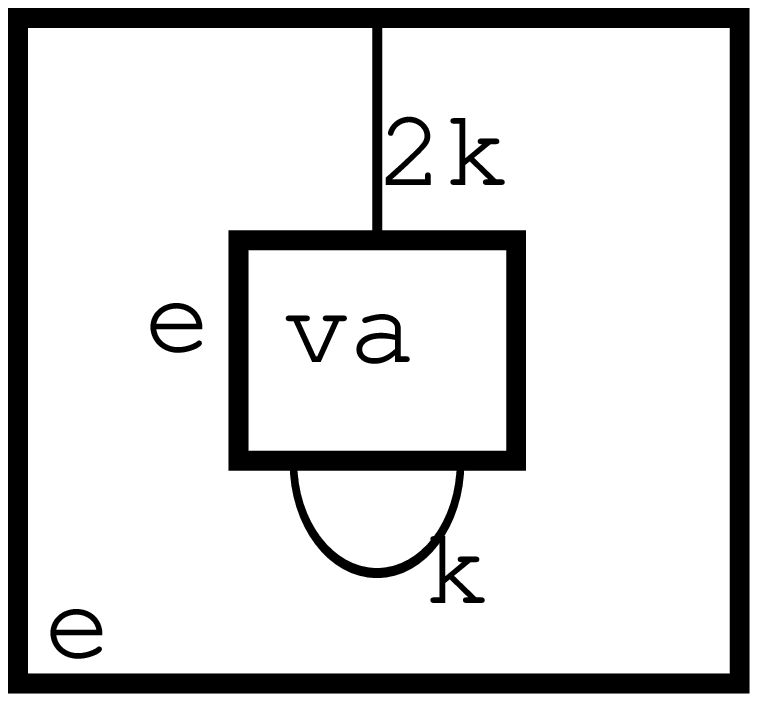}
} \in P_{\vlon k}$.
It easily follows from $v^*_\alpha v_\beta = \delta_{\alpha = \beta} p$ that $\{ w_\alpha \}_\alpha$ is an orthonormal subset of $P_{\vlon k}$ with respect to the canonical trace. On the other hand, $v^*_\alpha x v_\alpha = x^0_{\alpha , \alpha} p$ implies $P_{
\psfrag{w}{$w_\alpha$}
\psfrag{w*}{$w^*_\alpha$}
\psfrag{2k}{$2k$}
\psfrag{x}{$x$}
\includegraphics[scale=0.25]{figures/modwt0/oe.eps}
} = \delta^k x^0_{\alpha , \alpha} 1_{P_{\vlon 0}}$.
It only remains to show that $\{w_\alpha\}_{1 \leq \alpha \leq d_0}$ spans $P_{\vlon k}$. For this, we use Frobenius reciprocity for bimodules and get $d_0 = \text{dim} (P_{\vlon k})$.

Independence from the choice of an orthonormal basis of $P_{\vlon k}$, follows from the equation $\omega_\vlon \left( \psi^{2k}_{\vlon 0 , \vlon 0} (x) \right) = \us{\alpha}{\sum} \langle w_\alpha , f_x (w_\alpha ) \rangle_{P_{\vlon k}}$ where $f_x : P_{\vlon k} \ra P_{\vlon k}$ is the linear operator given by the action of the semi-labelled tangle
\psfrag{x}{$x$}
\psfrag{.}{}
\psfrag{2k}{$2k$}
\psfrag{m}{$m$}
\psfrag{e}{$\vlon$}
\includegraphics[scale=0.25]{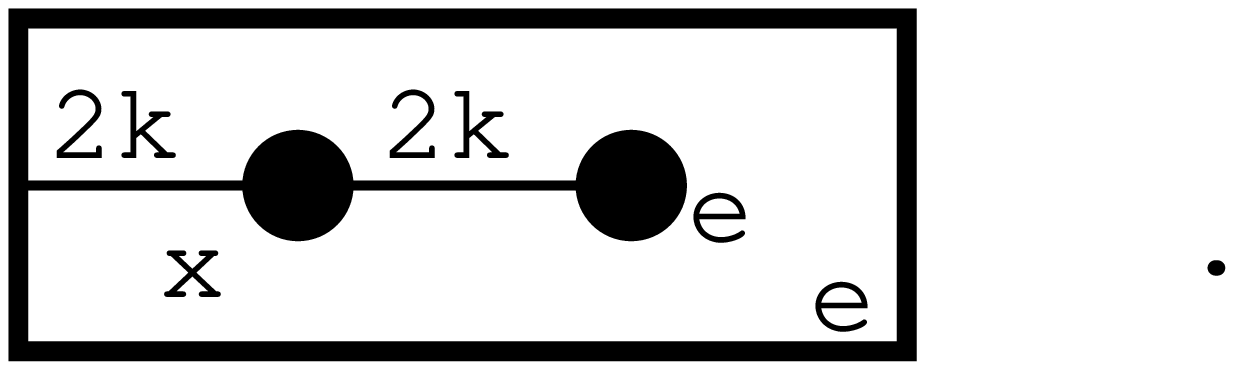}.
\end{pf}

We now define $H^{\vlon}_{\eta k} := AP_{\vlon 0 , \eta k}$ for all $\eta k \in \text{Col}$.
$H^\vlon = \left\{ H^\vlon_{\eta k} \right\}_{\eta k \in \t{Col}}$ forms an affine $P$-module with action of affine morphisms given by composition.
Define a sesquilinear form on the affine module $H^\vlon$ in the following way:
\[
\langle h_1 , h_2 \rangle := \omega_\vlon (h^*_1 \circ h_2) \text{ where } h_1 , h_2 \in H^{\vlon}_{\eta k} \t{ and } \eta k \in \text{Col}.
\]
\begin{thm}\label{hethm}
$H^\vlon$ is a bounded $*$-affine $P$-module with inner product given by the above form. Hence, its completion will be a Hilbert affine $P$-module.
\end{thm}
\begin{pf}
We first need to check whether the form is positive definite, that is, $\omega_\vlon (h^* \circ h) > 0$ for $0 \neq h \in H^\vlon_{\eta k} = AP_{\vlon 0 , \eta k}$. 
For each $m\in(2\N_0 + \delta_{\vlon \neq \eta})$, set $\Phi^m_{\vlon , \eta k} :=
\psfrag{AR}{$AR_{\vlon k}$}
\psfrag{y}{$\eta$}
\psfrag{2k-1}{$2k-1$}
\psfrag{2k}{$2k$}
\psfrag{2l}{$2k$}
\psfrag{m}{$m$}
\psfrag{e}{$\varepsilon$}
\psfrag{-e}{$-\varepsilon$}
\psfrag{psi}{$\Psi^m_{\varepsilon k, \eta l}$}
\psfrag{1ek}{$A1_{\varepsilon k}$}
\includegraphics[scale=0.30]{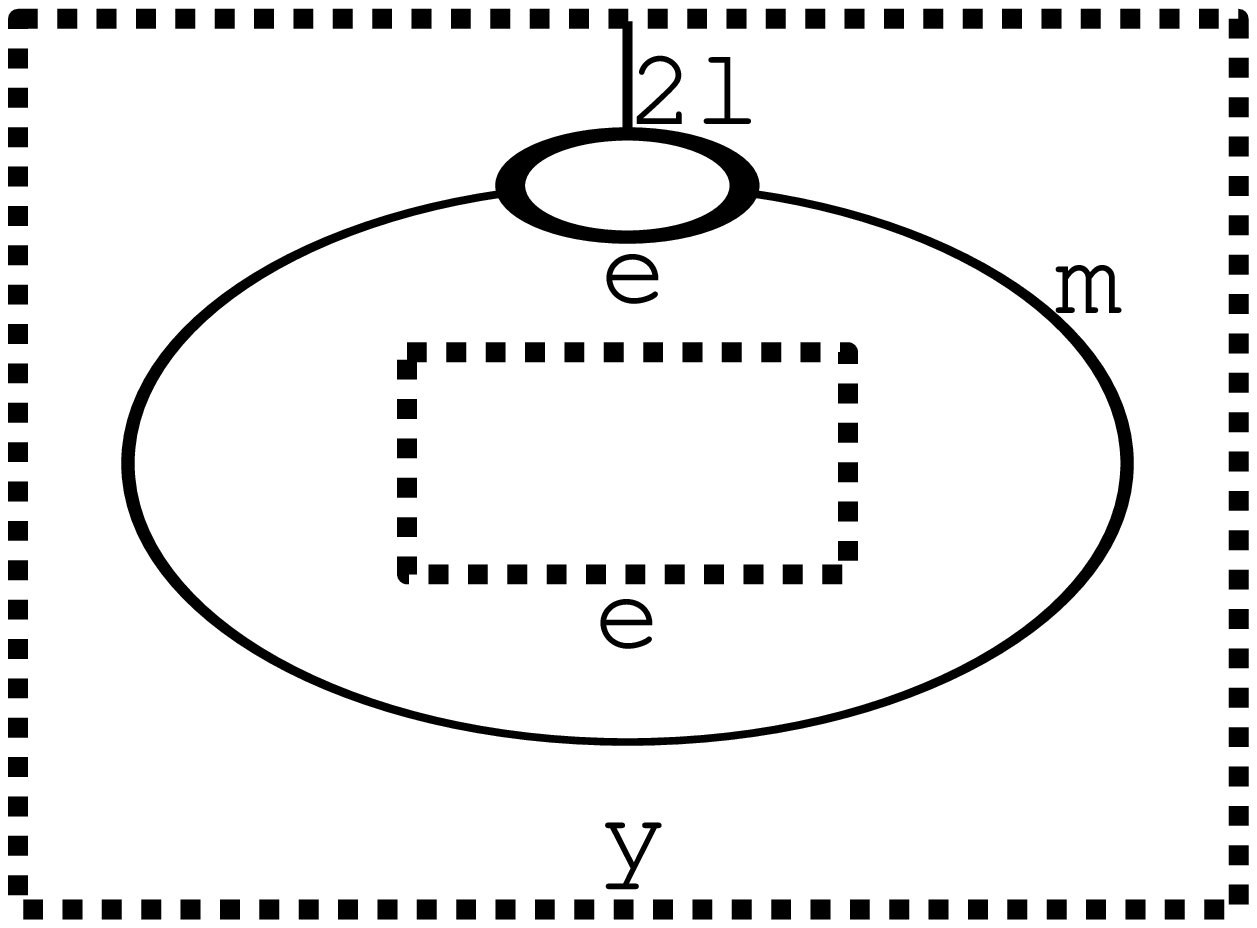}
$ which is the same as the unlabelled affine tangle $\Psi^m_{\vlon 0, \eta k}$ (defined in Figure \ref{aff-tangles}) except there is a certain rotation on the internal disc.
Let $\vphi^m_{\vlon , \eta k} : P_{\vlon (m+k)} \ra AP_{\vlon 0 , \eta k}$ be the linear map induced by the affine tangle $\Phi^m_{\vlon , \eta k}$.
Note that using affine isotopy, we can obtain 
\begin{equation}\label{commarnd}
\vphi^m_{\vlon , \eta k} \left(P_{RI^k_{\vlon m}} (y) x \right) = \vphi^m_{\vlon , \eta k} \left(x P_{RI^k_{\vlon m}} (y) \right) \t{ for all } x \in P_{\vlon (m+k)}
\end{equation}
and $y \in P_{\vlon m}$ where $RI^k_{\vlon m} : \vlon m \ra \vlon (k+m)$ is the tangle obtained from $RI_{\vlon m}$ (in Figure \ref{tangles}) by replacing the straight vertical string on the right by $k$ parallel strings.
Considering a path algebra model of $\C \cong P_{\vlon 0} \hookrightarrow P_{\vlon m} \os{P_{RI^k_{\vlon m}}}{\hookrightarrow}  P_{\vlon (m+k)}$ and using Equation \ref{commarnd}, we may deduce $\t{Range} \, \vphi^m_{\vlon , \eta k} = \vphi^m_{\vlon , \eta k} \left( P^\prime_{\vlon m} \cap P_{\vlon (m+k)} \right)$.
This along with Remark \ref{psi-remark} implies that there exist $m \in \N_0$ and $0 \neq x \in P^\prime_{\vlon m} \cap P_{\vlon (m+k)}$ such that $h = \vphi^m_{\vlon , \eta k } (x)$.
By Lemma \ref{oe}, we get
\[
\omega_\vlon (h^* \circ h)
= \delta^{-m} \us{\alpha}{\sum} P_{\!\!\!\!\!\!\!\!
\psfrag{w*}{$w_\alpha$}
\psfrag{w}{$w^*_\alpha$}
\psfrag{2k}{$2k$}
\psfrag{m}{$m$}
\psfrag{x}{$x$}
\psfrag{x*}{$x^*$}
\includegraphics[scale=0.25]{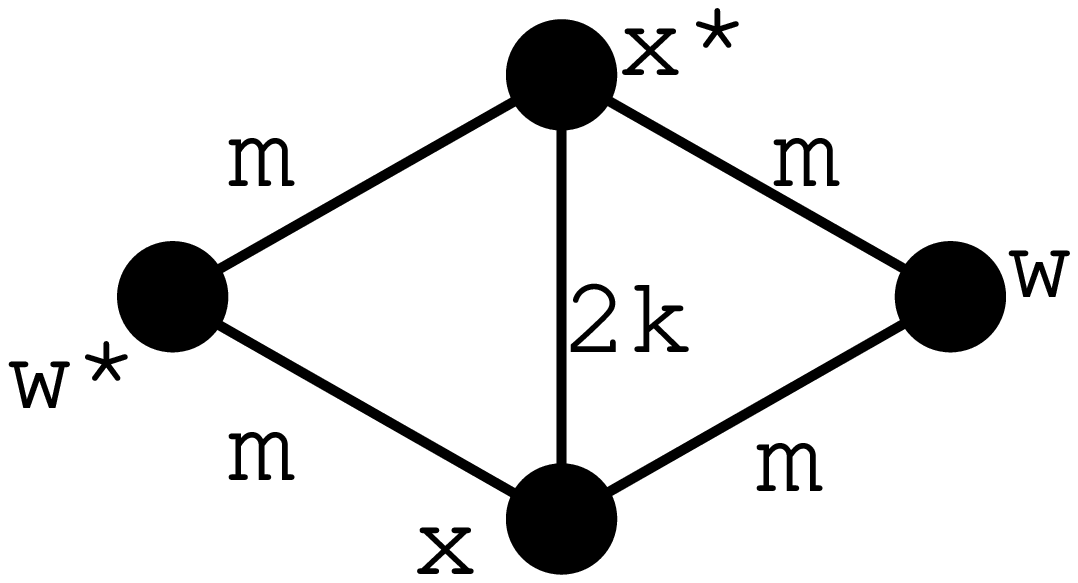}
}
= \delta^{-m} \us{\alpha}{\sum} P_{\!\!\!\!\!\!\!\!\!\!\!\!\!\!\!\!\!\!
\psfrag{w}{$w_\alpha w^*_\alpha$}
\psfrag{2k}{$2k+m$}
\psfrag{m}{$m$}
\psfrag{x}{$x$}
\psfrag{x*}{$x^*$}
\includegraphics[scale=0.25]{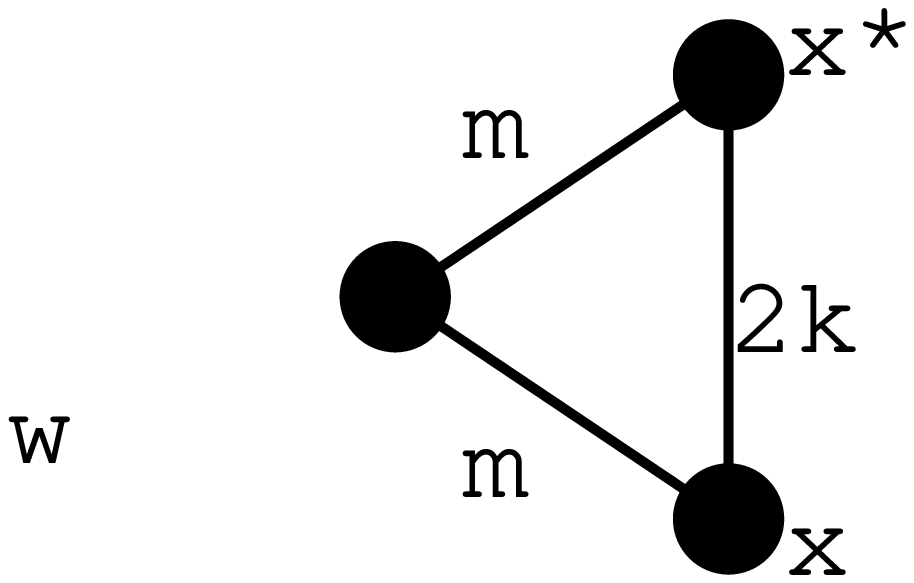}
}
\;\;\;= \delta^{-m} \us{\alpha}{\sum} P_{TR^r_{\vlon m}} (w_\alpha w^*_\alpha y)
= \us{\alpha}{\sum} \left\langle w_\alpha , y w_\alpha \right\rangle_{P_{\vlon m}}
\]
where $y = P_{
\psfrag{x}{$x$}
\psfrag{x*}{$x^*$}
\psfrag{2k}{$2k+m$}
\psfrag{m}{$m$}
\psfrag{e}{$\vlon$}
\includegraphics[scale=0.25]{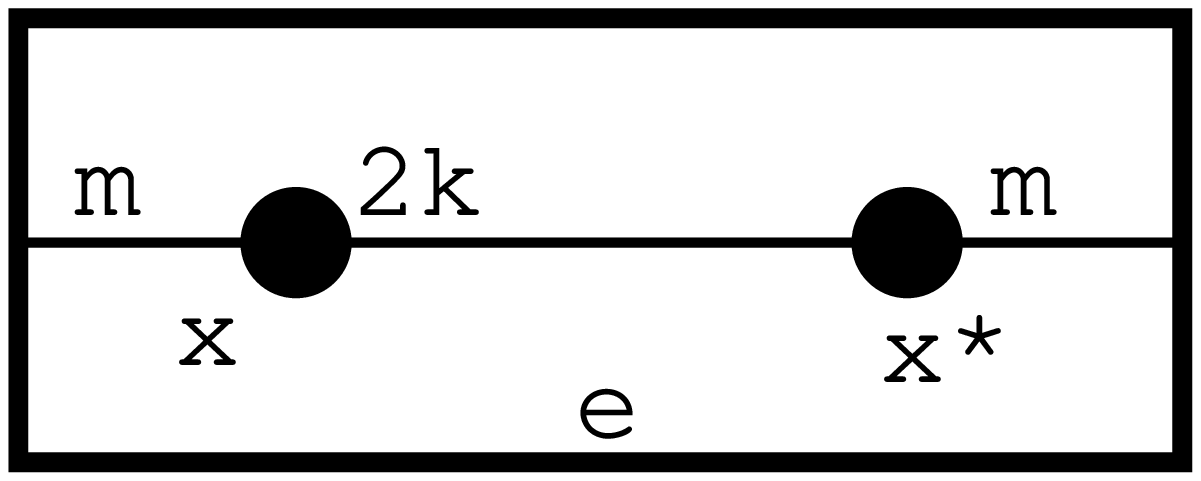}
}$ and $\{w_\alpha\}_\alpha$ is an orthonormal basis of $P_{\vlon m}$ with respect to the canonical trace.
The second equality follows from $[x,P_{\vlon m}] = 0$.
Note that $y$ is a positive element of $P_{\vlon m}$ and nonzero too since $P_{TR^r_{\vlon m}} (y) = P_{TR^r_{\vlon (m+k)}} (x x^*) \neq 0$. Thus, $\omega_\vlon (h^* \circ h) = \us{\alpha}{\sum} \left\| y^{1/2} w_\alpha \right\|^2 > 0$.

The $*$-preserving condition $\langle a \circ h_1 , h_2 \rangle = \langle h_1 , a^* \circ h_2 \rangle$ holds trivially. Hence, $H^\vlon$ is a $*$-affine module.

\noindent{\em Boundedness of the action of affine morphisms:}
This part is relevant only if depth of $P$ is infinite since for finite depth planar algebras, $H^\vlon$ will be locally finite (see \cite[Proof of Theorem 6.11]{Gho}).
Let $a = \psi^m_{\eta k , \nu l} (x) \in AP_{\eta k , \nu l}$ and $h = \vphi^n_{\vlon , \eta k} (y) \in H^\vlon_{\eta k} = AP_{\vlon 0 , \eta k}$ where $x \in P_{\nu (k+l+m)}$ and $y \in P_{\vlon (k+n)}$.
Now, $\| a \circ h \|^2 = \omega_\vlon (h^* \circ a^* \circ a \circ h)$ which, using Lemma \ref{oe}, can be expressed as
\[
\delta^{-(m+n)} \us{\alpha}{\sum} P_{
\psfrag{w}{$w_\alpha$}
\psfrag{w*}{$w^*_\alpha$}
\psfrag{x}{$x$}
\psfrag{x*}{$x^*$}
\psfrag{y}{$y$}
\psfrag{y*}{$y^*$}
\psfrag{2k}{$2k$}
\psfrag{2l}{$2l$}
\psfrag{m}{$m$}
\psfrag{n}{$n$}
\includegraphics[scale=0.20]{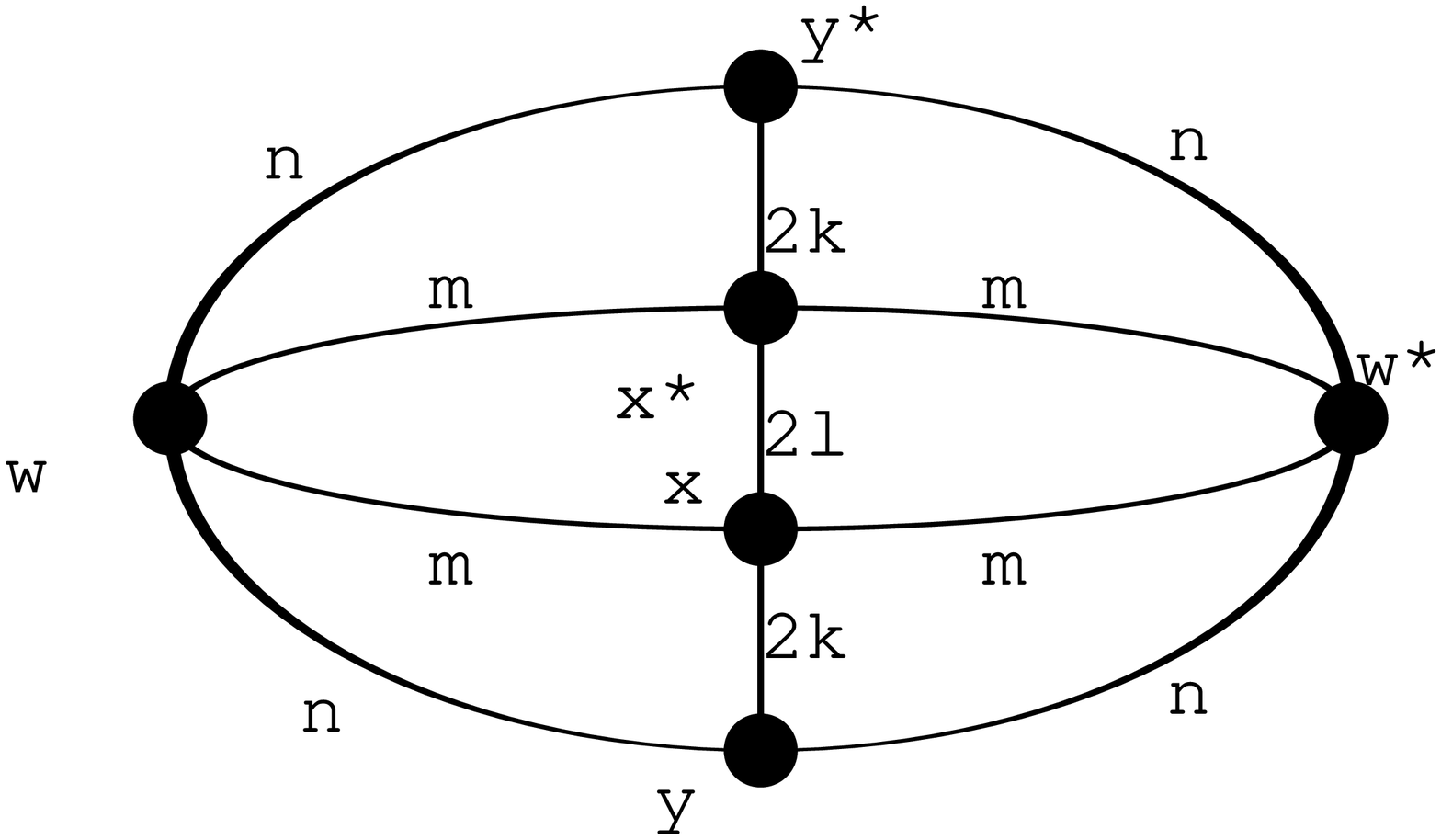}
}
= \gamma (s)
\]
where $\{w_\alpha\}_\alpha$ is an orthonormal basis of $P_{\vlon (m+n)}$ with respect to the canonical trace, $s$ is the element $P_{
\psfrag{x}{$x$}
\psfrag{x*}{$x^*$}
\psfrag{2l}{$2l$}
\psfrag{m}{$2m+2k$}
\psfrag{e}{$\nu$}
\includegraphics[scale=0.2]{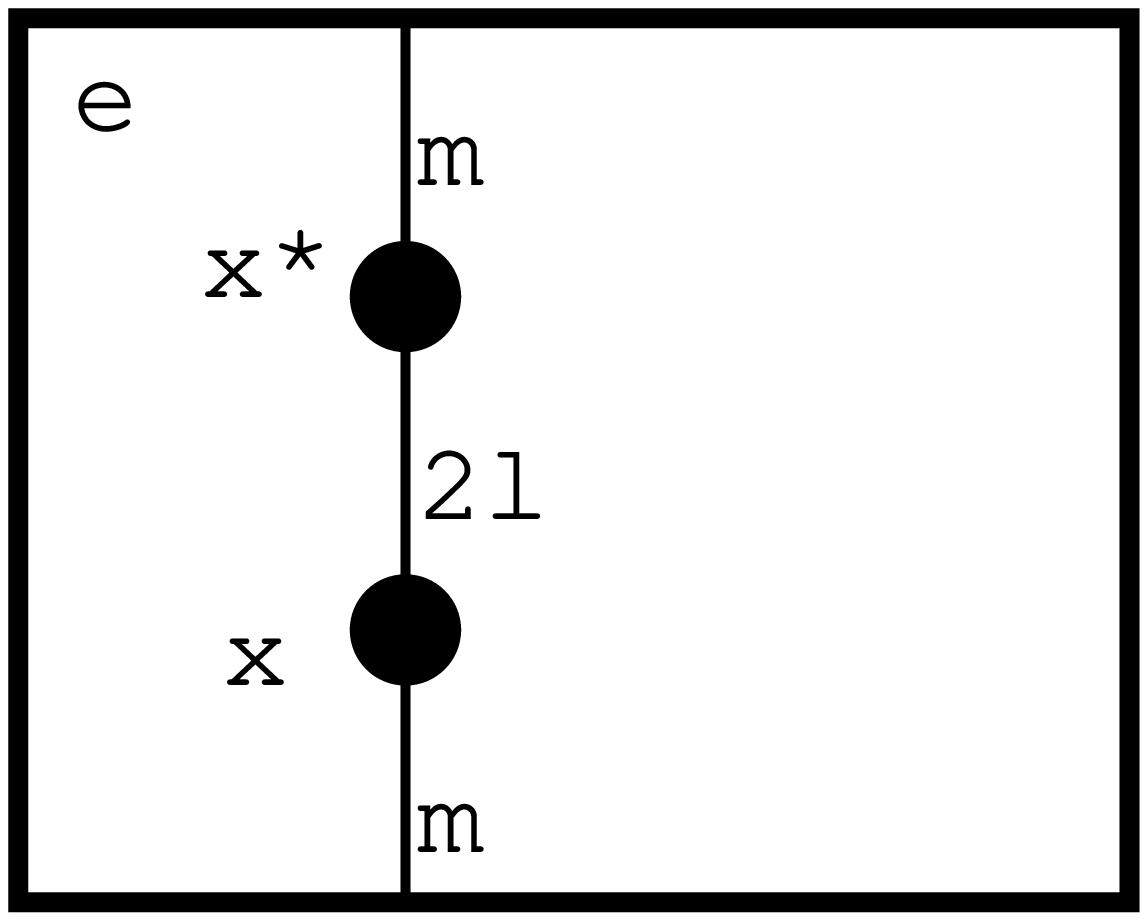}
} \in P_{\nu 2(m+k)}$ and $\gamma : P_{\nu 2(m+k)} \ra \C$ is the linear functional induced by the $P$-action of (the linear combination of semi-labelled tangles)
$\delta^{-(m+n)} \us{\alpha}{\sum} P_{\!\!\!\!\!\!\!
\psfrag{w}{$w_\alpha$}
\psfrag{w*}{$w^*_\alpha$}
\psfrag{s}{$\nu$}
\psfrag{y}{$y$}
\psfrag{y*}{$y^*$}
\psfrag{2k}{$2k$}
\psfrag{2l}{$2l$}
\psfrag{m}{$m$}
\psfrag{2m}{$2m$}
\psfrag{n}{$n$}
\includegraphics[scale=0.2]{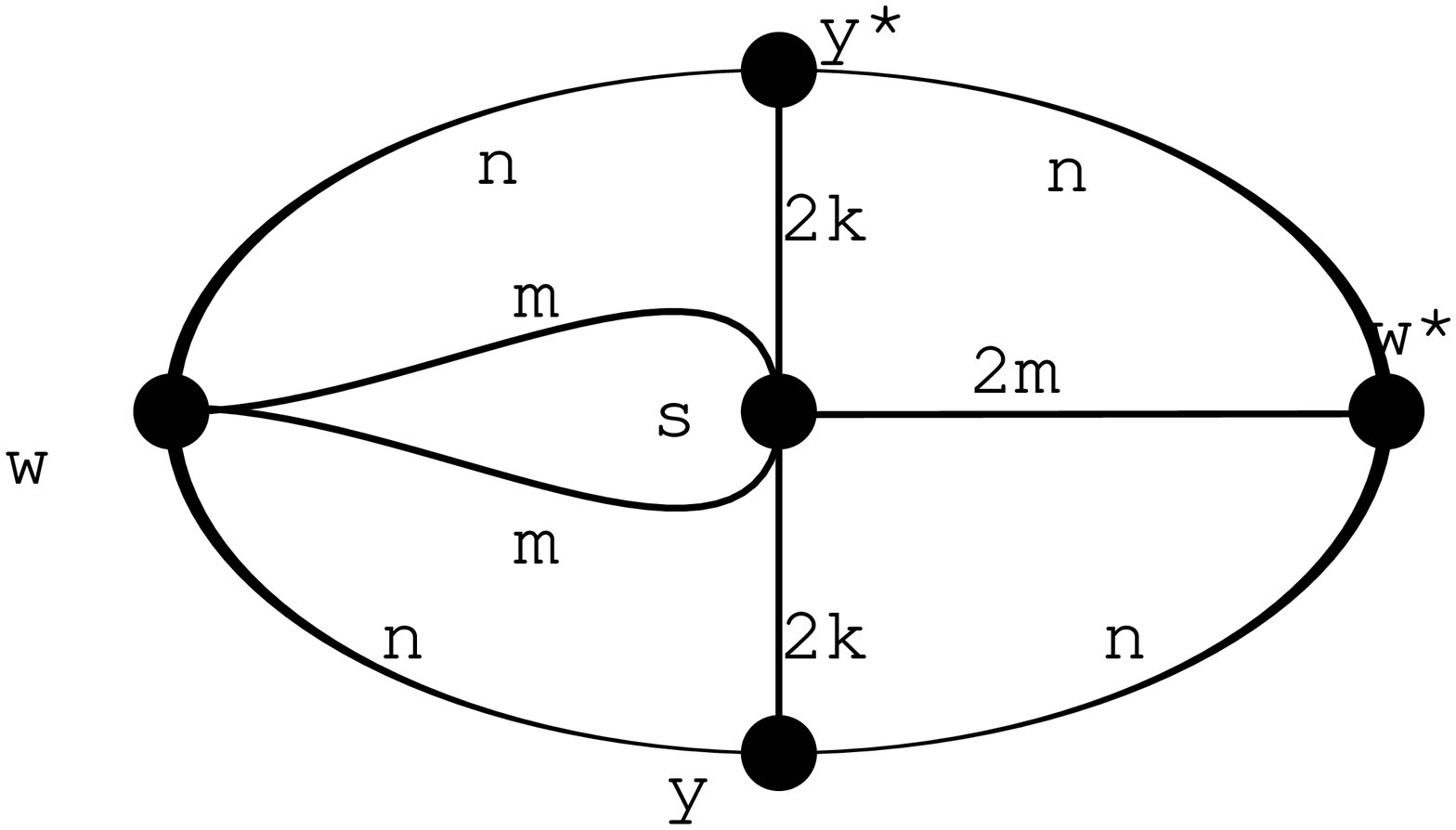}
}$.
Note that $s$ is a positive element of $P_{\nu 2(m+k)}$.
Also, $\gamma$ is positive semi-definite because for a positive $t \in P_{\nu 2(m+k)}$, we have $\gamma (t)
=\delta^{-(m+n)} \us{\alpha}{\sum} P_{
\psfrag{w}{$w_\alpha$}
\psfrag{w*}{$w^*_\alpha$}
\psfrag{x}{$t^{\frac{1}{2}}$}
\psfrag{x*}{$t^{\frac{1}{2}}$}
\psfrag{y}{$y$}
\psfrag{y*}{$y^*$}
\psfrag{2k}{$2k$}
\psfrag{2l}{$2k+2m$}
\psfrag{m}{$m$}
\psfrag{n}{$n$}
\includegraphics[scale=0.20]{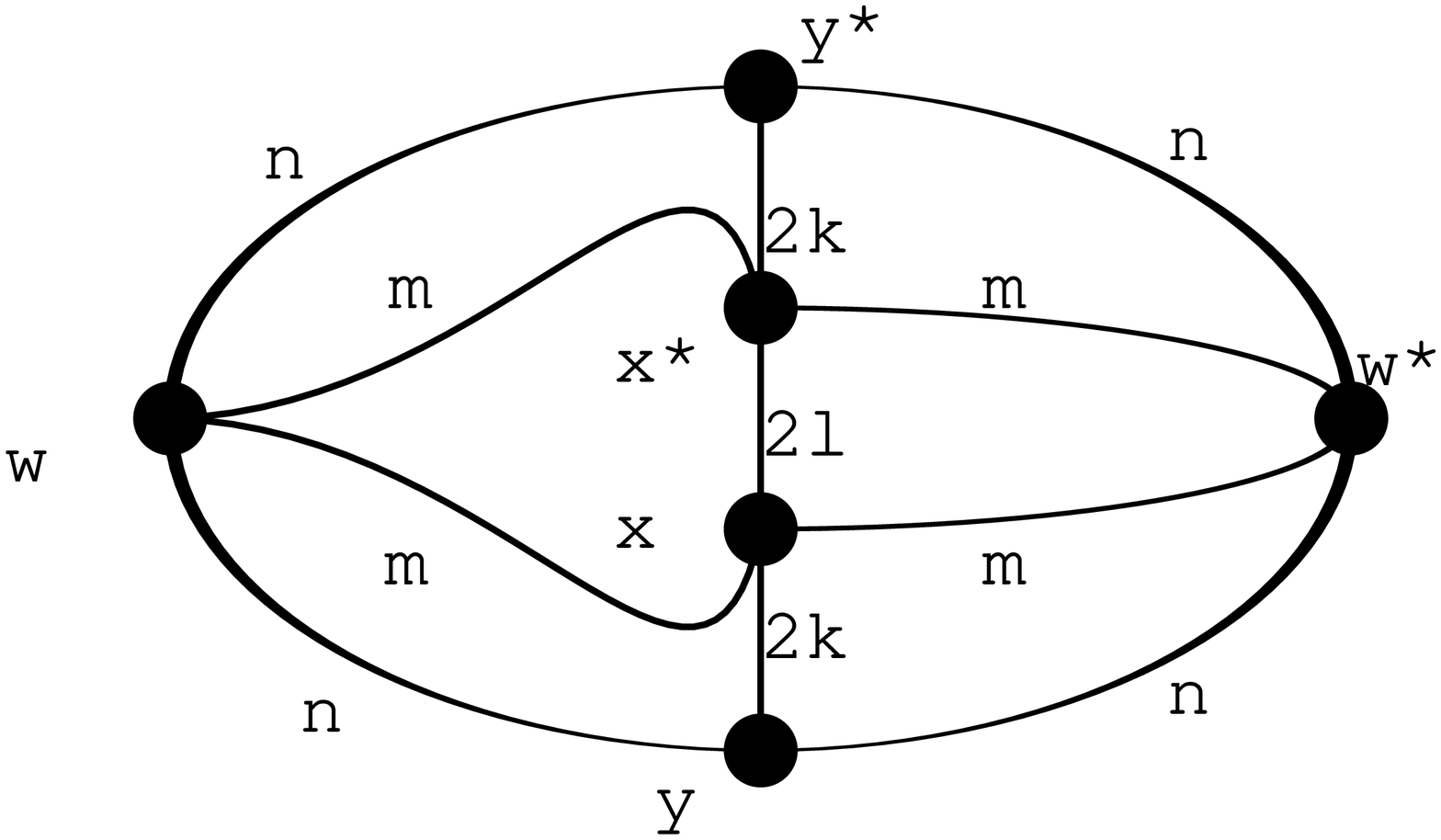}
}
= \left\| \psi^m_{\eta k , \nu (m+k)} \left( t^{1/2} \right) \circ h \right\|^2 \geq 0$ where the norm comes from the inner product in the first part.
Thus,
\[
\| a \circ h \|^2 \leq \| s \| \gamma (1) = \| s \| \delta^{-(m+n)} \us{\alpha}{\sum} P_{
\psfrag{w}{$w_\alpha$}
\psfrag{w*}{$w^*_\alpha$}
\psfrag{s}{$\nu$}
\psfrag{y}{$y$}
\psfrag{y*}{$y^*$}
\psfrag{2k}{$2k$}
\psfrag{2l}{$2l$}
\psfrag{m}{$m$}
\psfrag{2m}{$2m$}
\psfrag{n}{$n$}
\includegraphics[scale=0.20]{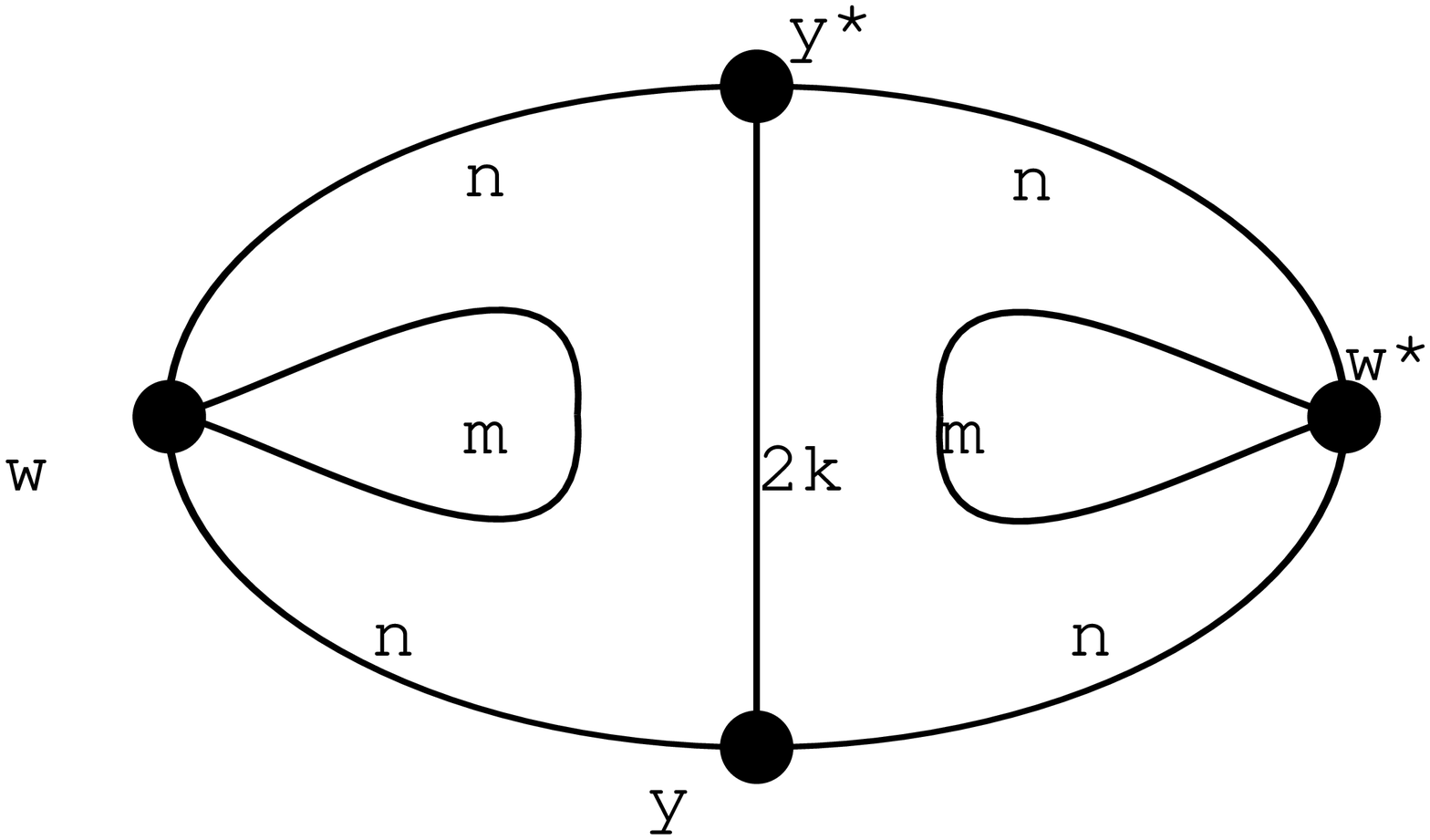}
}.
\]
We will now choose a special orthonormal basis of $P_{\vlon (m+n)}$. 
Let $\left\{ E^i_{\beta , \gamma} : 0 \leq i \leq n , 1 \leq \beta , \gamma \leq d_i \right\}$ be a system of matrix units of the finite dimensional $C^*$-algebra $P_{\vlon n}$.
Note that $\left\{ v^i_{\beta , \gamma} := c_i E^i_{\beta , \gamma} \right\}_{i,\beta,\gamma}$ is an orthonormal basis in $P_{\vlon n}$ where $c_i$'s are normalizing scalars, and thereby, $\left\{ P_{RI^m_{\vlon n}} \left( v^i_{\beta , \gamma} \right) \right\}_{i,\beta,\gamma}$ forms an orthonormal set in $P_{\vlon (m+n)}$.
On the other hand, any $w \in P_{\vlon (m+n)}$ which is orthogonal to this set, must satisfy $\left\langle v^i_{\beta,\gamma} , P_{RE^m_{\vlon (m+n)}} (w) \; v^{i'}_{\beta' , \gamma'} \right\rangle = 0$ where $RE^m_{\vlon (m+n)} : \vlon (m+n) \ra \vlon n$ is the tangle obtained from the `right conditional expectation tangle' $RE_{\vlon (m+n)}$ (described in Figure \ref{tangles}) replacing the single string with both endpoints attached to the internal disc, by $m$ many parallel strings; thus, $P_{RE^m_{\vlon (m+n)}} (w) = 0$.
This implies
\[
\|a \circ h \|^2 \leq \| s \| \delta^m \delta^{-n} \us{i,\beta,\gamma}{\sum} P_{
\psfrag{w}{$v^i_{\beta , \gamma}$}
\psfrag{w*}{$(v^i_{\gamma , \beta})^*$}
\psfrag{s}{$\nu$}
\psfrag{y}{$y$}
\psfrag{y*}{$y^*$}
\psfrag{2k}{$2k$}
\psfrag{2l}{$2l$}
\psfrag{m}{$m$}
\psfrag{2m}{$2m$}
\psfrag{n}{$n$}
\includegraphics[scale=0.2]{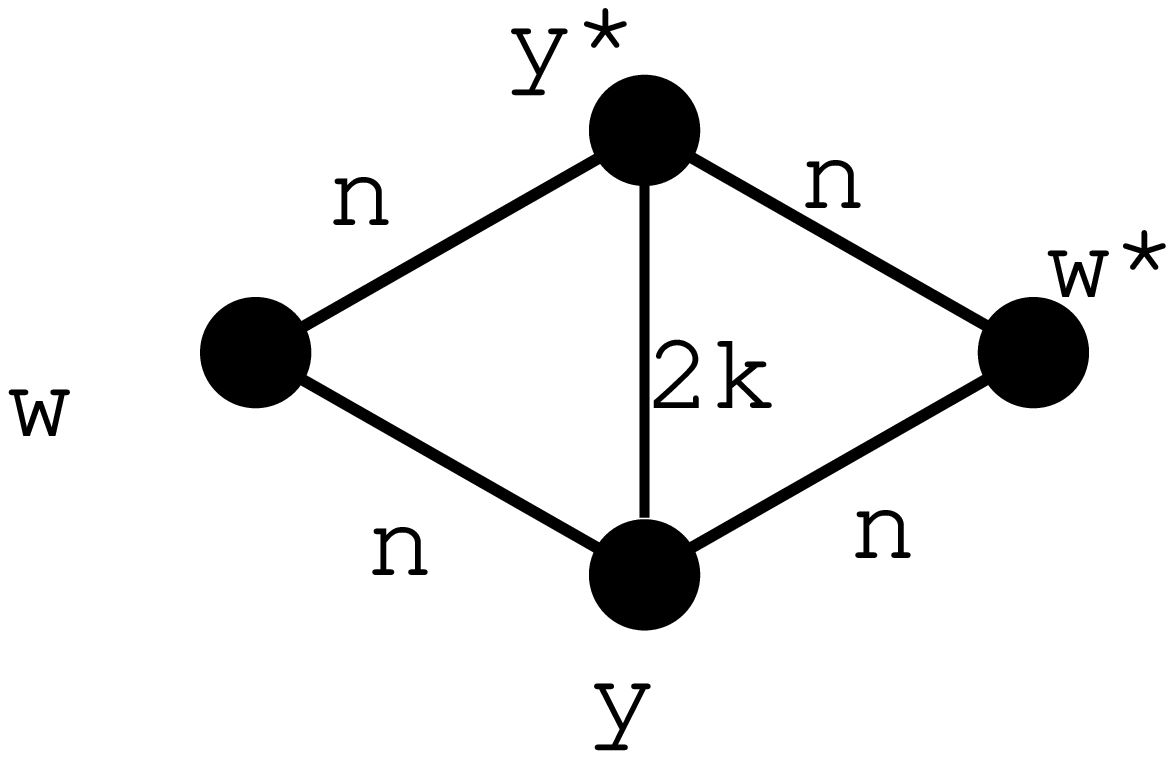}
} \;\;\;\;\;\;\; = \delta^m \|s\| \|h\|^2.
\]
Clearly, $\| s \|$ is independent of $h$. Hence, the action of $a$ is bounded.
\end{pf}
\begin{cor}\label{findep}
If $P$ has finite depth, then for every irreducible $AP_{\vlon 0 , \vlon 0}$-module $G$, there exists a unique (upto affine module isomorphism) irreducible Hilbert affine submodule of $H^\vlon$, with the $\vlon 0$ space being isomorphic to $G$ as an $AP_{\vlon 0 , \vlon 0}$-module.
Moreover, any  irreducible Hilbert affine $P$-module with weight zero is isomorphic to a submodule of $H^+$ or $H^-$.
\end{cor}
\begin{pf}
Finiteness of the depth of $P$ and positive definiteness of $\omega_\vlon$ provide $AP_{\vlon 0 , \vlon 0}$ with a finite dimensional $C^*$-algebra structure (using \cite[Proof of Theorem 6.11]{Gho}).
Now, $H^\vlon_{\vlon 0} = AP_{\vlon 0 , \vlon 0}$ is the regular $AP_{\vlon 0 , \vlon 0}$-module, and by Wedderburn-Artin, $H^\vlon_{\vlon 0}$ contains all irreducible $AP_{\vlon 0 , \vlon 0}$-modules (and hence, $G$ too) as submodules.
By Remark \ref{irrmod} (2), the submodule $[G]$ of $H^\vlon$, generated by $G$, is irreducible.
Uniqueness follows from Remark \ref{affmodmor}.

For the second statement, consider an irreducible Hilbert affine $P$-module $V$ with weight zero. Without loss of generality, let $V_{+ 0} \neq \{0\}$ which is also irreducible $AP_{+1,+1}$-module. By Remarks \ref{irrmod} and \ref{affmodmor} and the first part, $V=[V_{+ 0}]$ sits inside $H^+$ as a submodule.
\end{pf}

\vspace{2mm}
Next, we will investigate Hilbert affine $P$-modules which are generated by their $(+0)$- or $(-0)$- spaces where depth of $P$ is not necessarily finite.
The finite depth case is completely determinded by Corollary \ref{findep} which will not work in infinite depth because any irreducible $AP_{\vlon 0 , \vlon 0}$-module might not be isomorphic to a submodule of $H^\vlon_{\vlon 0}$.
However, the easiest example of an irreducible Hilbert affine $P$-module, namely, the planar algebra $P$ itself, does sit inside both $H^+$ and $H^-$ as a submodule.
It is the submodule of $H^\vlon$ generated by the one-dimensional orthogonal complement of the kernel of the linear homomorphism (which actually gives the dimension function via the isomorphism in Theorem \ref{main-theorem})
\[
AP_{\vlon 0, \vlon 0} \ni \psi^{2k}_{\vlon 0,\vlon 0} (x) \mapsto P_{TR^r_{\vlon 2k}} (x) \in P_{\vlon 0} \cong \C
\]
for $k\in \N_0$, $x \in P_{\vlon 2k}$.
\comments{
It can also be viewed as the extension (along the lines of Corollary \bf{Pending}) of the one-dimensional ($AP_{\vlon 0, \vlon 0} , \omega_\vlon)$-module given by the linear homomorphism (which is actually the dimension function)
\[
AP_{\vlon 0, \vlon 0} \ni \psi^{2k}_{\vlon 0,\vlon 0} (x) \mapsto P_{TR^r_{\vlon 2k}} (x) \in P_{\vlon 0} \cong \C
\]
for $k\in \N_0$, $x \in P_{\vlon 2k}$.
}

Let us denote the completion of $H^{\varepsilon}$ by $\mcal H^{\varepsilon}$ with $\mcal H_{a}^{\varepsilon}$ being the unique extension of $H_{a}^{\varepsilon}$ for all affine morphisms $a$ (see Theorem \ref{hethm}).
Then, $L^{\varepsilon}:=\left(\mcal H_{AP_{\varepsilon0,\varepsilon0}}^{\varepsilon}\right)^{\prime\prime}\subset\mcal B\left(\mcal H_{\vlon0}^{\varepsilon}\right)$ becomes a finite von Neumann algebra on which $\omega_{\varepsilon}$ extends to a faithful normal tracial state given by $\tilde \omega_\vlon := \left\langle \hat 1 , \cdot \left( \hat 1 \right) \right \rangle : L^\vlon \ra \C$.
Note that $H^\vlon_{\eta k}$ has a right $AP_{\vlon 0 , \vlon 0}$-module structure.
Now, for all $a \in AP_{\vlon 0 , \vlon 0}$, $b \in AP_{\vlon 0 , \eta k}$, we have
\[
\norm {b \circ a}^2 = \omega_\vlon (a^* \circ b^* \circ b \circ a) = \tilde \omega_\vlon (y \mcal H^{\vlon}_{a \circ a^*} y) \leq \norm {{\mcal H}^{\vlon}_a}^2 \tilde \omega_\vlon (y^2) = \norm {{\mcal H}^{\vlon}_a}^2 \norm b^2 
\]
where $y \in L^\vlon$ is the positive square root of $\mcal H^{\vlon}_{b^* \circ b}$.
So, for all $\eta k \in \t{Col}$, the right action of any element $a \in AP_{\vlon 0 , \vlon 0}$ on $H^\vlon_{\eta k}$ is bounded as well; let $\rho^\vlon_{\eta k} (a) \in \mcal B (\mcal H^\vlon_{\eta k})$ denote its unique extension.
\begin{lem}\label{ract}
For all $\eta k\in\mbox{Col}$, the anti-algebra $*$-homomorphism $\rho^\vlon_{\eta k} : AP_{ \vlon0 , \vlon0}  \ra \mcal B ( \mcal H_{\eta k}^{\vlon} )$ extends to a normal anti-algebra $*$-homomorphism from $L^{\vlon}$ to $\mcal B ( \mcal H_{\eta k}^{\vlon} )$. Moreover, it is  faithful for all $\eta k \neq -\vlon 0$.
\end{lem}
\begin{pf}
Note that $\mcal H_{a}^{\varepsilon}\circ\rho^\vlon_{\eta k}(b)=\rho^\vlon_{\nu l}(b)\circ\mcal H_{a}^{\varepsilon}$ for all $a\in AP_{\eta k,\nu l}$, $b\in AP_{\varepsilon0,\varepsilon0}$.
Set $W_{\eta k}:=\left(\rho^\vlon_{\eta k}\left(AP_{\varepsilon0,\varepsilon0}\right)\right)^{\prime\prime}\subset\mcal B\left(\mcal H_{\eta k}^{\varepsilon}\right)$.
Since $L^{\varepsilon}\supset\mcal H_{AP_{\varepsilon0,\varepsilon0}}^{\varepsilon}\ni\mcal H_{a}^{\varepsilon}\mapsto\rho^\vlon_{\varepsilon0}\left(a\right)=J\mcal H_{a^{*}}^{\varepsilon}J\in\rho^\vlon_{\varepsilon0}\left(AP_{\varepsilon0,\varepsilon0}\right)\subset W_{\varepsilon0}$ is an anti-algebra $*$-isomorphism (where $J$ is the canonical conjugate linear unitary on $\mcal H^\vlon_{\vlon 0}$), it is enough to show that
\[
W_{\varepsilon0}\supset\rho^\vlon_{\varepsilon0}\left(AP_{\varepsilon0,\varepsilon0}\right)\ni\rho^\vlon_{\varepsilon0}\left(a\right)\os{\alpha_{\eta k}}{\longmapsto}\rho^\vlon_{\eta k}\left(a\right)\in\rho^\vlon_{\eta k}\left(AP_{\varepsilon0,\varepsilon0}\right)\subset W_{\eta k}
\]
extends to a surjective normal $*$-homomorphism which is also injective for all $\eta k \neq -\vlon 0$.

\noindent \textbf{Case 1:}
Suppose $\eta k=\varepsilon 0$. This case is trivial.

\noindent \textbf{Case 2:}
Suppose $k>0$.
Let $c_{\varepsilon0,\eta k} \in AP_{\varepsilon0,\eta k}$ denote the affine tangle
\psfrag{e}{$\vlon$}
\psfrag{-e}{$-\vlon$}
\psfrag{k}{$k$}
\psfrag{k-1}{$k-1$}
\includegraphics[scale=0.25]{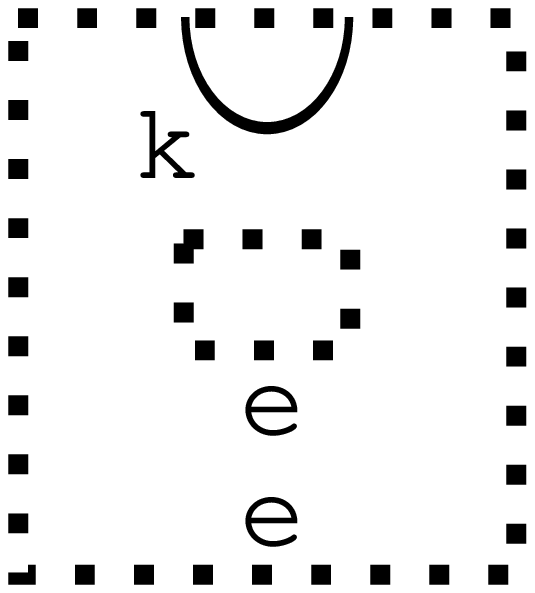}
or
\psfrag{e}{$\vlon$}
\psfrag{-e}{$-\vlon$}
\psfrag{k}{$k$}
\psfrag{k-1}{$k-1$}
\includegraphics[scale=0.25]{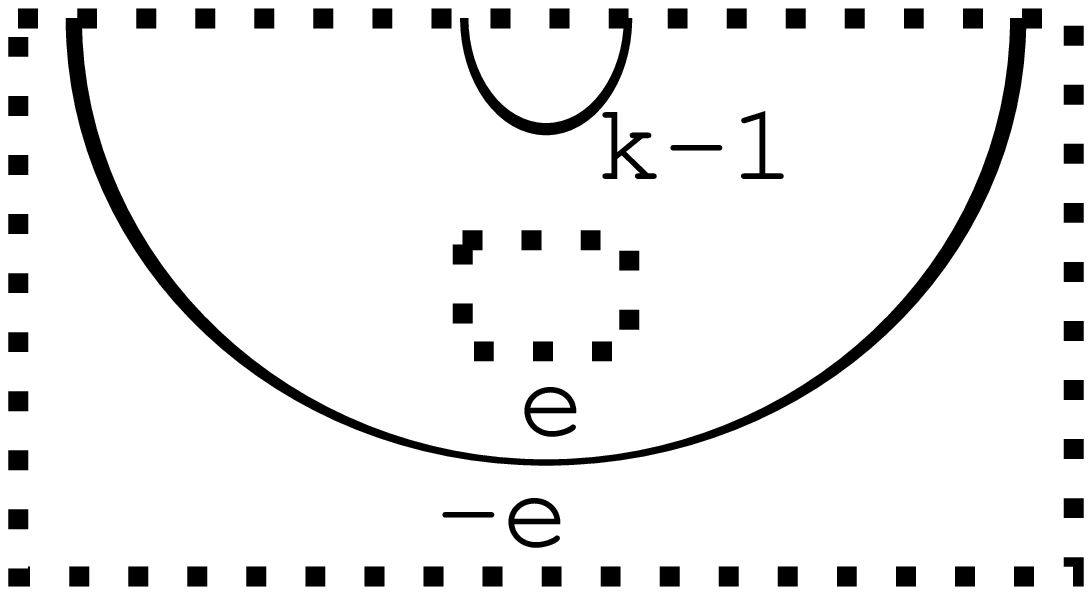}
according as $\eta = \vlon$ or $\eta = -\vlon$.
Note that $U:=\delta^{-k/2}\mcal H_{c_{\varepsilon0,\eta k}}^{\vlon}:\mcal H_{\vlon0}^{\vlon}\ra\mcal H_{\eta k}^{\vlon}$ is an isometry.
Let $p:=UU^{*}=\delta^{-k}\mcal H_{c_{\varepsilon0,\eta k}\circ c_{\varepsilon0,\eta k}^{*}}^{\vlon}\in \mscr P\left(\mcal H_{\eta k}^{\vlon}\right)$; clearly, $p\in W_{\eta k}^{\prime}$. It is easy to check that the central support of $p$ in $W_{\eta k}$ is $1$ (using the fact $\mcal H_{AP_{\eta k,\eta k}}^{\vlon}\rho^\vlon_{\eta k}(AP_{\vlon0,\vlon0})p(H_{\eta k}^{\vlon})=H_{\eta k}^{\vlon}$).
Thus, $W_{\eta k}\ni x\mapsto xp\in pW_{\eta k}$ is an isomorphism.
This gives us an injective $*$-algebra homomorphism $W_{\eta k}\ni x\os{\alpha}{\mapsto}U^{*}xU\in\mcal B\left(\mcal H_{\varepsilon0}^{\varepsilon}\right)$.
$\t{Range} \,\alpha$ is a von Neumann algebra since $\alpha$ is normal.
On the other hand, $\alpha\left(\rho^\vlon_{\eta k}(a)\right)=\rho^\vlon_{\varepsilon0}(a)$ for all $a\in AP_{\vlon0,\varepsilon0}$.
So, $\mbox{Range} \,\alpha=W_{\varepsilon0}$.
Hence, $\alpha_{\eta k}$ is given by $\alpha^{-1}$.

\noindent \textbf{Case 3:} Suppose $\eta k=-\varepsilon0$. It is
enough to show that
\[
W_{\varepsilon1}\supset\rho^\vlon_{\varepsilon1}\left(AP_{\varepsilon0,\varepsilon0}\right)\ni\rho^\vlon_{\varepsilon1}\left(a\right)\mapsto\rho^\vlon_{-\varepsilon0}\left(a\right)\in W_{-\varepsilon0}
\]
extends to a normal $*$-homomorphism.
For this, set $c_{-\varepsilon0,\varepsilon1} :=
\psfrag{e}{$\vlon$}
\psfrag{-e}{$-\vlon$}
\psfrag{k}{$k$}
\psfrag{k-1}{$k-1$}
\includegraphics[scale=0.25]{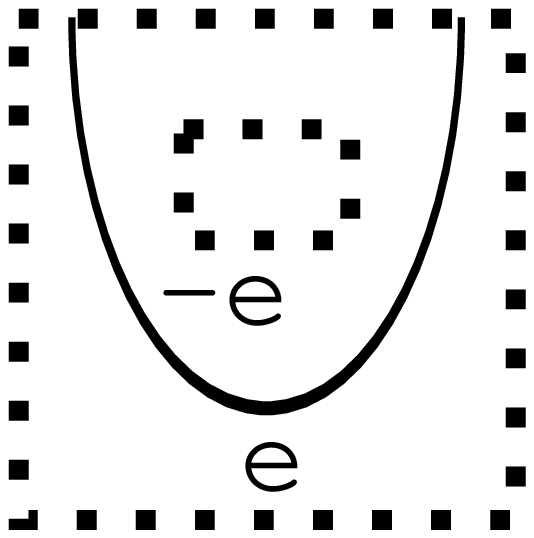}
\in AP_{-\varepsilon0,\varepsilon1}$ and $U:=\delta^{-1/2} \mcal H_{ c_{-\varepsilon0 , \varepsilon1}}^{\varepsilon} : \mcal H_{-\varepsilon0}^{\vlon}\ra\mcal H_{\varepsilon1}^{\vlon}$.
Note that $U^{*}U=1$.
Let $p:=UU^{*}=\delta^{-1} \mcal H_{c_{-\varepsilon0,\varepsilon1}\circ c_{-\varepsilon0,\varepsilon1}^{*}} \in \mscr P \left( \mcal H_{\varepsilon1}\right)$; clearly, $p\in W_{\varepsilon1}^{\prime}$.
So, there exists a normal $*$-homomorphism $W_{\varepsilon1}\ni x\os{\alpha}{\mapsto}U^{*}xU\in\mcal B \left( \mcal H_{-\varepsilon0}\right)$.
Note that $\alpha\left(\rho^\vlon_{\varepsilon1}(a)\right)=\rho^\vlon_{-\varepsilon0}\left(a\right)$ for all $a\in AP_{\varepsilon0,\varepsilon0}$; this implies $\mbox{Range} \,\alpha=W_{-\varepsilon0}$.
\comments{Now, the range of the central support of $p$ (say $z_p$) in $W_{\varepsilon1}$, contains at least $\ol{\rho^\vlon_{\varepsilon1}\left(I^\vlon \right)H^\vlon_{\varepsilon1}}$ since
\begin{align*}
\mcal H_{AP_{\vlon1,\vlon1}} \rho^\vlon_{\vlon1}\left(AP_{\varepsilon0,\varepsilon0}\right) p\left[\left(AP_{\varepsilon0,\varepsilon1}\right)^{\widehat{}}\right] & =\left(AP_{\varepsilon1,\varepsilon1}\circ c_{-\varepsilon0,\varepsilon1}\circ c_{-\varepsilon0,\varepsilon1}^{*}\circ AP_{\varepsilon0,\varepsilon1}\circ AP_{\vlon0,\vlon0}\right)^{\widehat{}}\\
& =\left(AP_{-\varepsilon0,\varepsilon1}\circ AP_{\vlon0,-\vlon0}\right)^{\widehat{}} \supset \rho^\vlon_{\varepsilon1}\left(I^\vlon \right)\left[\left(AP_{\varepsilon0,\varepsilon1}\right)^{\widehat{}}\right] = \rho^\vlon_{\varepsilon1}\left(I^\vlon \right)H^\vlon_{\vlon 1}.\end{align*}
This along with $\t{Range} \, z_{\vlon 1} = \ol{\rho^\vlon_{\vlon 1} \left( I^\vlon \right) H^\vlon_{\vlon 1}}$, implies $z_p \geq z_{\vlon 1}$.
For the reverse inequality, note that
Note that $\alpha\left(z_{\varepsilon1}\right)=1$.
Thus, $\mbox{ker }\alpha$ contains $\left(1-z_{\varepsilon1}\right)W_{\varepsilon1}$.
Equality follows from the injectivity of $z_{\varepsilon1}W_{\varepsilon1}\ni x\mapsto xp\in pW_{\varepsilon1}$.
}
\end{pf}

We now proceed towards finding the kernel of the extension of $\rho^\vlon_{-\vlon 0}$ to $L^\vlon$, which we denote with the same symbol.
For this, consider the $*$-closed two sided ideal $I^{\vlon} := AP_{-\vlon0, \vlon0} \circ AP_{\vlon0 , - \vlon 0}$ in $AP_{\vlon0,\vlon0}$.
Thus, $\ol{\mcal H^\vlon_{I^\vlon} }$ (with respect to weak operator topology (WOT) in $\mcal B \mcal (H^\vlon_{\vlon 0} )$) becomes a $*$-closed, WOT-closed two-sided ideal in $L^\vlon$; let  $z_\vlon$ be the central projection of $L^\vlon$ such that $\ol{\mcal H^\vlon_{I^\vlon} } = z_\vlon L^\vlon$.
\begin{lem}
$\mbox{ker} \, \rho^\vlon_{- \vlon 0} = (1 - z_\vlon ) L^\vlon$.
\end{lem}
\begin{pf}
If $x\in L^\vlon$, then $\rho^\vlon_{-\vlon 0} (x) = 0$ if and only if
\begin{align*}
0 &
= \langle \hat c , \rho^\vlon_{-\vlon 0} (x) \hat d \rangle
= \langle \hat 1 , \mcal H^\vlon_{c^*} \rho^\vlon_{-\vlon 0} (x) \mcal H^\vlon_d \hat 1 \rangle
= \langle \hat 1 , \mcal H^\vlon_{c^* \circ d} \rho^\vlon_{\vlon 0} (x) \hat 1 \rangle
= \langle \hat 1 , \mcal H^\vlon_{c^* \circ d} \rho^\vlon_{\vlon 0} (x) \hat 1 \rangle
= \langle \hat 1 , \mcal H^\vlon_{c^* \circ d} x \hat 1 \rangle\\
& = \tilde \omega_\vlon (\mcal H^\vlon_{c^* \circ d} \, x)
\end{align*}
for all $c , d \in AP_{\vlon 0 , - \vlon 0}$.
Thus, by WOT-continuity of $\tilde \omega_\vlon$, we get $x \in \t{ker} \, \rho^\vlon_{-\vlon 0}$ if and only if $\tilde \omega_\vlon ( z_\vlon x x^*) = 0$ which is equivalent to $x z_\vlon = 0$ (using faithfulness of $\tilde \omega_\vlon$).
This give the required equation.
\end{pf}
\begin{thm}
For every left $L^{\vlon}$-module $\mcal K$, 
there exists a unique Hilbert affine $P$-module, say $[ \mcal K ]$, such that $\mcal K$ and $[ \mcal K ]_{\vlon 0}$ are isometrically isomorphic as $AP_{\vlon 0 , \vlon 0}$-modules and $\left[ [\mcal K]_{\vlon 0} \right] = [\mcal K]$.
Further, $[\mcal K ]_{-\vlon 0} = \{0\}$ if and only if the action of $z_\vlon$ on $\mcal K$, is zero.
\end{thm}
\begin{pf}
Uniqueness easily follows from Remark \ref{affmodmor}.
Next, we consider the space of bounded vectors, $\mcal K^{0}$ which will be dense in $\mcal K$ and have a left $L^{\vlon}$-valued inner product $\lrsuf{L^\vlon}{\langle}{} \cdot , \cdot \rangle$ satisfying $\tilde \omega_\vlon \circ \lrsuf{L^\vlon}{\langle}{} \cdot , \cdot \rangle = \langle \cdot , \cdot \rangle$.
On the other hand, $\mcal H^\vlon_{\eta k}$ gets a right $L^\vlon$-module structure from Lemma \ref{ract}; so, $\left( \mcal H^\vlon_{\eta k} \right)^0$ (the space of bounded vectors of $\mcal H^\vlon_{\eta k}$) will have a right $L^\vlon$-valued inner product compatible with $\tilde \omega_\vlon$.
We now use Connes-fusion techniques to build $[\mcal K ]_{\eta k} := \mcal H^\vlon_{\eta k} \us{L^\vlon}{\otimes} \mcal K$.
The action of $a \in AP_{\eta k , \nu l}$ is given by $[\mcal K]_a := \mcal H^\vlon_a \us{L^\vlon}{\otimes} \t{id}_{\mcal K} : [\mcal K ]_{\eta k} \ra [\mcal K ]_{\nu l}$.
This makes $[\mcal K]$ a Hilbert affine $P$-module.

For the remaining part, first note that $\left( AP_{\vlon 0 , \eta k} \right)\hat{}$ sits inside $\left( \mcal H^\vlon_{\eta k} \right)^0$ and is dense in $\mcal H^\vlon_{\eta k}$. Thereby, $\t{span} \left\{ \hat{a} \us{L^\vlon}{\otimes} \zeta = [\mcal K]_a \left( \hat{1} \us{L^\vlon}{\otimes} \zeta \right) : a \in AP_{\vlon 0,\eta k} , \zeta \in \mcal K^0 \right\}$ becomes a dense subset in $[\mcal K]_{\eta k}$.
Thus, $\left[ [\mcal K]_{\vlon 0} \right] = [\mcal K]$.
The map $\mcal K^0 \ni \zeta \mapsto \hat{1} \us{L^\vlon}{\otimes} \zeta \in [\mcal K]_{\vlon 0}$ extends to a surjective $AP_{\vlon 0 , \vlon 0}$-linear isometry from $\mcal K$ to $[\mcal K]_{\vlon 0}$.

Observe that $[\mcal K]_{-\vlon 0} = 0$ if and only if $
0 = \langle \hat{a} \us{L^\vlon}{\otimes} \zeta , \hat{b} \us{L^\vlon}{\otimes} \zeta \rangle = \langle \zeta , (a^* \circ b) \zeta \rangle$ for all $a , b \in AP_{\vlon 0, -\vlon 0}$, $\zeta \in \mcal K^0$.
Since the representation $L^\vlon \ra \mcal B (\mcal K)$ is normal, this is equivalent to $z_\vlon (\mcal K^0) = \{ 0 \}$ and hence, we get the required result.
\end{pf}

From the above theorem, we wonder whether every $*$-affine $P$-module $V$ which is generated by $V_{\vlon 0}$, can be obtained in this way of extending an $L^\vlon$-module.
The trivial module $P$ is the extension of the one-dimensional $AP_{\vlon 0 , \vlon 0}$-module given by the dimension function.
Another question along this line is whether we can do similar analysis for $AP_{\vlon k , \vlon k}$ where $k > 0$.
\begin{rem}
Note that the spaces of affine morphisms and `annular morphisms' (see \cite{Jon01}) with the color of internal or external rectangle being $\pm0$, are canonically isomorphic (because there will not be any difference between affine isotopy and the usual planar isotopy in such cases).
So, all results on affine category over $P$ and affine $P$-modules, obtained in Sections \ref{mt} and \ref{regmod}, also hold for annular category over $P$ and annular representations.
\end{rem}
\section{Affine modules of irreducible depth two planar algebras}\label{JWconjirrdep2}
This section deals with irreducible depth two subfactor planar algebras.
Such planar algebras are the ones associated to the subfactors arising from action of finite dimensional Kac algebras.
We will establish an equivalence between the category of affine modules over such a planar algebra and the representation category of the quantum double of the corresponding Hopf algebra, and thereby, confirming Jones-Walker conjecture in this case.

Throughout this section, $P$ will denote an irreducible depth two subfactor planar algebra and $\vlon = \pm$.
\subsection{Affine morphisms at level one}\label{affmorlevel1}
By \cite{Sz}, $H_\vlon := P_{\vlon 2}$ has a Kac algebra structure.
We will first briefly recall this structure in the language of planar algebras (see \cite{KLS}, \cite{DK} for details).
Suppose $G_{\vlon k} =
\psfrag{D1}{$D_1$}
\psfrag{D2}{$D_2$}
\psfrag{Dk}{$D_k$}
\psfrag{...}{$\cdots$}
\psfrag{.}{.}
\psfrag{2}{$2$}
\psfrag{e}{$\vlon$}
\psfrag{gek}{$G_{\vlon k}$}
\vcenter{
\includegraphics[scale=0.20]{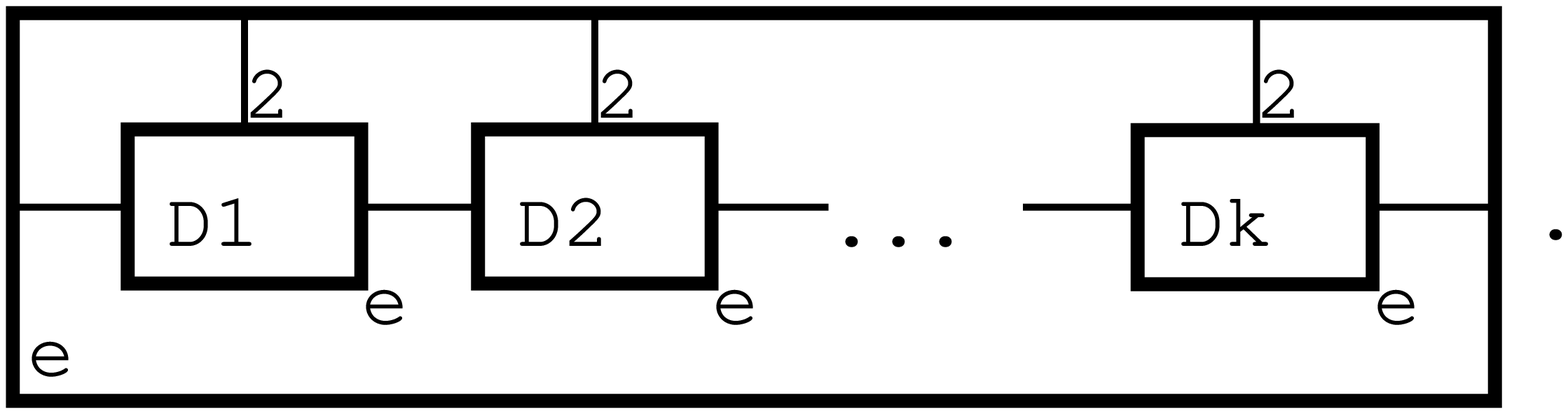}
}$
\begin{rem}\label{irrdep2gen} 
Since depth of $P$ is $2$, $\text{Range} \, P_{G_{\vlon k}} = P_{\vlon(k+1)}$ for all $k \geq 1$.
This along with irreducibility of $P$ gives $\t{dim}_\C (P_{\vlon (k+1)}) = \left[ \t{dim}_\C ( P_{\vlon 2}) \right]^k$ which implies that $P_{G_{\vlon k}} : \left( P_{\vlon 2} \otimes P_{\vlon 2} \otimes \cdots k \t{ factors} \right) \ra P_{\vlon (k+1)}$ is an isomorphism for all $k\geq 1$.
\end{rem}
We already know the $C^*$-algebra structure on $H_\vlon$. We now define the comultiplication map $\Delta_\vlon : H_\vlon \ra H_\vlon \otimes H_\vlon$; we will use Sweedler's notation, namely, $\Delta_\vlon (x) = x_{(1)} \otimes x_{(2)}$ which is determined by the equations
\[
P_{
\psfrag{x}{$x$}
\psfrag{+}{$+$}
\psfrag{=}{$=$}
\includegraphics[scale=0.3]{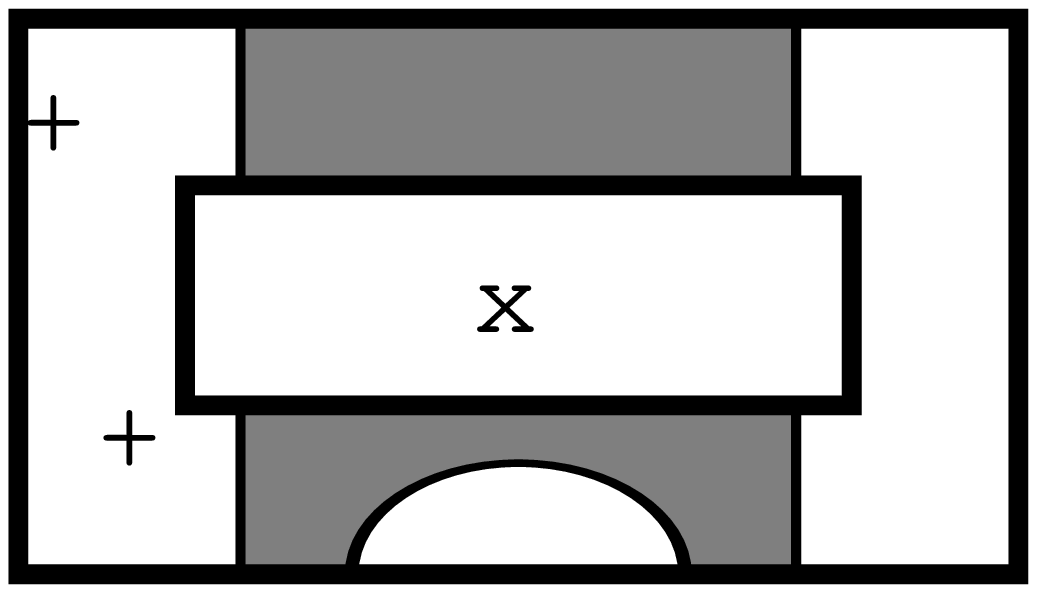}
}
=P_{
\psfrag{x1}{$x_{(1)}$}
\psfrag{x2}{$x_{(2)}$}
\psfrag{x}{$x$}
\psfrag{+}{$+$}
\psfrag{=}{$=$}
\includegraphics[scale=0.3]{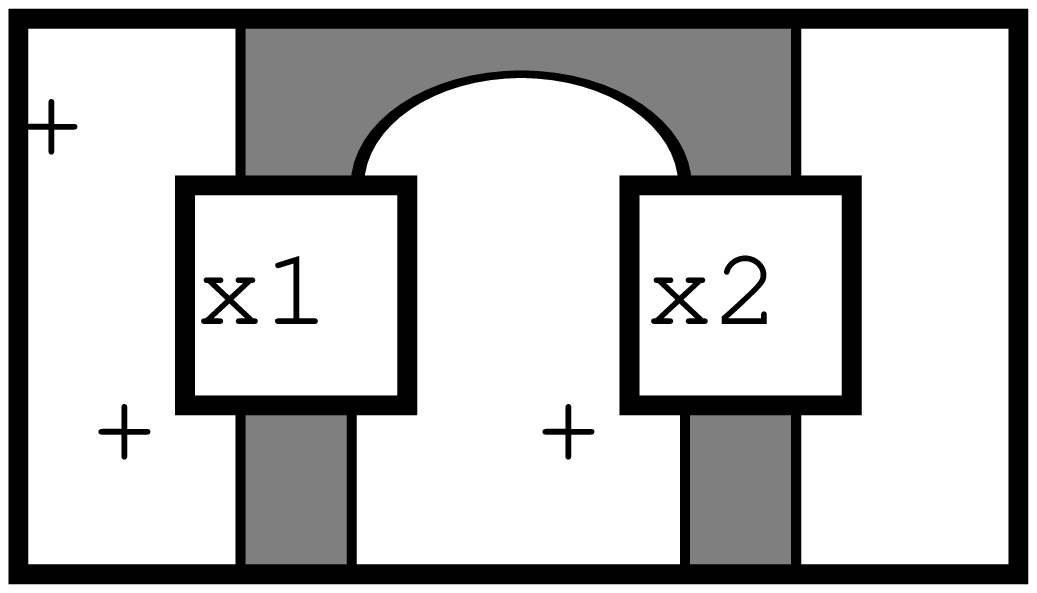}
}
\;\;\;\;\;\t{ and }\;\;\;\;\;
P_{
\psfrag{x}{$x$}
\psfrag{-}{$-$}
\psfrag{=}{$=$}
\includegraphics[scale=0.3]{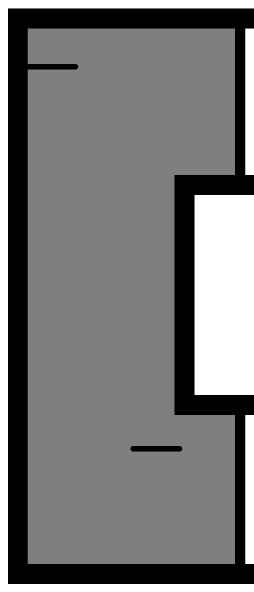}
}
=
P_{
\psfrag{x1}{$x_{(1)}$}
\psfrag{x2}{$x_{(2)}$}
\psfrag{x}{$x$}
\psfrag{-}{$-$}
\psfrag{=}{$=$}
\includegraphics[scale=0.3]{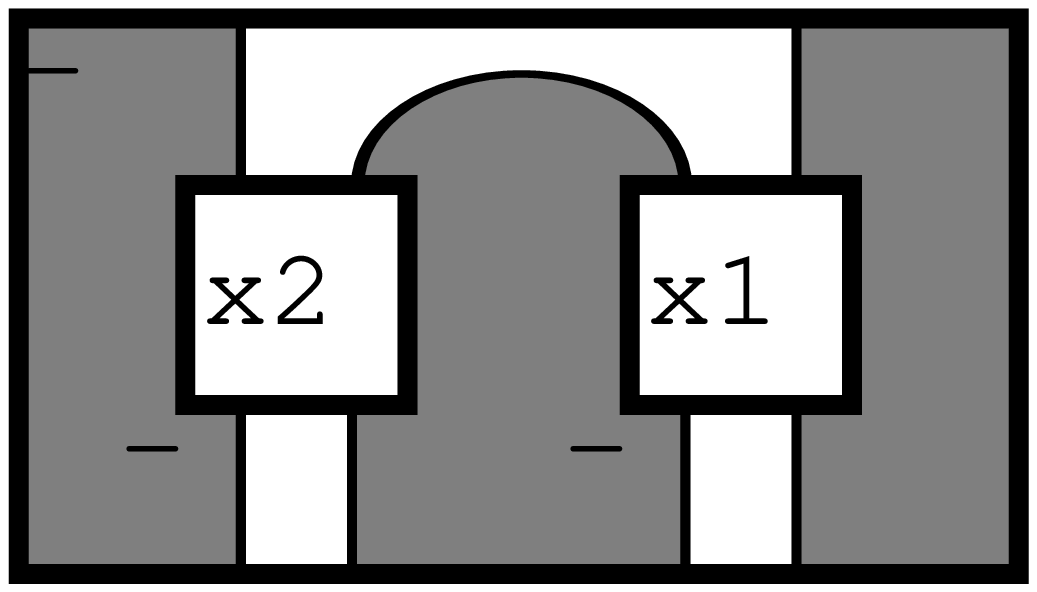}
}
\]
for $x \in H_\vlon$. The counit is given by $H_\vlon \ni x \os{\chi_\vlon}{\mapsto} \delta^{-1} P_{
\psfrag{x}{$x$}
\includegraphics[scale=0.15]{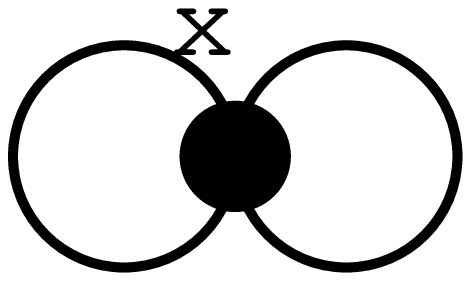}
} \in P_{\vlon 0} \cong \C$ and the antipode is $H_\vlon \ni x \os{S_\vlon}{\mapsto} P_{
\psfrag{x}{$x$}
\psfrag{e}{$\vlon$}
\psfrag{2}{$2$}
\includegraphics[scale=0.15]{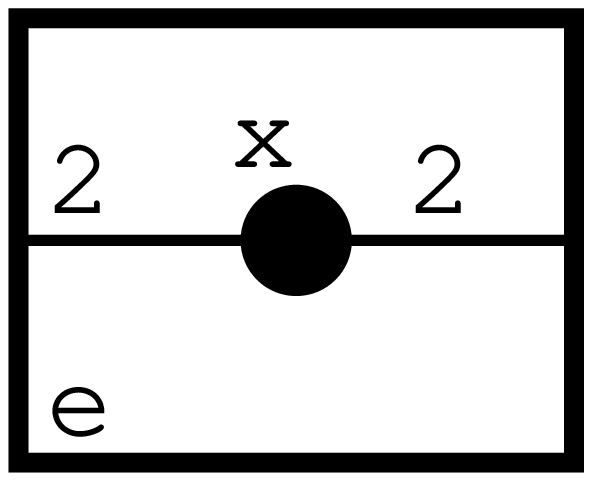}
} \in H_{\vlon}$. With these structural maps, $H_\vlon$ becomes a finite dimensional Kac algebra.
The following two relations will be very useful:
\[
P_{
\psfrag{x}{$x$}
\psfrag{+}{$+$}
\psfrag{=}{$=$}
\includegraphics[scale=0.30]{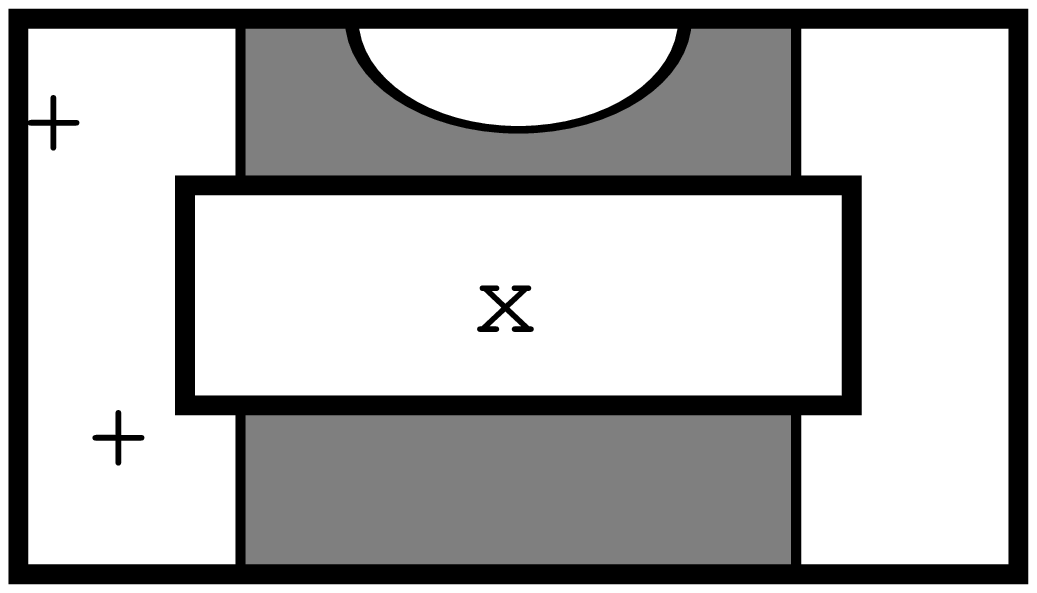}
}
=P_{
\psfrag{x1}{$x_{(1)}$}
\psfrag{x2}{$x_{(2)}$}
\psfrag{x}{$x$}
\psfrag{+}{$+$}
\psfrag{=}{$=$}
\includegraphics[scale=0.30]{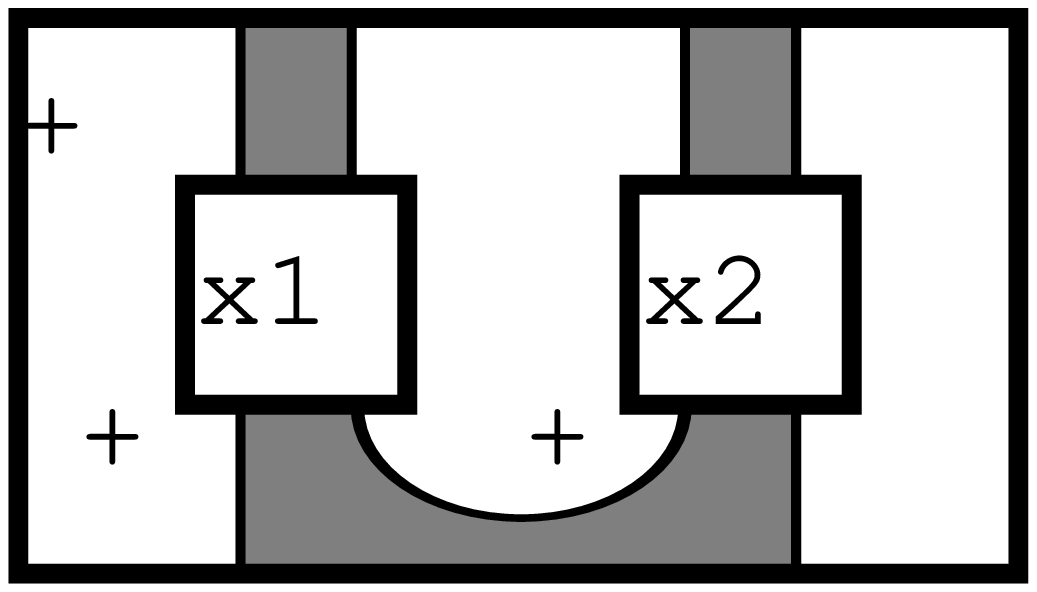}
}
\;\;\;\;\;\t{ and }\;\;\;\;\;
P_{
\psfrag{x}{$x$}
\psfrag{-}{$-$}
\psfrag{=}{$=$}
\includegraphics[scale=0.30]{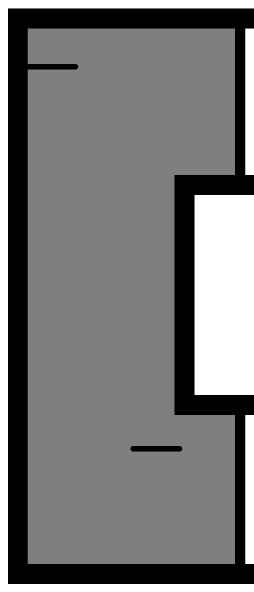}
}
=
P_{
\psfrag{x1}{$x_{(1)}$}
\psfrag{x2}{$x_{(2)}$}
\psfrag{x}{$x$}
\psfrag{-}{$-$}
\psfrag{=}{$=$}
\includegraphics[scale=0.30]{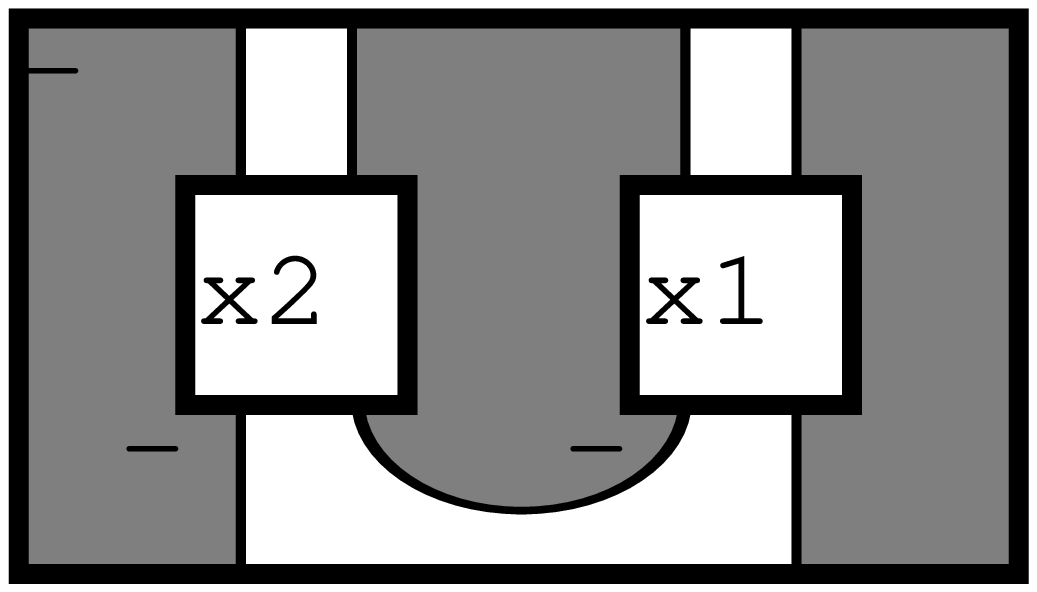}
}.
\]
\begin{lem}\label{duality}
$H_- \cong \left( H_+^{\t{op}} \right)^*$ as Kac algebras.
\end{lem}
\begin{pf}
Define a bilinear form $H_+ \times H_- \ni (p,a) \os{\langle \cdot , \cdot \rangle}{\longmapsto} \langle p,a \rangle := \delta^{-1} P_{
\psfrag{p}{$p$}
\psfrag{a}{$a$}
\psfrag{+}{$+$}
\psfrag{-}{$-$}
\includegraphics[scale=0.15]{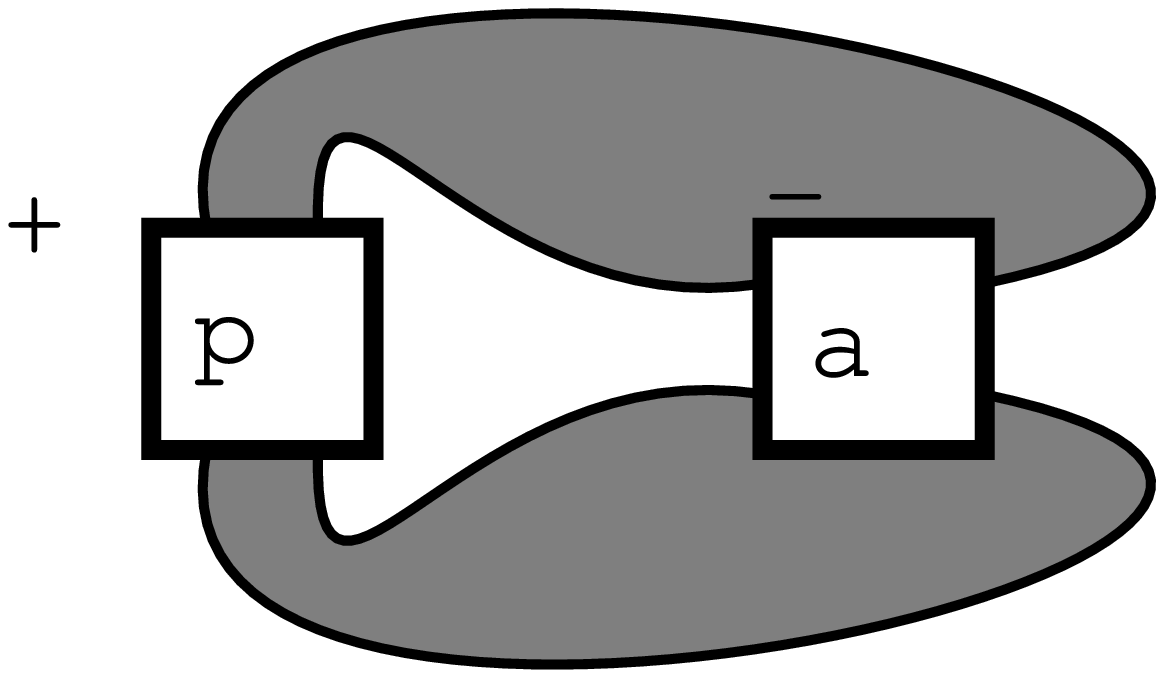}
}\in P_{+ 0} \cong \C$.
Non-degeneracy of the actions of the trace tangles, implies that $\langle \cdot , \cdot \rangle$ is non-degenerate. From the definition of the structural maps and the above formulae, it is easy to verify
\[
\langle p , a b \rangle = \langle p_{(1)} , a \rangle \langle p_{(2)} , b \rangle \;\;\;\;\; \t{;} \;\;\;\;\; \langle q p , a \rangle = \langle p , a_{(1)} \rangle \langle q , a_{(2)} \rangle \;\;\;\;\; \t{;} \;\;\;\;\; \langle p , a^* \rangle = \ol{\langle S_+ (p^*) , a \rangle }
\]
where $p , q \in H_+$ and $a,b \in H_-$. 
\end{pf}

\vspace{2mm}
\noindent Next, we recall the definition of the quantum double from \cite{Kas}. Let $H$ be a finite dimensional Hopf algebra.
The {quantum double of $H$} is the Hopf algebra $\left( H^{\t{op}} \right)^* \bowtie H$ (also denoted by $DH$) which is $\left( H^{\t{op}} \right)^* \otimes H$ as a vector space with structural maps given by:
\begin{itemize}
\item \emph{Multiplication:} $\left(f_{1}\bowtie h_{1}\right)\left(f_{2}\bowtie h_{2}\right)=f_{1}\left[f_{2}\left(S^{-1}\left(\left(h_{1}\right)_{(3)}\right)\cdot\left(h_{1}\right)_{(1)}\right)\right]\bt\left(h_{1}\right)_{(2)}h_{2},$
\item \emph{Unit:} $\chi_{H}\bowtie1$ ($\chi_H$ is the counit of $H$),
\item \emph{Comultiplication:} $\Delta\left(f\bowtie h\right)=f_{(1)}\bowtie h_{(1)}\otimes f_{(2)}\bowtie h_{(2)}$,
\item \emph{Counit:} $\chi\left(f\bowtie h\right)=f(1)\chi_{H}(h)$,
\item \emph{Antipode: $S\left(f\bowtie h\right)=f\left(h_{(3)}S^{-1}(\cdot)S^{-1}\left(h_{(1)}\right)\right)\bowtie S\left(h_{(2)}\right).$}
\end{itemize}
Moreover, if $H$ is a Hopf $*$-algebra, then $DH$ also has a $*$-structure given by 
\[
\left(f\bt h\right)^{*}=\ol{f} \left(h_{(3)}\left[S\circ*\left(\cdot\right)\right]S\left(h_{(1)}\right)\right)\bt h_{(2)}^{*} = f^* \left( S^{-1} (h^*_{(3)}) \cdot h^*_{(1)} \right) \bt h^*_{(2)}.
\]
\vspace{2mm}
Getting back to our context, by Lemma \ref{duality}, $DH_+$ can be considered as $H_- \bt H_+$.
\begin{rem}\label{DH=H-H+}
Using the duality defined in the proof of Lemma \ref{duality}, the $*$-algebra structure of $DH_+$ can be expressed as:
\begin{itemize}
\item Multiplication: $\left(a \bt p\right) \left(b \bt q\right) = P_{
\psfrag{-}{$-$}
\psfrag{+}{$+$}
\psfrag{a}{$a$}
\psfrag{b}{$b$}
\psfrag{p3}{$p_{(3)}$}
\psfrag{p1}{$p_{(1)}$}
\includegraphics[scale=0.2]{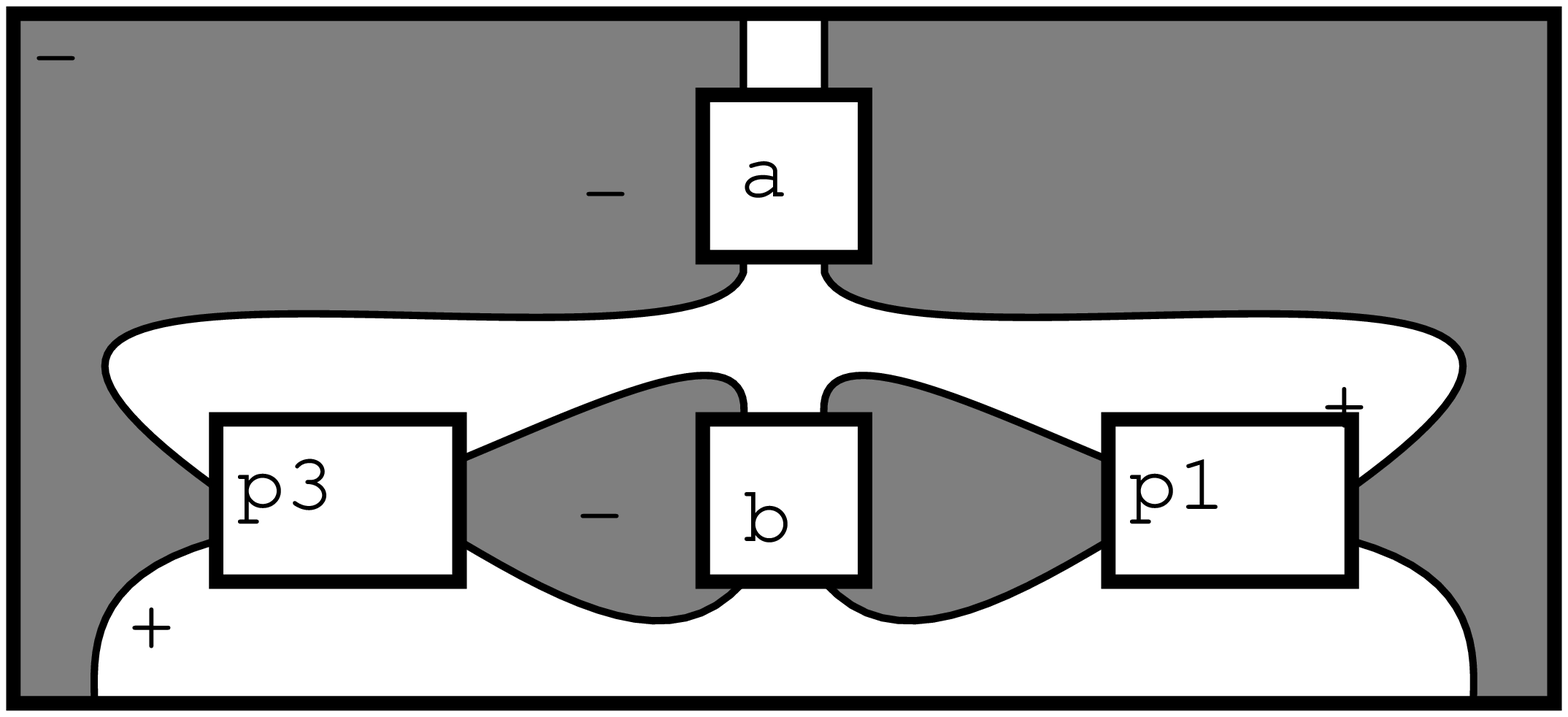}
} \bt p_{(2)} q$,
\item Unit: $1\bt 1$,
\item $*$-structure: $\left( a \bt p \right)^* = P_{
\psfrag{-}{$-$}
\psfrag{+}{$+$}
\psfrag{a}{$a^*$}
\psfrag{p3}{$p^*_{(3)}$}
\psfrag{p1}{$p^*_{(1)}$}
\includegraphics[scale=0.2]{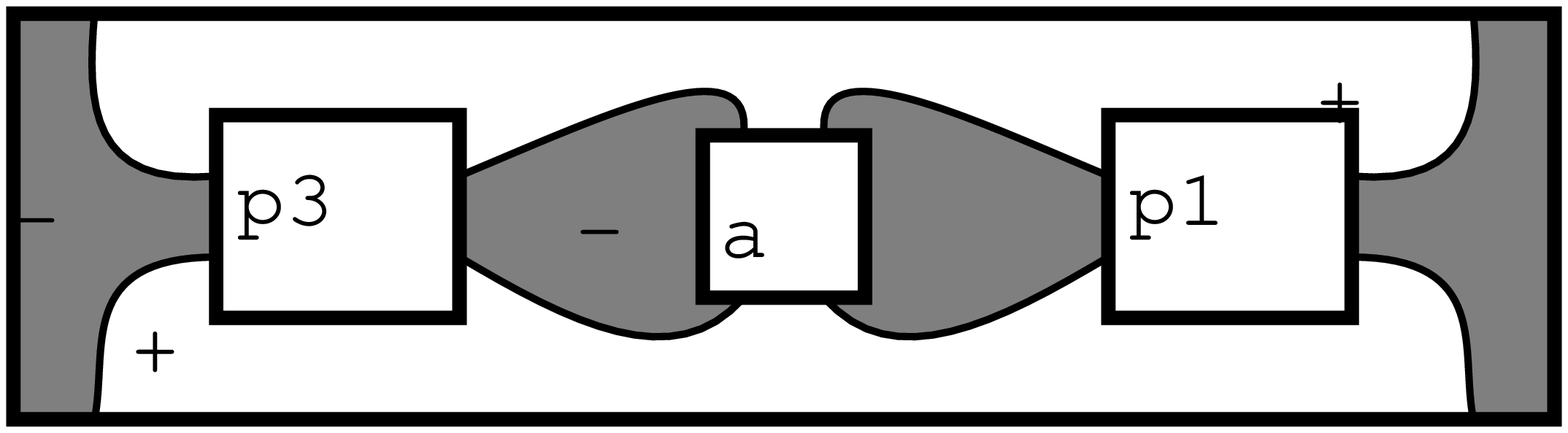}
} \bt p^*_{(2)}$.
\end{itemize}
\end{rem}
In order to establish a link between the quantum double of $H_+$ and the affine morphisms, we consider the tangles
$T^k_{l,m} :=
\psfrag{2l}{$2l$}
\psfrag{2m}{$2m$}
\psfrag{-}{$-$}
\psfrag{+}{$+$}
\psfrag{D1}{$D_{1}$}
\psfrag{D2}{$D_{2}$}
\psfrag{Dk}{$D_{k}$}
\psfrag{Dk1}{$D_{k+1}$}
\psfrag{Dk2}{$D_{k+2}$}
\psfrag{Dk3}{$D_{k+3}$}
\includegraphics[scale=0.2]{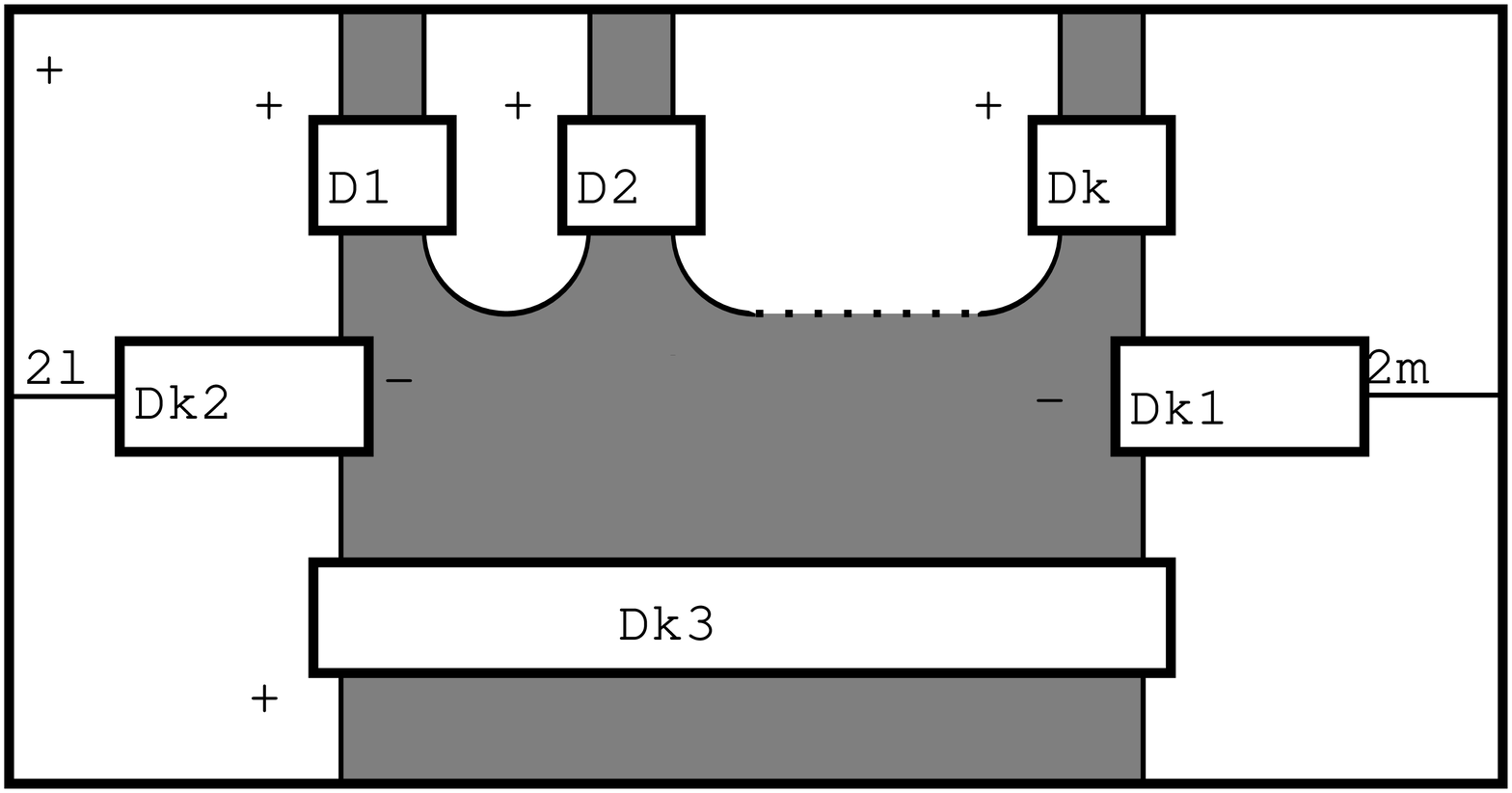}
$
and
$U:=
\psfrag{-}{$-$}
\psfrag{+}{$+$}
\psfrag{D1}{$D_{1}$}
\psfrag{D2}{$D_{2}$}
\includegraphics[scale=0.2]{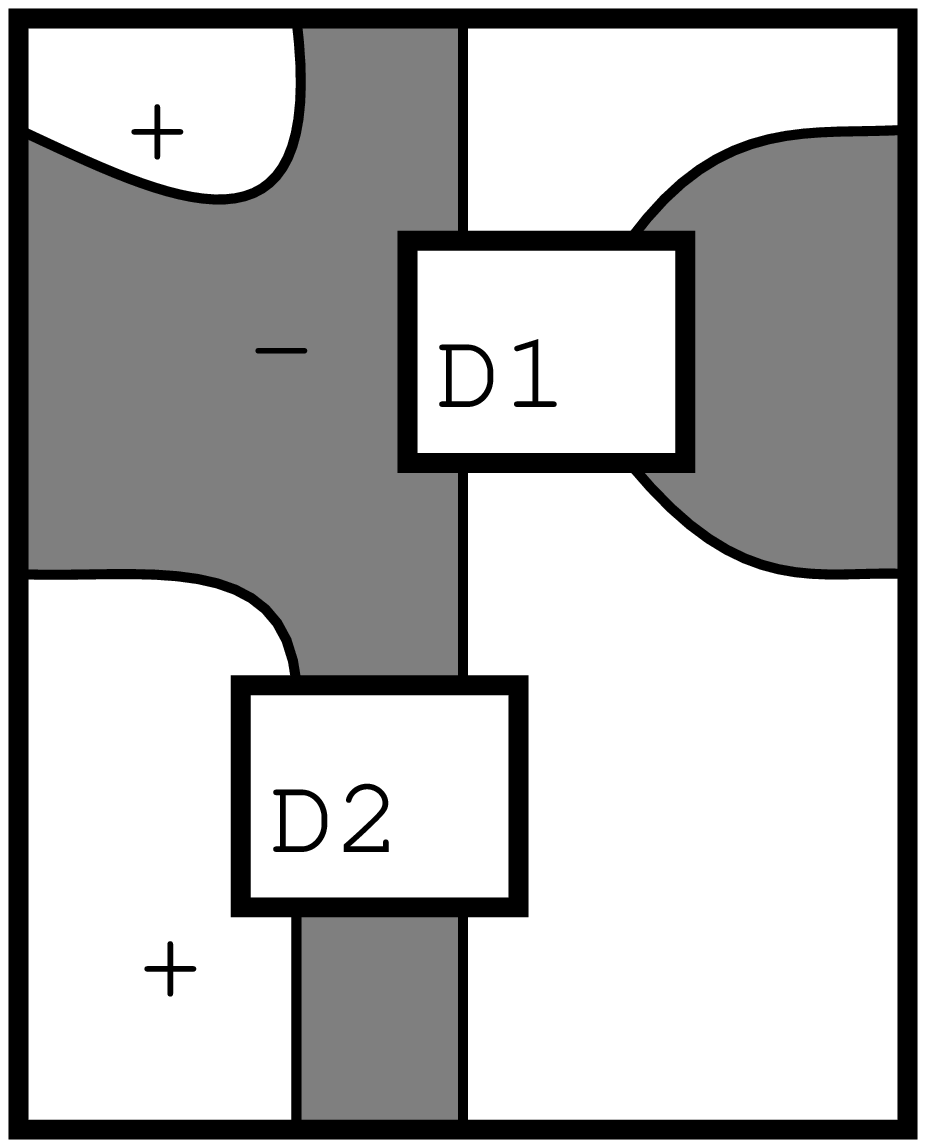}
$
where $k,l,m \geq 1$. When $l$ (resp., $m$) is zero, then the $T^k_{l,m}$ denotes the tangle obtained by composing the above tangle with $1_{-1}$ (defined in Figure \ref{tangles}) over the internal disc $D_{k+2}$ (resp., $D_{k+1}$).
Note that $P_{T^1_{1,1}} (1_{H_+} , \cdot, P_{E_{-1}} ,  \cdot )  = P_U$.
\begin{prop}\label{DHAP11}
The map $DH_+ \ni (a \bt p) \os{\Gamma}{\mapsto} \psi^2_{+1,+1} \left( P_{U} (a , p) \right) \in AP_{+1,+1}$ is a surjective $*$-algebra homomorphism.
\end{prop}
\begin{pf}
Using the structural maps of $H_\pm$ and $DH_+$, and affine isotopy, it is completely routine to check that $\Gamma$ preserves multiplication and $*$.

\noindent {\em Surjectivity of $\Gamma$:}
Consider an element $\psi^{2l}_{+1,+1} (x) \in AP_{+1,+1}$ for $x \in P_{+2(l+1)}$.
Now, Remark \ref{irrdep2gen} implies that $\t{Range} \, P_{T^1_{l,l}} (1_{H_+} , \cdot , \cdot , \cdot) = P_{+2(l+1)}$; so, without loss of generality, we can assume  $x = P_{T^{1}_{l,l}} (1_{H_+} , a, b , p)$ for $p \in H_+$ and $a , b \in P_{-(l+1)}$. 
Applying affine isotopy, we can write $\psi^{2l}_{+1,+1} (x) = \psi^2_{+1,+1} \left( P_U \left( (a \odot b) , p \right) \right)$ where $P_{-(l+1)} \times P_{-(l+1)} \ni (a,b) \os{\odot}{\longmapsto} P_{
\psfrag{-}{$-$}
\psfrag{2l}{$2l$}
\psfrag{a}{$a$}
\psfrag{b}{$b$}
\includegraphics[scale=0.15]{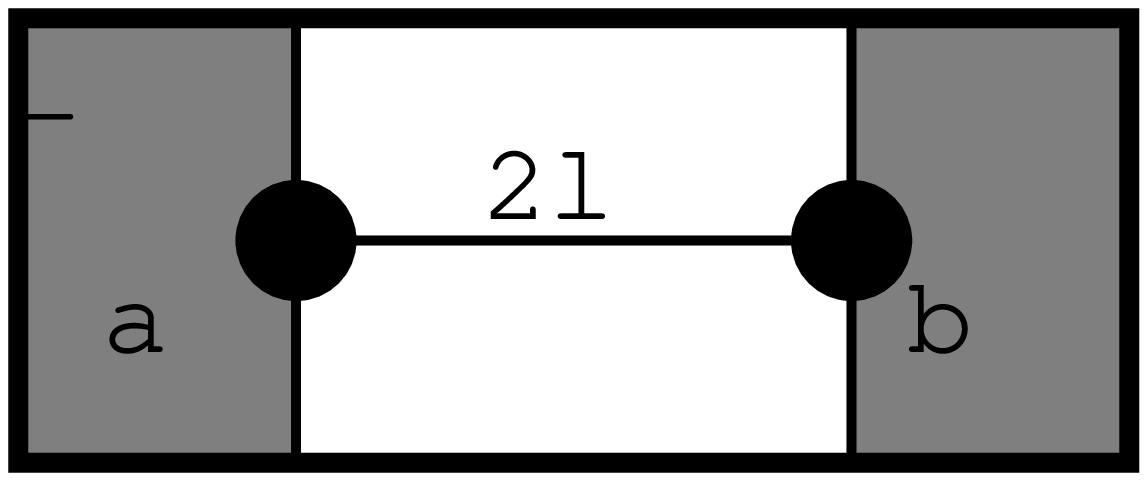}
} \in H_-$.
\end{pf}

\vspace{2mm}
Next, we proceed towards proving injectivity of $\Gamma$.
Set $V := \t{Range} \, P_U$ which (by Remark \ref{irrdep2gen}) is isomorphic to $H_- \otimes H_+$. Proposition \ref{DHAP11} implies $\t{dim}_\C \left( AP_{+1,+1} \right) \leq \t{dim}_\C (H_+) \t{dim}_\C (H_-)$.
So, it is enough to construct a surjective linear map from $AP_{+1,+1}$ to $V$.

For all $l\in \N$, consider the maps
\begin{align*}
P_{-(l+1)} \otimes P_{-(l+1)} \otimes H_+ \ni (a \otimes b \otimes p) & \os{\sigma_l}{\longmapsto} P_U \left( (a \odot b) , p \right) \in V \t{, and}\\
P_{-(l+1)} \otimes P_{-(l+1)} \otimes H_+ \ni (a \otimes b \otimes p) & \os{\tau_l}{\longmapsto} P_{T^1_{l,l}} (1_{H_+} , a , b , p) \in  P_{+2(l+1)}.
\end{align*}
By Remark \ref{irrdep2gen}, $\tau_l$ is an isomorphism and $\sigma_l$ is surjective.
Define the linear maps 
\[
\mcal P_{+2(l+1)} (P) \ni X \os{\gamma_l}{\longmapsto} \sigma_l \left( \tau^{-1}_l (P_X) \right) \in V \t{ for }l \geq 1 \;\;\;\;\;\;\; \t{ and } \;\;\;\;\;\;\; \mcal P_{+2} (P) \ni X \os{\gamma_0}{\longmapsto} P_U (1_{H_-} , P_X) \in V.
\]
We construct a map $\mcal{AT}_{+1,+1} (P) \ni A \os{\tilde{\gamma}}{\longmapsto} \gamma_l (T) \in V$ where $A= \Psi^{2l}_{+1,+1} (T)$ for some $T \in \mcal T_{+2(l+1)} (P)$, $l \geq 0$.
Then, the obvious question is whether $\tilde{\gamma}$ is well-defined.
If so, then we extend it linearly to $\tilde{\gamma}: \mcal A_{+1,+1} (P) \ra V$ which also becomes surjective and satisfies $\mcal W_{+1,+1} \subset \t{ker} \, \tilde{\gamma}$.
Thus, $\tilde{\gamma}$ factors through the quotient $AP_{+1,+1}$ and thereby, $\Gamma$ becomes injective.

\begin{figure}[h!]
\psfrag{+}{$+$}
\psfrag{~}{$\sim$}
\psfrag{2}{$2$}
\psfrag{2k}{$2k$}
\psfrag{2l}{$2l$}
\psfrag{S}{$S$}
\psfrag{T}{$T$}
\includegraphics[scale=0.2]{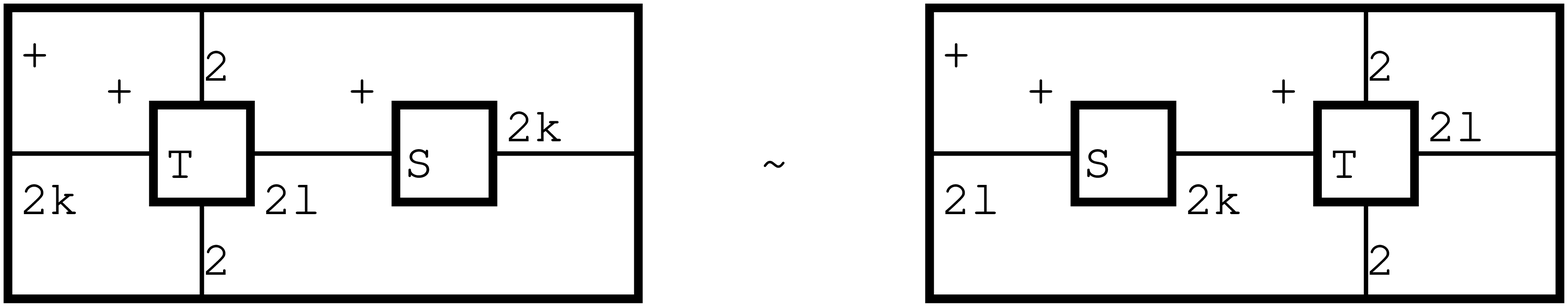}
\caption{$S \in \mcal T_{+(k+l)} (P)$, $T \in \mcal T_{+(k+l+2)} (P)$, $k,l \in \N_0$}\label{tangeqap11}
\end{figure}

\vspace{2mm}
\noindent{\em Well-definedness of $\tilde{\gamma}$:}
We will follow the treatment as in Section \ref{eqtang}.
Set $\mcal T := \us{l \in \N_0}{\sqcup} \mcal T_{+2(l+1)} (P)$.
We define an equivalence relation $\sim$ on $\mcal T$ in Figure \ref{tangeqap11}.
Analogous to Lemma \ref{psi-equivalence}, we have the following useful straight forward adaptation of \cite[Proposition 2.8]{Gho06} to the setting of morphisms in the affine category over a planar algebra.
\begin{lem}
For $X \in \mcal T_{+2(k+1)} (P)$, $Y \in \mcal T_{+2(l+1)} (P)$, $k,l \in \N_0$, we have:

(i) $\Psi^{2k}_{+1,+1} (X) = \Psi^{2l}_{+1,+1} (Y)$ if and only if $X \sim Y$, and

(ii) $X \sim Y$ implies $\gamma_k (X) = \gamma_l (Y)$.
\end{lem}
\begin{pf}
(i) The `if' part can easily be seen using affine isotopy.
The `only if' part can be proved by following the arguments used in the proof of the `only if' part in Lemma \ref{psi-equivalence}.

(ii) Suppose $X$ and $Y$ are the tangles on the left and right sides of $\sim$ in Figure \ref{tangeqap11} respectively.
Remark \ref{irrdep2gen} implies $\t{Range} \, P_{T^1_{k,l}} (1_{H_+}, \cdot, \cdot, \cdot) = P_{+(k+l+2)}$; so, there exists $\{a_i\}_i \subset P_{-(l+1)}$, $\{b_i\}_i \subset P_{-(k+1)}$ and $\{p_i\}_i \subset H_+$ such that $P_T = \us{i}{\sum} P_{T^1_{k,l}} (1_{H_+}, a_i, b_i, p_i)$.
Now, if $k , l > 0$, then
\[
\gamma_k (X) = \us{i}{\sum} P_U \left(P_{
\psfrag{-}{$-$}
\psfrag{2k}{$2k$}
\psfrag{2l}{$2l$}
\psfrag{a}{$a_i$}
\psfrag{b}{$b_i$}
\psfrag{PS}{$P_S$}
\includegraphics[scale=0.15]{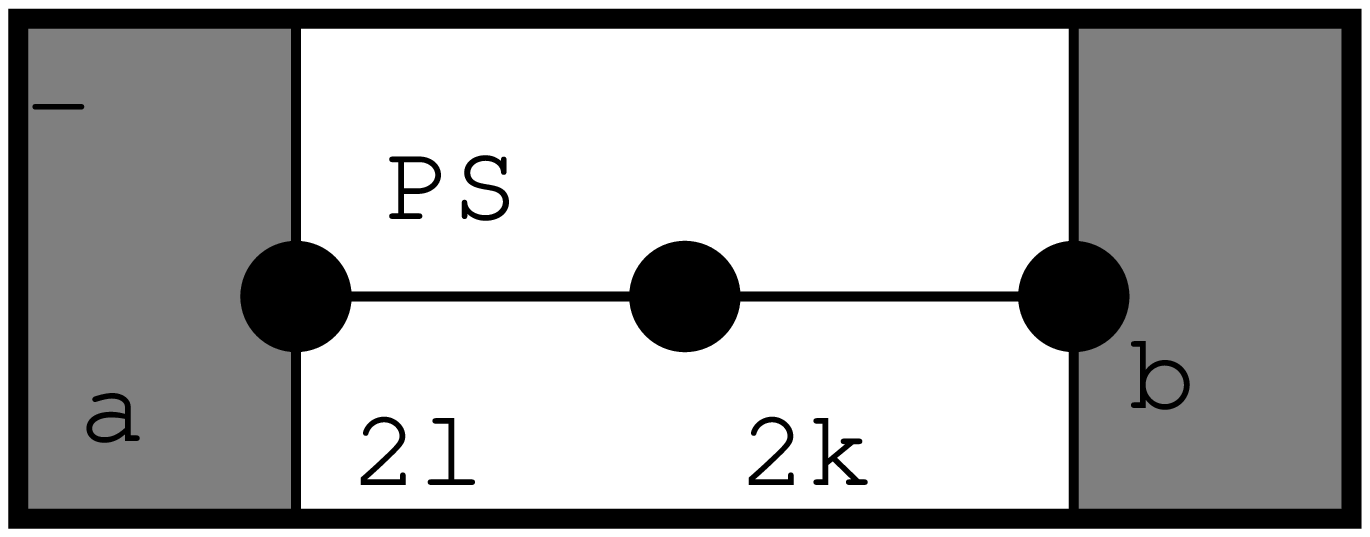}
} , p_i \right) = \gamma_l (Y).
\]

The case when $k = 0 = l$, the equation holds trivially. 

Suppose $k=0 < l$.
Again, by Remark \ref{irrdep2gen}, there exists $\{a_i\}_i \subset P_{-(l+1)}$ and $\{p_i\}_i \subset H_+$ such that $P_T = \us{i}{\sum} P_{T^1_{0,l}} (1_{H_+}, a_i, p_i)$.
Note that $P_Y = \us{i}{\sum} P_{T^1_{l,l}} (1_{H_+}, a_i , P_{LI_{+l}} (P_S) , p_i)$.
Thus, $\gamma_l (Y) = \us{i}{\sum} P_U (a_i \odot P_{LI_{+l}} (P_S) , p_i)$.
Since $P$ is irreducible, there exists $c_i \in \C$ such that $a_i \odot P_{LI_{+l}} (P_S) = c_i 1_{H_-}$ which also implies $P_X = \us{i}{\sum} c_i \, p_i$. Hence, $\gamma_0 (X) = \us{i}{\sum} c_i P_U (1_{H_-} , p_i) = \gamma_l (Y)$.
Similar arguments yeild the case $k > 0 = l$.
\end{pf}

Hence, we have proved the following proposition.
\begin{prop}\label{giso}
$\Gamma$, as in Proposition \ref{DHAP11}, is an isomorphism.
\end{prop}
\subsection{The affine modules of $P$}
Let $t_\vlon : H_\vlon \ra \C$ denote the normalized action of the trace tangle on $H_\vlon$, that is, $t_\vlon = \delta^{-2} P_{TR^r_{\vlon 2}}$.
Consider the linear functional $DH_+ \ni a \bt p \os{t}{\longmapsto} t_- (a) t_+ (p) \in \C$.
From Remark \ref{DH=H-H+} and the structural maps in the beginning of Section \ref{affmorlevel1}, it easily follows $t \left( (a\bt p)^* (a\bt p) \right) = t_- (a^* a) t_+(p^* p)$; thus, $t_- \bt t_+$ is positive definite and $DH_+$ becomes a finite dimensional $C^*$-algebra. Set $\tilde{t} := t \circ \Gamma^{-1} : AP_{+1,+1} \ra \C$.
\begin{thm}\label{conj}
If $N \subset M$ is an irreducible subfactor with depth two and planar algebra $P$, then the category of Hilbert affine $P$-modules is equivalent to the center of the category of $N$-$N$-bimodules generated by $\lrsuf{N}{L^2 (M)}{}_M$ as additive categories.
\end{thm}
\begin{pf}
From \cite{Sz}, one can deduce that the category of $N$-$N$-bimodules (appearing in the standard invariant) is contravariantly equivalent to the representation category of the Kac algebra $H_+$; thus, its center then becomes contravariantly equivalent to the representation category of $DH_+$ (see \cite{Kas}, Theorem XIII.5.1).
So, using Remark \ref{eaff}, it is enough to establish a one-to-one correspondence between the isomorphism classes of irreducible Hilbert $+$-affine $P$-modules and that of irreducible $DH_+$-modules.
The key step towards this will be given by the following construction of an Hilbert $+$-affine $P$-module  generated by $AP_{+1,+1}$.

Set $V_k := AP_{+1,+k}$ for all $k \geq 1$ and $V_{\vlon 0} := AP_{+1 , {\vlon 0}}$.
Note that by \cite[Proof of Theorem 6.11]{Gho}, $V_k$'s are all finite dimensional.
We define a sesquilinear form $\langle v , w \rangle := \tilde{t}(v^* \circ w)$ for all $v , w \in V_k$, $k\in \{ \pm 0\}\cup \N$.

\vspace{2mm}
\noindent{\em Positivity of $\langle \cdot , \cdot \rangle$:}
The case $k = 1$ is already covered by Proposition \ref{giso}.

\noindent{\bf Case 1:} Suppose $k > 1$. Note that $V_k = \us{l \in \N}{\bigcup} \psi^{2l}_{+1,+k} (P_{+(2l+k+1)})$.
Now, Remark \ref{irrdep2gen} implies $\t{Range} \, P_{T^k_{l,l}} (1_{H_+} , \cdot, \ldots ,  \cdot) = P_{+(2l+k+1)}$.
Applying affine isotopy, we get
\begin{align*}
& \; \psi^{2l}_{+1,+k} \left( P_{T^k_{l,l}} (1_{H_+} , x_2 , \ldots , x_k , a , b , p ) \right)
= \psi^2_{+1,+k} \left( P_{T^k_{1,1}} (1_{H_+} , x_2 , \ldots , x_k , a \odot b , P_{E_{-1}} , p ) \right)\\
= & \; \psi^0_{+1,+k} \left( P_{T^k_{0,0}} (1_{H_+} , x_2 , \ldots , x_k , 1_{H_+}) \right)
\circ \psi^2_{+1,+1} \left( P_U (a \odot b , p ) \right)
\end{align*}
(which is independent of $l$) for all $x_2 , \dots x_k \in H_+$, $a,b \in P_{-(l+1)}$ and $p \in H_+$.
Thus, the linear map defined by 
\[
\left[ (H_+)^{\otimes(k-1)} \otimes AP_{+1,+1} \right] \ni x_2 \otimes \cdots \otimes x_k \otimes w \os{\zeta}{\mapsto} \psi^2_{+1,+k} \left( P_{T^k_{0,0}} (1_{H_+} , x_2 , \ldots , x_k , 1_{H_+}) \right) \circ w \in V_k
\]
is surjective.
Since $P$ is irreducible, we have
\begin{align*}
&\left[ \psi^2_{+1,+k} \left( P_{T^k_{0,0}} (1_{H_+} , x_2 , \ldots , x_k , 1_{H_+}) \right) \right]^* \circ \psi^2_{+1,+k} \left( P_{T^k_{0,0}} (1_{H_+} , x_2 , \ldots , x_k , 1_{H_+}) \right)\\
= \; & \delta^k \left[ \us{2\leq n \leq k}{\prod} t_+ (x^*_n x_n) \right] 1_{AP_{+1,+1}}
\end{align*}
which implies
$\langle \zeta(x_2 \otimes \cdots \otimes x_k \otimes w) , \zeta(y_2 \otimes \cdots \otimes y_k \otimes v) \rangle = \delta^k \left[ \us{2\leq n \leq k}{\prod} t_+ (x^*_n y_n) \right] t(w^* \circ v)$ for all $x_2, \ldots , x_k, y_2, \ldots, y_k \in H_+$ and $w,v \in AP_{+1,+1}$.
Hence $\langle \cdot , \cdot \rangle$ is positive definite on $V_k$.

\noindent{\bf Case 2:} Suppose $k = \vlon 0$.
Consider the affine morphism $c_\vlon \in AP_{\vlon 0,+1}$ given by the affine tangle with a single string attached to the two marked points on the boundary of the external rectangle.
Now, since  $c_\vlon \circ v \neq 0$ for all $0 \neq v \in V_{\vlon 0}$, we have $\langle v , v \rangle_{V_{\vlon 0}} = \delta^{-1} \left\langle c_\vlon \circ v , c_\vlon \circ v \right\rangle_{V_1} > 0$.

\vspace{2mm}
Hence, $V$ is a Hilbert affine $P$-module.
Now, $V_1$ is the regular $AP_{+1,+1}$-module; so, it contains every irreducible $AP_{+1,+1}$-module as a submodule.
By Remark \ref{irrmod}, the affine submodule $\tilde{W}$ of $V$ generated by each of these irreducible $AP_{+1,+1}$-submodule $W$ of $V_1$, will be irreducible; moreover, $W_1 \cong W_2$ if and only if $\tilde{W}_1 \cong \tilde{W}_2$.
On the other hand, if we start with an irreducible Hilbert affine $P$-module $U$, then $U_1$ is nonzero (since weight of every affine $P$-modules cannot exceed  $1$ by \cite[Theorem 6.10]{Gho}) and is an irreducible $AP_{+1,+1}$-module (see Remark \ref{irrmod}).
So, there exists a submodule $W$ of $V_1$, which is isomorphic to $U_1$.
Using Remark \ref{affmodmor}, we may conclude $U \cong \tilde{W}$.

Hence, we have established a one-to-one correspondence between the isomorphism class of irreducible Hilbert affine $P$-modules and that of irreducible $AP_{+1,+1}$-modules.
This ends the proof.\end{pf}

Note that Theorem \ref{conj} confirms Jones-Walker conjecture (stated in the introduction) for the case of irreducible depth two subfactors.
\vspace{4mm}
\begin{center}
{\bf Some questions.}
\end{center}

In Section \ref{regmod}, we provided an explicit way of constructing a large class of Hilbert affine $P$-modules generated by their zero spaces for any subfactor planar algebra $P$.
The natural question to ask is whether all $*$-affine $P$-modules (not necessarily bounded) generated by their zero spaces, arise in this way for infinite depth $P$'s.
It will also be interesting to analyze the affine $P$-modules with weight greater than zero.

In Section \ref{JWconjirrdep2}, we used irreducibility of $P$ quite crucially in affirming the Jones-Walker conjecture. The next obvious thing to check will be whether we can make this work in the absence of irreducibility, that is, the `weak Hopf algebra' case.
An important drawback of the category of the Hilbert affine $P$-modules, is the lack of a monoidal structure, let alone braiding; note that the equivalence established in Theorem \ref{conj}, is an equivalence of additive categories.
One would guess some kind of comultiplication structure on the affine category, might yeild an appropriate monoidal structure on the category of affine $P$-modules.

We will address and answer some of these questions in a forthcoming article.
\bibliographystyle{alpha}

\end{document}